\documentclass[a4paper,10pt]{article}

\usepackage[OT2,T1]{fontenc}
\usepackage[french, english]{babel}

\usepackage{amsfonts}
\usepackage{amsmath,amssymb}
\usepackage{pstricks}
\usepackage{amsthm}
\usepackage[all]{xy}

\title{Th\'eor\`emes de dualit\'e pour les complexes de tores}
\author{Cyril Demarche}
\date{\today}

\newtheorem{theo}{Th\'eor\`eme}[section]
\newtheorem*{thm}{Th\'eor\`eme}

\newtheorem{defi}[theo]{D\'efinition}
\newtheorem{prop}[theo]{Proposition}

\newtheorem{cor}[theo]{Corollaire}
\newtheorem{lem}[theo]{Lemme}
\newtheorem{rem}[theo]{Remarque}
\newenvironment{dem}{\noindent
  \textit{{D\'emonstration}} : }
  {\hfill \qedsymbol\newline}

\DeclareSymbolFont{rsfs}{U}{rsfs}{m}{n}
\DeclareSymbolFontAlphabet{\mathcal}{rsfs}

\newcommand{\ensemble}[1]{\ensuremath{\mathbf #1} \xspace}
  \renewcommand{\H}{\ensemble H}
  \newcommand{\N}{\ensemble N}

  \newcommand{\Z}{\ensemble Z}
  \newcommand{\D}{\ensemble D}
  \newcommand{\Q}{\ensemble Q}
  \newcommand{\R}{\ensemble R}
  \renewcommand{\P}{\ensemble P}
  \newcommand{\C}{\ensemble C}
  
  \renewcommand{\L}{\ensemble L}
  \newcommand{\F}{\ensemble F}
  \renewcommand{\L}{\ensemble L}
  \newcommand{\G}{\ensemble G}
  \newcommand{\A}{\ensemble A}
  \renewcommand{\C}{\ensemble C}
  
  \newcommand\cyr[1]{{\fontencoding{OT2}\fontfamily{wncyr}\selectfont #1}}

\newenvironment{changemargin}[2]{
 \begin{list}{}{
  \setlength{\topsep}{0pt}
  \setlength{\leftmargin}{#1}
  \setlength{\rightmargin}{#2}
  \setlength{\listparindent}{\parindent}
  \setlength{\itemindent}{\parindent}
  \setlength{\parsep}{\parskip}
 }
\item[]}{\end{list}}

\newcommand{\fct}[4]{\begin{displaymath}
\begin{array}{ccc}
#1 & \longrightarrow & #2 \\
#3 & \longmapsto & #4
\end{array}
\end{displaymath} }

\begin{document}	

\selectlanguage{french}

\maketitle

\begin{abstract}
On consid\`ere un complexe de tores de longueur $2$ d\'efini sur un corps de nombres $k$. On \'etablit des r\'esultats de dualit\'e locale et globale pour l'hypercohomologie (\'etale ou galoisienne) de ce complexe. On obtient notamment une suite de Poitou-Tate pour de tels complexes, g\'en\'eralisant les suites de Poitou-Tate pour les modules galoisiens finis ou les tores. En particulier, on d\'emontre l'existence d'une telle suite pour les $k$-groupes de type multiplicatifs. Les r\'esultats g\'en\'eraux obtenus ici pour les complexes de tores ont par exemple des applications dans des r\'esultats r\'ecents concernant le d\'efaut d'approximation forte dans les groupes lin\'eaires connexes et des th\'eor\`emes de dualit\'e sur la cohomologie galoisienne (non-ab\'elienne) de tels groupes.
\end{abstract}

\selectlanguage{english}

\begin{abstract}
We consider a complex of tori of length $2$ defined over a number field $k$. We establish here some local and global duality theorems for the (\'etale or Galois) hypercohomology of such a complex. We prove the existence of a Poitou-Tate exact sequence for such a complex, which generalizes the Poitou-Tate exact sequences for finite Galois modules and tori. In particular, we obtain a Poitou-Tate exact sequence for $k$-groups of multiplicative type. The general results proven here lie at the root of recent results about the defect of strong approximation in connected linear algebraic groups and about some arithmetic duality theorems for the (non-abelian) Galois cohomology of such groups.
\end{abstract}

\selectlanguage{french}

\section{Introduction}

Les th\'eor\`emes de dualit\'e pour la cohomologie galoisienne des groupes alg\'ebriques commutatifs sur les corps locaux et globaux constituent parmi les plus importants r\'esultats en arithm\'etique.

Un permier r\'esultat bien connu est le th\'eor\`eme de dualit\'e locale pour les modules galoisiens finis (voir par exemple \cite{Mil}, corollaire I.2.3).

On dispose aussi des r\'esultats de Tate-Nakayama et Tate pour les tores et les vari\'et\'es ab\'eliennes sur des corps locaux, r\'esultats g\'en\'eralis\'es par Harari et Szamuely (pour des vari\'et\'es semi-ab\'eliennes et plus g\'en\'eralement les 1-motifs sur un corps local) : si $K$ est un corps local et $M$ un 1-motif sur $K$, et $M^*$ le 1-motif dual, alors on dispose d'accouplements naturels, d\'efinis par le cup-produit
$$\H^i(K,M) \times \H^{1-i}(K,M^*) \rightarrow \Q / \Z$$
($-1 \leq i \leq 2$) qui induisent des dualit\'es parfaites entre ces deux groupes (apr\`es une compl\'etion profinie pour $i = -1$ ou $0$) : voir \cite{HSz}, th\'eor\`eme 2.3.

Concernant les th\'eor\`emes de dualit\'e globale, on dispose du th\'eor\`eme de Tate pour les modules galoisiens finis, \`a savoir : si $F$ est un module galoisien fini sur un corps de nombres $k$, on dispose d'accouplements naturels non-d\'eg\'en\'er\'es entre groupes finis
$$\textup{\cyr{SH}}^i(F) \times \textup{\cyr{SH}}^{3-i}(F^*) \rightarrow \Q / \Z$$
o\`u $F^*$ est le dual de Cartier de $F$, et $i = 1,2$ (voir \cite{Mil}, th\'eor\`eme 4.10.(a)).

Un autre exemple de r\'esultat de dualit\'e globale est donn\'e par le th\'eor\`eme de Cassels-Tate pour les vari\'et\'es ab\'eliennes, ainsi que par le th\'eor\`eme analogue pour les tores, souvent attribu\'e \`a Kottwitz (voir par exemple l'appendice de \cite{KoS}, o\`u les preuves ne sont parfois pas compl\`etes). Ces deux r\'esultats ont \'et\'e g\'en\'eralis\'es par Harari et Szamuely, sous la forme suivante : si $k$ est un corps de nombres et $M$ est un 1-motif sur $k$, alors on dispose d'un accouplement canonique
$$\textup{\cyr{SH}}^i(M) \times \textup{\cyr{SH}}^{2-i}(M^*) \rightarrow \Q / \Z$$
($i = 0, 1$) qui est non-d\'eg\'en\'er\'e  modulo les sous-groupes divisibles. Les auteurs d\'emontrent d'ailleurs diff\'erentes versions de ce r\'esultat, en rempla\c cant parfois les groupes de Tate-Shafarevich par des variantes faisant intervenir des compl\'etions profinies : voir \cite{HSz}, th\'eor\`eme 4.8, corollaire 4.9, propositions 4.12 et 5.1.

Enfin, on peut rassembler les th\'eor\`emes de dualit\'e locale et globale dans une suite exacte dite de Poitou-Tate : voir \cite{HSz}, th\'eor\`eme 5.6 pour les 1-motifs, et \cite{Mil}, th\'eor\`eme 4.10 pour les modules galoisiens finis.

On citera enfin les r\'esultats de Gonzales-Aviles (voir \cite{Gon}) \`a propos des th\'eor\`emes de dualit\'e sur les corps locaux et globaux de caract\'eristique positive. 

L'objectif de ce texte est de d\'emontrer de nouveaux th\'eor\`emes de dualit\'e, locale et globale, pour l'hypercohomologie des complexes de tores de longueur $2$, qui g\'en\'eralisent notamment les r\'esultats pour les tores et les modules finis rappel\'es plus haut. On retrouve en particulier (par une m\'ethode diff\'erente) certains des r\'esultats de l'appendice de \cite{KoS}. Ces r\'esultats de dualit\'e ont notamment des applications pour le calcul de la cohomologie galoisienne des groupes lin\'eaires connexes : on citera par exemple les travaux de Borovoi (notamment \cite{BorAMS} et \cite{Bor3} section 4), ceux de Kottwitz et Shelstad (voir \cite{KoS}), et les r\'esultats de l'auteur dans \cite{Dem2}.

Rappelons quelques notations avant d'\'enoncer les r\'esultats de ce texte. On utilise les notations usuelles suivantes : si $A$ est un groupe topologique ab\'elien, on note $A^D$ le groupe des morphismes de groupes continus $A \rightarrow \Q / \Z$. On munit ce groupe $A^D$ de la topologie compacte-ouverte. On note $A^{\wedge}$ le compl\'et\'e de $A$ pour la topologie des sous-groupes ouverts d'indice fini, et $A_{\wedge} := \varprojlim_n A / {n A}$.

Si $k$ est un corps, on note $\Gamma_k$ son groupe de Galois absolu. Si $C$ est un complexe de modules galoisiens sur $k$, et $i \in \N$, on note $\H^i(k, C) := \H^i(\Gamma_k, C(\overline{k}))$ le $i$-\`eme groupe d'hypercohomologie galoisienne.
Si $k$ est un corps de nombres et $\Omega_k$ l'ensemble des places de $k$, $\P^i(k,C)$ d\'esigne le produit restreint des groupes $\H^i(\widehat{k}, C)$ par rapport aux groupes $\H^i(\widehat{\mathcal{O}}_v, \mathcal{C})$, o\`u $\widehat{\mathcal{O}}_v$ d\'esigne l'anneau des entiers de $\widehat{k}_v$. On d\'efinit aussi
$$\textup{\cyr{SH}}^i(k, C) = \textup{\cyr{SH}}^i(C) := \textup{Ker}\left( \H^i(k, C) \rightarrow \P^i(k,C) \right)$$
o\`u $\widehat{k}_v$ d\'esigne le compl\'et\'e de $k$ \`a la place $v$. 
On aura \'egalement besoin d'une version modifi\'ee de ces groupes, \`a savoir :
$$\textup{\cyr{SH}}^i_{\wedge}(C) := \textup{Ker} \left( \H^i(k, C)_{\wedge} \rightarrow \P^i(k, C)_{\wedge} \right)$$

On peut d\'esormais r\'esumer les r\'esultats principaux de ce texte.

\'Etant donn\'e un complexe de tores $C = [T_1 \xrightarrow{\rho} T_2]$ (en degr\'es $-1$ et $0$) sur un corps local $K$, les groupes $\H^i(K,C)$ sont munis d'une topologie naturelle, et l'on dispose d'un accouplement canonique entre les groupes d'hypercohomologie
$$\H^i(K,C) \times \H^{1-i}(K, \widehat{C}) \rightarrow \Q / \Z$$
qui induit une dualit\'e parfaite apr\`es certaines compl\'etions profinies (voir th\'eor\`eme \ref{dualite locale}), o\`u $\widehat{C}$ d\'esigne le complexe (concentr\'e en degr\'es $-1$ et $0$) de modules galoisiens $[\widehat{T_2} \xrightarrow{\widehat{\rho}} \widehat{T_1}]$, $\widehat{T_i}$ \'etant le module des caract\`eres du tore $T_i$.

On obtient \'egalement le r\'esultat de dualit\'e globale suivant (voir th\'eor\`emes \ref{theo SH 1}, \ref{theo SH 2}, \ref{dualite SH 0}, \ref{dualite SH 0 tm} et proposition \ref{prop 2}) :
\begin{thm}
Soit $k$ un corps de nombres et $C = [T_1 \xrightarrow{\rho} T_2]$ un complexe de $k$-tores. 
Alors il existe un accouplement canonique
$$\textup{\cyr{SH}}^1(C) \times \textup{\cyr{SH}}^1(\widehat{C}) \rightarrow \Q / \Z$$
qui est une dualit\'e parfaite de groupes finis.
Il existe un accouplement canonique, fonctoriel en $C$
$$\textup{\cyr{SH}}^2(C) \times \textup{\cyr{SH}}^0_{\wedge}(\widehat{C}) \rightarrow \Q / \Z$$
qui est une dualit\'e parfaite. 
\begin{itemize}
 	\item Si $\textup{Ker } \rho$ est fini, alors il existe des accouplements canoniques, fonctoriels en $C$
$$\textup{\cyr{SH}}^1(C) \times \textup{\cyr{SH}}^1(\widehat{C}) \rightarrow \Q / \Z$$
$$\textup{\cyr{SH}}^2(C) \times \textup{\cyr{SH}}^0(\widehat{C}) \rightarrow \Q / \Z$$
et
$$\textup{\cyr{SH}}^0_{\wedge}(C) \times \textup{\cyr{SH}}^2(\widehat{C}) \rightarrow \Q / \Z$$
qui sont des dualit\'es parfaites entre groupes finis.
	\item Si le morphisme $\rho$ est surjectif, il existe des accouplements canoniques, fonctoriels en $C$
$$\textup{\cyr{SH}}^1(C) \times \textup{\cyr{SH}}^1(\widehat{C}) \rightarrow \Q / \Z$$
et
$$\textup{\cyr{SH}}^0(C) \times \textup{\cyr{SH}}^2(\widehat{C}) \rightarrow \Q / \Z$$
qui sont des dualit\'es parfaites entre groupes finis.
\end{itemize}
\end{thm}

\begin{rem}
{\rm
L'hypoth\`ese de surjectivit\'e du morphisme $\rho$ implique que le complexe $C$ est quasi-isomorphe (dans la cat\'egorie d\'eriv\'ee associ\'ee \`a la cat\'egorie des complexes born\'es de modules galoisiens sur $k$) \`a $(\textup{Ker } \rho)[1]$, $\textup{Ker } \rho$ \'etant un $k$-groupe de type multiplicatif. R\'eciproquement, tout $k$-groupe de type multiplicatif est le noyau d'un morphisme surjectif de $k$-tores : par cons\'equent, le th\'eor\`eme pr\'ec\'edent contient un th\'eor\`eme de dualit\'e globale pour les groupes de type multiplicatif.
}
\end{rem}

On obtient \'egalement des suites de type Poitou-Tate pour certains complexes de tores, comme par exemple le r\'esultat suivant, o\`u les diff\'erents morphismes proviennent des th\'eor\`emes de dualit\'e locale et globale, et $(.)^{\wedge}$ d\'esigne la compl\'etion pour la topologie des sous-groupes ouverts d'indice fini :
\begin{thm}[Th\'eor\`eme \ref{theo PT}]
Soit $C = [T_1 \xrightarrow{\rho} T_2 ]$ un complexe de tores d\'efini
sur $k$, avec $\textup{Ker}(\rho)$ fini. On a alors une suite exacte de groupes topologiques,
fonctorielle en $C$ :
\begin{displaymath}
\xymatrix{
0 \ar[r] & \H^{-1}(k, C) \ar[r] & \P^{-1}(k, C) \ar[r] & \H^2(k, \widehat{C})^D \ar[d] \\
& \H^1(k, \widehat{C})^D \ar[d] & \P^0(k, C)^{\wedge} \ar[l] & \H^0(k, C)^{\wedge} \ar[l] \\
& \H^1(k,C) \ar[r] & \P^1(k, C) \ar[r] & \H^0(k, \widehat{C})^D \ar[d] \\
0 & \H^{-1}(k, \widehat{C})^D \ar[l] & \P^2(k, C) \ar[l] & \H^2(k, C) \ar[l]
}
\end{displaymath}
On dispose \'egalement, sous les m\^emes hypoth\`eses, de la suite exacte duale :
\begin{displaymath}
\xymatrix{
0 \ar[r] & \H^{-1}(k, \widehat{C})^{\wedge} \ar[r] & \P^{-1}(k, \widehat{C})^{\wedge} \ar[r] & \H^2(k, C)^D \ar[d] \\
& \H^1(k, C)^D \ar[d] & \P^0(k, \widehat{C}) \ar[l] & \H^0(k, \widehat{C}) \ar[l] \\
& \H^1(k,\widehat{C}) \ar[r] & \P^1(k, \widehat{C})_{\textup{tors}} \ar[r] & \left( \H^0(k, C)^D \right)_{\textup{tors}} \ar[d] \\
0 & \H^{-1}(k, C)^D \ar[l] & \P^2(k, \widehat{C}) \ar[l] & \H^2(k, \widehat{C}) \ar[l]
}
\end{displaymath}
\end{thm}

\begin{rem}
{\rm
Comme on le montre dans \cite{Dem2}, ce th\'eor\`eme est un outil important pour \'etudier le d\'efaut d'approximation forte dans les groupes lin\'eaires connexes, ainsi que pour obtenir une suite de Poitou-Tate non-ab\'elienne pour de tels groupes, \`a l'aide des applications d'ab\'elianisation d\'efinies notamment par Borovoi dans le chapitre 3 de \cite{BorAMS}.
}
\end{rem}

On obtient aussi une suite de Poitou-Tate pour les groupes de type multiplicatif, g\'en\'eralisant les suites usuelles de Poitou-Tate pour les modules galoisiens finis (voir \cite{Mil}, th\'eor\`eme I.4.10) et pour les tores (voir le th\'eor\`eme 5.6 de Harari et Szamuely dans \cite{HSz}). Cette suite n'\'etait pas connue pour des groupes de type multiplicatif g\'en\'eraux : voir notamment la remarque apr\`es le corollaire I.4.21 de \cite{Mil}, ainsi que la suite de Poitou-Tate partielle du th\'eor\`eme 8.6.14 de \cite{NSW}.

\begin{thm}[Th\'eor\`eme \ref{theo PT tm}]
Soit $M$ un $k$-groupe de type multiplicatif. On a alors une suite
exacte de groupes topologiques, fonctorielle en $M$ :
\begin{displaymath}
\xymatrix{
0 \ar[r] & H^{0}(k, M)^{\wedge} \ar[r] & P^{0}(k, M)^{\wedge} \ar[r] & H^2(k, \widehat{M})^D \ar[d] & \\
& H^1(k, \widehat{M})^D \ar[d] & P^1(k, M) \ar[l] & H^1(k, M) \ar[l] & \\
& H^2(k,M) \ar[r] & P^2(k, M) \ar[r] & H^0(k, \widehat{M})^D \ar[r] & 0 
}
\end{displaymath}
On dispose \'egalement de la suite exacte duale :
\begin{displaymath}
\xymatrix{
0 \ar[r] & H^{0}(k, \widehat{M})^{\wedge} \ar[r] & P^{0}(k, \widehat{M})^{\wedge} \ar[r] & H^2(k, M)^D \ar[d] & \\
& H^1(k, M)^D \ar[d] & P^1(k, \widehat{M}) \ar[l] & H^1(k, \widehat{M}) \ar[l] & \\
& H^2(k,\widehat{M}) \ar[r] & P^2(k, \widehat{M})_{\textup{tors}} \ar[r] & \left( H^0(k, M)^D \right)_{\textup{tors}} \ar[r] & 0
}
\end{displaymath}
\end{thm}

Le plan du texte est le suivant : on montre d'abord les th\'eor\`emes de dualit\'e locale \`a la section \ref{section locale}, puis on \'etablit des r\'esutats de dualit\'e en cohomologie \'etale \`a l'aide du th\'eor\`eme d'Artin-Verdier (section \ref{section etale}). Ensuite, la section \ref{section Galois} est consacr\'ee \`a la d\'emonstration des th\'eor\`emes de dualit\'e globale. Enfin, on obtient les suites exactes de type Poitou-Tate \`a la section \ref{section PT}, et on fait le lien avec une suite en hypercohomologie obtenue par Borovoi.

\paragraph{Remerciements}
Je remercie tr\`es chaleureusement David Harari pour son aide et sa
patience. Je remercie \'egalement Mikhail Borovoi pour ses pr\'ecieux commentaires.

\section{Quelques pr\'eliminaires sur les complexes de tores}

Si $A$ est un groupe topologique ab\'elien, on note $A^{\wedge}$ son compl\'et\'e pour la topologie des sous-groupes ouverts d'indice fini, et $A_{\wedge}$ d\'esigne le groupe $\varprojlim_n A / n$ (o\`u $A / n := A / n A$ par d\'efinition). En outre, pour un groupe ab\'elien $A$, on note $_n A$ le sous-groupe de $n$-torsion de $A$. On consid\`ere \'egalement le module de Tate de $A$, \`a savoir le groupe $T(A) := \varprojlim_n {_n A}$; de m\^eme, si $G$ est un sch\'ema en groupes de type multiplicatif et de type fini sur la base, on note $_n G$ le sous-groupe de type multiplicatif noyau de la multiplication par $n$ sur $G$ (voir \cite{SGA3}, proposition 2.2). Enfin, pour un groupe ab\'elien $A$ et un nombre premier $l$, on note $A\{l\}$ le sous-groupe de torsion $l$-primaire de $A$, $\overline{A}$ le quotient de $A $par son sous-groupe divisible maximal, $\overline{A}\{l\}$ le quotient de $A\{l\}$ par son sous-groupe divisible maximal, et $A^{(l)}$ la limite projective $\varprojlim_n A / l^n$.

Introduisons d\'esormais quelques notations concernant les complexes de tores de longueur $2$. Le contexte est le suivant : soit $S$ un sch\'ema. On se donne deux $S$-tores (au sens de \cite{SGA3}, Expos\'e IX, d\'efinition 1.3) $T_1$ et $T_2$, et un morphisme de $S$-tores $\rho : T_1 \rightarrow T_2$. On note $C := \left[ T_1 \xrightarrow{\rho} T_2 \right]$ le complexe de $S$-tores ainsi obtenu, o\`u $T_1$ est en degr\'e $-1$ et $T_2$ en degr\'e $0$. 

On note aussi $\widehat{T_i}$ le faisceau $\Z$-constructible sur $S$ dual de $T_i$  (on rappelle qu'un faisceau $\mathcal{F}$ pour la topologie \'etale sur $S$ est dit $\Z$-constructible s'il existe un rev\^etement \'etale fini d'un ouvert de $S$ sur lequel $\mathcal{F}$ est le faisceau constant associ\'e \`a un groupe ab\'elien de type fini et les tiges $\mathcal{F}$ hors de cet ouvert sont de type fini comme groupes ab\'eliens), et $\widehat{C}$ le complexe de faisceaux \'etales (localement constants) $\left[\widehat{T_2} \xrightarrow{\widehat{\rho}} \widehat{T_1} \right]$, o\`u $\widehat{\rho}$ est le morphisme dual de $\rho$. On travaille dans la cat\'egorie des faisceaux en groupes ab\'eliens sur le site fppf de la base $S$, et dans la cat\'egorie d\'eriv\'ee associ\'ee \`a la cat\'egorie des complexes born\'es de faisceaux fppf. 

On construit alors un accouplement naturel $C \otimes^{\L} \widehat{C}
\rightarrow \G_m[1]$, fonctoriel en $C$, qui prolonge l'accouplement bien connu $T \otimes \widehat{T} \rightarrow \G_m$ pour un tore $T$ : pour cela, on remarque que $\widehat{C}$ \'etant un complexe de faisceaux plats sur $S$, le produit tensoriel d\'eriv\'e $C \otimes^{\L} \widehat{C}$ co\"incide avec le complexe "produit tensoriel total" (voir \cite{Wei}, 10.5.5 et 10.6.2). L'objet $C \otimes^{\L} \widehat{C}$ est donc repr\'esent\'e par le complexe $\left[ T_1 \otimes \widehat{T_2} \rightarrow (T_1 \otimes \widehat{T_1}) \oplus (T_2 \otimes \widehat{T_2}) \rightarrow T_2 \otimes \widehat{T_1} \right]$, la premi\`ere fl\^eche \'etant $(t_1, \widehat{t_2}) \mapsto t_1 \otimes \widehat{\rho}(\widehat{t_2}) - \rho(t_1) \otimes \widehat{t_2}$ et la seconde $t_1 \otimes \widehat{t_1} + t_2 \otimes \widehat{t_2} \mapsto \rho(t_1) \otimes \widehat{t_1} + t_2 \otimes \widehat{\rho}(\widehat{t_2})$. On dispose alors du morphisme canonique $ (T_1 \otimes \widehat{T_1}) \oplus (T_2 \otimes \widehat{T_2}) \rightarrow \G_m \oplus \G_m$, qui compos\'e avec le morphisme \fct{\G_m \oplus \G_m}{\G_m}{(t, t')}{t.t'} fournit un morphisme $ (T_1 \otimes \widehat{T_1}) \oplus (T_2 \otimes \widehat{T_2}) \rightarrow \G_m $. 

Celui-ci induit clairement un morphisme de complexes 
$$C \otimes^{\L} \widehat{C} \rightarrow \G_m[1]$$
 On v\'erifie ais\'ement que dans le cas o\`u $C = [0 \rightarrow T]$ ou $C = [T \rightarrow 0 ]$, on retrouve bien l'accouplement usuel entre un tore $T$ et son dual $\widehat{T}$.

\begin{rem}
{\rm
Cet accouplement correspond \`a une biextension naturelle de $(C, \widehat{C})$ par $\G_m$ (au sens de \cite{Del2}, 10.2.1), via la bijection canonique $\textup{Biext}(C, \widehat{C} ; \G_m) \cong \textup{Hom}(C \otimes^{\L} \widehat{C}, \G_m[1])$. 
}
\end{rem}

\paragraph{R\'ealisations $n$-adiques}

On d\'efinit ici les r\'ealisations $n$-adiques d'un complexe $C = [T_1 \xrightarrow{\rho} T_2]$ sur $S$.
\begin{defi}
On pose, pour $n \geq 1$, 
$T_{\Z / n}(C) := H^0(C[-1] \otimes^{\L} \Z / n)$ et $T_{\Z / n}(\widehat{C}) := H^0(\widehat{C}[-1] \otimes^{\L} \Z / n)$.
\end{defi}

Montrons alors le lemme suivant :
\begin{lem}
\label{real}
Le faisceau fppf $ T_{\Z / n}(C) $ est repr\'esentable par
un sch\'ema en groupes de type multiplicatif fini sur $S$, et $ T_{\Z
  / n}(\widehat{C})$ est repr\'esentable par le groupe
constant tordu fini dual
$\underline{\textup{Hom}}_{S\textup{-gr}}(T_{\Z / n}(C) ,
\G_m)$, et le produit tensoriel d\'eriv\'e $C \otimes^{\L}
\Z / n$ s'ins\`ere dans les triangles exacts suivants, fonctoriels en $C$, dans la cat\'egorie d\'eriv\'ee des faisceaux ab\'eliens fppf sur $S$ :
$$_n (\textup{Ker } \rho) [2] \rightarrow C \otimes^{\L} \Z / n \rightarrow T_{\Z / n}(C) [1] \rightarrow _n (\textup{Ker }\rho) [3] $$
et
$$T_{\Z / n}(\widehat{C}) [1] \rightarrow \widehat{C} \otimes^{\L} \Z / n \rightarrow \widehat{_n (\textup{Ker } \rho)} \rightarrow T_{\Z / n}(\widehat{C}) [2]$$
\end{lem}

\begin{dem}
En utilisant la r\'esolution plate $(\Z \xrightarrow{n} \Z)$ de $\Z / n$, on voit que $C \otimes^{\L} \Z / n$ s'identifie au complexe $\left[ T_1 \xrightarrow{n \oplus \rho} T_1 \oplus T_2 \xrightarrow{(t_1, t_2) \mapsto \rho(t_1) / t_2^n} T_2 \right]$, o\`u $T_1$ est en degr\'e $-2$ (voir par exemple \cite{Wei}, lemme 10.6.2). On note $\rho - n$ le second morphisme. On voit imm\'ediatemment que $\rho - n$ est surjectif, par cons\'equent ce complexe est quasi-isomorphe \`a $\left[ T_1 \xrightarrow{n \oplus \rho} \textup{Ker}(\rho-n) \right]$. Or on dispose du diagramme commutatif suivant, dont les lignes sont exactes :
\begin{displaymath}
\xymatrix{
0 \ar[r] & _n (\textup{Ker } \rho)  \ar[r] \ar[d]^{n \oplus \rho} & T_1 \ar[r]^{n \oplus \rho} \ar[d]^{n \oplus \rho} & \textup{Im}(n \oplus \rho) \ar[r] \ar[d] & 0 \\
0 \ar[r] & 0 \ar[r] & \textup{Ker}(\rho-n) \ar[r]^{=} & \textup{Ker}(\rho-n) \ar[r] & 0
}
\end{displaymath}
Or par d\'efinition, $T_{\Z / n}(C) = \textup{Ker}(\rho-n) / \textup{Im}(n \oplus \rho)$, et la troisi\`eme fl\`eche verticale est injective, donc le complexe $\left[ \textup{Im}(n \oplus \rho) \rightarrow \textup{Ker}(\rho-n) \right]$ est quasi-isomorphe \`a $T_{\Z / n}(C)$, donc ce diagramme s'identifie au premier triangle exact du lemme, dans la cat\'egorie d\'eriv\'ee (voir par exemple \cite{Wei}, 10.4.9). 
En ce qui concerne la finitude de $T_{\Z / n}(C)$, deux applications successives du lemme du serpent assurent que l'on a une suite exacte de faisceaux $_n T_1 \xrightarrow{\rho} _n T_2 \rightarrow \textup{Ker}(\rho-n) / \textup{Im}(n \oplus \rho) \rightarrow 0$. Or la cat\'egorie des $S$ sch\'emas en groupes de type multiplicatif de type fini est une cat\'egorie ab\'elienne (voir \cite{SGA3}, Expos\'e 9, corollaire 2.8), et $_n T_i$ est un-$S$ sch\'ema en groupes de type multiplicatif fini (voir \cite{SGA3}, Expos\'e 9, proposition 2.2), donc $T_{\Z / n}(C)$, qui est isomorphe au quotient $_n T_2 / \rho(_n T_1)$, est bien un $S$-sch\'ema en groupes de type multiplicatif fini. \\

Montrons d\'esormais le second triangle exact :
le complexe $\widehat{C} \otimes^{\L} \Z / n$ s'identifie cette fois \`a $\left[\widehat{T_2} \xrightarrow{n \oplus \widehat{\rho}} \widehat{T_2} \oplus \widehat{T_1} \xrightarrow{\widehat{\rho} - n} \widehat{T_1} \right]$, avec $\widehat{T_2}$ en degr\'e $-2$. Or $\widehat{T_2}$ est sans torsion, donc $n \oplus \widehat{\rho}$ est injectif. Donc $\widehat{C} \otimes^{\L} \Z / n$ est quasi-isomorphe au complexe $\left[ \left(\widehat{T_2} \oplus  \widehat{T_1} \right) / \textup{Im}(n \oplus \widehat{\rho}) \xrightarrow {\widehat{\rho} - n} \widehat{T_1} \right]$. On en d\'eduit que $T_{\Z / n}(\widehat{C}) = \textup{Ker}(\widehat{\rho} - n) / \textup{Im}(n \oplus \widehat{\rho})$, ainsi que le diagramme commutatif suivant, \`a lignes exactes :
\begin{displaymath}
\xymatrix{
0 \ar[r] & T_{\Z / n}(\widehat{C}) \ar[r] \ar[d] & \left(\widehat{T_2} \oplus  \widehat{T_1} \right) / \textup{Im}(n \oplus \widehat{\rho}) \ar[d]^{\widehat{\rho} - n} \ar[rr]^{\widehat{\rho} - n} & & \textup{Im}(\widehat{\rho} - n) \ar[r] \ar[d] & 0 \\
0 \ar[r] & 0 \ar[r] & \widehat{T_1} \ar[rr]^{=} & & \widehat{T_1} \ar[r] & 0
}
\end{displaymath}
qui fournit exactement le second triangle exact du lemme, \'etant donn\'e que le complexe $\left[ \textup{Im}(\widehat{\rho} - n) \hookrightarrow \widehat{T_1} \right]$ est quasi-isomorphe \`a $\widehat{T_1} / \textup{Im}(\widehat{\rho} - n)$, qui n'est autre que $\widehat{_n (\textup{Ker } \rho)}$ (car $\widehat{(\rho \oplus -( n))} = \widehat{\rho} - n$ et $\textup{Ker}(\rho \oplus (-n)) = _n \textup{Ker}(\rho)$).

Montrons d\'esormais que $T_{\Z / n}(C)$ et $T_{\Z/ n}(\widehat{C})$ se correspondent via l'\'equivalence de cat\'egorie entre $S$-sch\'emas en groupes de type multiplicatif de type fini et $S$-groupes constants tordus finiment engendr\'es (voir \cite{SGA3}, Expos\'e 10, corollaire 5.9) : on remarque que l'accouplement canonique $C \otimes^{\L} \widehat{C} \rightarrow \G_m[1]$ d\'efini plus haut induit naturellement un accouplement 
$$(C \otimes^{\L} \Z / n) \otimes^{\L} (\widehat{C} \otimes^{\L} \Z / n) \rightarrow \G_m[2]$$
d'o\`u en particulier un accouplement $T_{\Z / n}(C) \times T_{\Z / n}(\widehat{C}) \rightarrow \G_m$. Or on a montr\'e par application du lemme du serpent que l'on avait une suite exacte naturelle
$$_n T_1 \xrightarrow{\rho} _n T_2 \rightarrow T_{\Z / n}(C) \rightarrow 0$$
De m\^eme, on montre facilement la suite exacte suivante :
$$0 \rightarrow T_{\Z / n}(\widehat{C}) \rightarrow \widehat{T_2} / n \xrightarrow{\widehat{\rho}} \widehat{T_1} / n$$
Alors on conclut imm\'ediatement que $T_{\Z / n}(\widehat{C})$ est constant tordu, et qu'il s'identifie au dual de $T_{\Z / n}(C)$, en remarquant que $_n T_i$ est le dual de $\widehat{T_i} / n$ pour $i = 1, 2$.
\end{dem}

\begin{rem}
{\rm
Les calculs pr\'ec\'edents montrent en fait que $C \otimes^{\L} \Z / n$ s'identifie au c\^one $C(\rho_n)$ du morphisme de complexes $\left[T_1 \xrightarrow{n} T_1 \right] \xrightarrow{\rho_n} \left[ T_2 \xrightarrow{n} T_2 \right]$. Il se trouve que les deux complexes apparaissant ici sont respectivement quasi-isomorphes \`a $_n T_1[1]$ et $_n T_2 [1]$, on a donc un triangle exact dans la cat\'egorie d\'eriv\'ee :
$$_n T_1 [1] \xrightarrow{\rho_n} _n T_2 [1] \rightarrow C \otimes^{\L} \Z / n \rightarrow _n T_1 [2] $$
de m\^eme, dualement on obtient un triangle exact :
$$\widehat{T_2} / n \xrightarrow{\widehat{\rho}_n} \widehat{T_1} / n \rightarrow \widehat{C} \otimes^{\L} \Z /n \rightarrow \widehat{T_2} / n [1]$$
Ces triangles exacts sont des variantes des triangles exacts pr\'ec\'edents, les suites d'hypercohomologie ne sont pas les m\^emes, mais une partie de la suite du raisonnement peut se faire en utilisant ces variantes, les arguments de d\'evissage \'etant exactement les m\^emes, puisque les sch\'emas en groupes de type multiplicatif $_n T_i$, $_n \textup{Ker}(\rho)$ et $T_{\Z / n}(C)$ sont finis. On utilise dans la suite de cette section uniquement la version du lemme, qui a l'avantage de pr\'esenter une analogie avec \cite{HSz}. En outre, la version du lemme fait intervenir le noyau $\textup{Ker}(\rho)$, sur lequel vont porter des hypoth\`eses de finitude dans les sections suivantes, hypoth\`ese qui sera cruciale pour certains des r\'esultats de type Poitou-Tate qui vont suivre (voir section \ref{section PT}.
}
\end{rem}

\begin{rem}
{\rm
Dans sa th\`ese \cite{Jos}, P. Jossen a pr\'esent\'e un formalisme
tr\`es g\'en\'eral pour construire les modules de Tate $l$-adiques de
complexes (dits ``mod\'er\'es'') de faisceaux fppf sur une base
quelconque comme objets dans une cat\'egorie d\'eriv\'ee de faisceaux
$l$-divisibles localement consants (voir le chapitre 2 de
\cite{Jos}). En particulier, son travail s'applique \`a des complexes
de tores. Les constructions de cette section, ainsi que certains des
r\'esultats qui suivent, peuvent se reformuler dans le langage qu'il a
d\'evelopp\'e, et se voir ainsi de fa\c con plus naturelle.
}
\end{rem}

\section{Th\'eor\`emes de dualit\'e locale}
\label{section locale}
Soit $K$ un corps complet pour une valuation discr\`ete, \`a corps r\'esiduel fini $\F$, $\mathcal{O}$ son anneau des entiers. On se donne deux $K$-tores $T_1$ et $T_2$, et un morphisme de $K$-tores alg\'ebriques $\rho : T_1 \rightarrow T_2$. On note $C := \left[ T_1 \xrightarrow{\rho} T_2 \right]$ le complexe de tores ainsi obtenu, o\`u $T_1$ est en degr\'e $-1$ et $T_2$ en degr\'e $0$.

\paragraph{Topologie}

On munira $\H^i(K, C)$ de la topologie discr\`ete, sauf pour $i = -1, 0$ : pour $i = -1$, $\H^{-1}(K, C)$ est muni de la topologie induite par celle de $T_1(K)$ via l'identification $\H^{-1}(K, C) = \textup{Ker}(T_1(K) \rightarrow T_2(K))$. Pour $i = 0$, on consid\`ere la suite exacte de groupes ab\'eliens
$$T_1(K) \xrightarrow{\rho} T_2(K) \rightarrow \H^0(K, C) \rightarrow H^1(K, T_1)$$
On sait que l'image de $T_1(K)$ par $\rho$ s'identifie \`a un sous-groupe ferm\'e de $T_2(K)$, et par cons\'equent le quotient topologique $T_2(K) / \rho(T_1(K))$ est un groupe topologique s\'epar\'e. La suite exacte permet d'identifier ce quotient \`a un sous-groupe d'indice fini de $\H^0(K, C)$ (le groupe $H^1(K, T_2)$ est fini par \cite{Mil}, I, th\'eor\`eme 2.1), et on munit ce sous-groupe d'indice fini de la topologie quotient sur $T_2(K) / \rho(T_1(K))$, ce qui d\'efinit une topologie sur $\H^0(K, C)$. Par d\'efinition, le morphisme $T_2(K) \rightarrow \H^0(K, C)$ est alors continu et ouvert, et $\H^0(K, C)$ est s\'epar\'e.

Remarquons que l'accouplement $C \otimes^{\L} \widehat{C} \rightarrow \G_m[1]$ induit, via le cup produit, un morphisme sur les groupes d'hypercohomologie $\H^i(K, C) \times \H^{1-i}(K, \widehat{C}) \rightarrow H^1(K, \G_m[1]) \cong H^2(K, \G_m) \xrightarrow{j_K} \Q / \Z$
o\`u $j_K$ est l'invariant donn\'e par la th\'eorie du corps de classes local.
Montrons alors le r\'esultat suivant :

\begin{theo}
\label{dualite locale}
Le cup-produit $\H^i(K, C) \times \H^{1-i}(K, \widehat{C}) \rightarrow
\Q / \Z$ r\'ealise des dualit\'es parfaites, fonctorielles en $C$, entre les groupes suivants :
\begin{itemize}
	\item le groupe profini $\H^{-1}(K, C)^{\wedge}$ et le groupe discret $\H^2(K, \widehat{C})$.
	\item le groupe profini $\H^0(K, C)^{\wedge}$ et le groupe discret $\H^1(K, \widehat{C})$.
	\item le groupe discret $\H^1(K, C)$ et le groupe profini $\H^0(K, \widehat{C})^{\wedge}$.
	\item le groupe discret $\H^2(K, C)$ et le groupe profini $\H^{-1}(K, \widehat{C})^{\wedge}$.
\end{itemize}
\end{theo}

\begin{dem}
Ce r\'esultat s'obtient essentiellement par d\'evissage \`a partir de la dualit\'e locale de Tate-Nakayama pour les tores (voir \cite{Mil}, corollaires I.2.3 et I.2.4).
\begin{itemize}
	\item $i = -1$. \\
Consid\'erons la suite exacte suivante :
$$0 \rightarrow \H^{-1}(K, C) \rightarrow H^0(K,  T_1) \xrightarrow{\rho} H^0(K, T_2)$$
o\`u $\H^{-1}(K, C)$ s'identifie au sous-groupe topologique ferm\'e $\textup{Ker}(T_1(K) \xrightarrow{\rho} T_2(K))$ du groupe s\'epar\'e compactement engendr\'e, totalement discontinu et localement compact $T_1(K)$. De m\^eme, $T_2(K)$ est s\'epar\'e compactement engendr\'e, totalement discontinu et localement compact, et l'image de $T_1(K)$ s'identifie \`a un sous-groupe ferm\'e de $T_2(K)$, donc on est donc bien dans le cadre de la proposition de l'appendice de \cite{HSz}, et par cons\'equent la suite compl\'et\'ee
$$0 \rightarrow \H^{-1}(K, C)^{\wedge} \rightarrow H^0(K,  T_1)^{\wedge} \xrightarrow{\rho} H^0(K, T_2)^{\wedge}$$
reste exacte (l'inclusion $\H^{-1}(K, C) \rightarrow H^0(K,  T_1)$ est bien stricte par d\'efinition de la topologie sur $\H^{-1}(K, C)$).

On dispose du diagramme commutatif suivant, \`a lignes exactes (on peut dualiser la suite exacte pr\'ec\'edente car les groupes apparaissant sont profinis) :
\begin{displaymath}
\xymatrix{
H^2(K, \widehat{T_2}) \ar[r] \ar[d]^{\simeq} & H^2(K, \widehat{T_1}) \ar[r] \ar[d]^{\simeq} & \H^2(K, \widehat{C}) \ar[d] \ar[r] & H^3(K, \widehat{T_2})  \\
(H^0(K, T_2)^{\wedge})^{D} \ar[r] & (H^0(K, T_1)^{\wedge})^{D} \ar[r] & (\H^{-1}(K, C)^{\wedge})^{D} \ar[r] & 0
}
\end{displaymath}
Or $H^3(K, \widehat{T_2}) = 0$ car $K$ est de dimension cohomologique stricte \'egale \`a $2$ (voir \cite{Ser}, I.5.3, proposition 15), et les deux premi\`eres fl\`eches verticales sont des isomorphismes (voir \cite{Mil}, corollaire I.2.4). Donc une chasse au diagramme assure que l'on a un isomorphisme $\H^2(K, \widehat{C}) \simeq  (\H^{-1}(K, C)^{\wedge})^{D}$. Dans l'autre sens, on consid\`ere le diagramme exact (la seconde ligne est exacte car c'est la suite duale d'une suite exacte de groupes discrets) :
\begin{displaymath}
\xymatrix{
0 \ar[r] &\H^{-1}(K, C)^{\wedge} \ar[r] \ar[d] & H^0(K,  T_1)^{\wedge} \ar[r] \ar[d]^{\simeq} & H^0(K, T_2)^{\wedge} \ar[d]^{\simeq} \\
0 \ar[r] &  \H^2(K, \widehat{C}) ^{D} \ar[r] & H^2(K, \widehat{T_1})^{D} \ar[r] & H^2(K, \widehat{T_2}) ^{D}
}
\end{displaymath}
les deux isomorphismes provenant de la dualit\'e locale de Tate-Nakayama (voir \cite{Mil}, corollaire I.2.4). D'o\`u le point 1 du th\'eor\`eme.

	\item $i = 0$. \\
On dispose du diagramme suivant :

\begin{displaymath}
\xymatrix{
H^0(K,  T_1)^{\wedge} \ar[r] \ar[d]^{\simeq} & H^0(K, T_2)^{\wedge} \ar[d]^{\simeq} \ar[r] &  \H^0(K, C)^{\wedge} \ar[r] \ar[d] & H^1(K, T_1) \ar[r] \ar[d]^{\simeq} & H^1(K, T_2) \ar[d]^{\simeq} \\
H^2(K, \widehat{T_1})^{D} \ar[r] & H^2(K, \widehat{T_2}) ^{D} \ar[r] & \H^1(K, \widehat{C})^D \ar[r] & H^1(K, \widehat{T_1})^D \ar[r] & H^1(K, \widehat{T_2})^D
}
\end{displaymath}
La ligne inf\'erieure est exacte (puisque les groupes dont on prend le dual sont des groupes discrets), le morceau $H^0(K, T_2)^{\wedge} \rightarrow \H^0(K, C)^{\wedge}\rightarrow H^1(K, T_1) \rightarrow H^1(K, T_2)$ est exact puisque $H^1(K, T_1)$ est fini, et le morceau $H^0(K,  T_1)^{\wedge} \rightarrow H^0(K, T_2)^{\wedge} \rightarrow  \H^0(K, C)^{\wedge}$ est un complexe. Par cons\'equent, une chasse au diagramme assure que le morphisme $ \H^0(K, C)^{\wedge} \rightarrow \H^1(K, \widehat{C})^D$ est un isomorphisme. Dans l'autre sens, on consid\`ere le diagramme suivant :
\begin{displaymath}
\xymatrix{
H^1(K,  \widehat{T_2}) \ar[r] \ar[d]^{\simeq} & H^1(K, \widehat{T_1}) \ar[d]^{\simeq} \ar[r] &  \H^1(K, \widehat{C}) \ar[r] \ar[d] & H^2(K, \widehat{T_2}) \ar[r] \ar[d]^{\simeq} & H^2(K, \widehat{T_1}) \ar[d]^{\simeq} \\
H^1(K, T_2)^{D} \ar[r] & H^1(K, T_1) ^{D} \ar[r] & \H^0(K, C)^D \ar[r] & H^0(K, T_2)^D \ar[r] & H^0(K, T_1)^D
}
\end{displaymath}
La premi\`ere ligne est exacte. Concernant la seconde, la finitude des groupes $H^1(K, T_i)$, et la d\'efinition de la topologie sur $\H^0(K, C)$ assure que la suite suivante 
$$H^1(K, T_2)^{D} \rightarrow H^1(K, T_1) ^{D} \rightarrow \H^0(K, C)^D \rightarrow H^0(K, T_2)^D$$
est exacte, et le morceau $ \H^0(K, C)^D \rightarrow H^0(K, T_2)^D \rightarrow H^0(K, T_1)^D$ est un complexe. Alors une chasse au diagramme assure le point 2.

	\item $i = 1$. \\
	On fait le m\^eme raisonnement que pour $i = 0$ :
On dispose du diagramme suivant, commutatif \`a lignes exactes :
\begin{displaymath}
\xymatrix{
H^1(K,  T_1) \ar[r] \ar[d]^{\simeq} & H^1(K, T_2) \ar[d]^{\simeq} \ar[r] &  \H^1(K, C) \ar[r] \ar[d] & H^2(K, T_1) \ar[r] \ar[d]^{\simeq} & H^2(K, T_2) \ar[d]^{\simeq} \\
H^1(K, \widehat{T_1})^{D} \ar[r] & H^1(K, \widehat{T_2}) ^{D} \ar[r] & \H^0(K, \widehat{C})^D \ar[r] & H^0(K, \widehat{T_1})^D \ar[r] & H^0(K, \widehat{T_2})^D
}
\end{displaymath}
donc le morphisme $ \H^1(K, C) \rightarrow \H^0(K, \widehat{C})^D$ est un isomorphisme. Dans l'autre sens, on consid\`ere le diagramme suivant :
\begin{displaymath}
\xymatrix{
H^0(K,  \widehat{T_2})^{\wedge} \ar[r] \ar[d]^{\simeq} & H^0(K, \widehat{T_1})^{\wedge} \ar[d]^{\simeq} \ar[r] &  \H^0(K, \widehat{C})^{\wedge} \ar[r] \ar[d] & H^1(K, \widehat{T_2}) \ar[r] \ar[d]^{\simeq} & H^1(K, \widehat{T_1}) \ar[d]^{\simeq} \\
H^2(K, T_2)^{D} \ar[r] & H^2(K, T_1) ^{D} \ar[r] & \H^1(K, C)^D \ar[r] & H^1(K, T_2)^D \ar[r] & H^1(K, T_1)^D
}
\end{displaymath}
et la premi\`ere ligne est un complexe, et elle est exacte au niveau du morceau $H^0(K, \widehat{T_1})^{\wedge} \rightarrow  \H^0(K, \widehat{C})^{\wedge} \rightarrow H^1(K, \widehat{T_2}) \rightarrow H^1(K, \widehat{T_1})$ (puisque $H^1(K, \widehat{T_i})$ est fini), et on conclut par une chasse au diagramme.

	\item $i = 2$. \\
Comme dans le cas $i = -1$, on remarque que la suite suivante
$$0 \rightarrow \H^{-1}(K, \widehat{C})^{\wedge} \rightarrow H^0(K,  \widehat{T_2})^{\wedge} \xrightarrow{\widehat{\rho}} H^0(K,  \widehat{T_1})^{\wedge} $$
est exacte, puisque l'inclusion $ \H^{-1}(K, \widehat{C}) \rightarrow H^0(K,  \widehat{T_2})$ est stricte par d\'efinition de la topologie sur $\H^{-1}(K, \widehat{C}) $ et puisque les r\'eseaux $H^0(K,  \widehat{T_i})$ sont discrets et de type fini, donc on est bien dans le cadre la proposition de l'appendice de \cite{HSz}.
On dispose donc du diagramme commutatif suivant, \`a lignes exactes :
\begin{displaymath}
\xymatrix{
H^2(K, T_1) \ar[r] \ar[d]^{\simeq} & H^2(K, T_2) \ar[r] \ar[d]^{\simeq} & \H^2(K, C) \ar[d] \ar[r] & H^3(K, T_1) = 0 \\
(H^0(K, \widehat{T_1})^{\wedge})^{D} \ar[r] & (H^0(K,  \widehat{T_2})^{\wedge})^{D} \ar[r] & (\H^{-1}(K, \widehat{C})^{\wedge})^{D} \ar[r] & 0
}
\end{displaymath}
d'o\`u un isomorphisme $\H^2(K, C) \simeq  (\H^{-1}(K, \widehat{C})^{\wedge})^{D}$. Dans l'autre sens, on consid\`ere le diagramme \`a lignes exactes (les groupes $H^2(K, T_i)$ et $\H^2(K, C)$ sont discrets) :
\begin{displaymath}
\xymatrix{
0 \ar[r] &\H^{-1}(K, \widehat{C})^{\wedge} \ar[r] \ar[d] & H^0(K,  \widehat{T_2})^{\wedge} \ar[r] \ar[d]^{\simeq} & H^0(K,  \widehat{T_1})^{\wedge} \ar[d]^{\simeq} \\
0 \ar[r] &  \H^2(K, C) ^{D} \ar[r] & H^2(K, T_2)^{D} \ar[r] & H^2(K, T_1) ^{D}
}
\end{displaymath}
les deux isomorphismes provenant de la dualit\'e locale de Tate-Nakayama. D'o\`u le point 4 du th\'eor\`eme.
\end{itemize}
\end{dem}

\'Etudions maintenant le comportement des groupes de cohomologie non ramifi\'ee vis-\`a-vis de cet accouplement.
On se place dans la situation suivante : soient $\mathcal{T}_1$ et $\mathcal{T}_2$ deux tores sur $\textup{Spec } \mathcal{O}$, et $\rho : \mathcal{T}_1 \rightarrow \mathcal{T}_2$ un morphisme de $\textup{Spec } \mathcal{O}$-sch\'emas en groupes. On note $\mathcal{C}$ le complexe de tores associ\'e, et on note aussi $T_1$ (resp. $T_2$, resp. $C$) la fibre g\'en\'erique de $\mathcal{T}_1$ (resp. $\mathcal{T}_2$, resp. $\mathcal{C}$, i.e. le complexe de $K$-tores $\left[ T_1 \xrightarrow{\rho} T_2 \right]$). On d\'efinit alors $\H^i_{\textup{nr}}(K, C)$ comme l'image de $\H^i(\mathcal{O}, \mathcal{C})$ dans $\H^i(K, C)$. On rappelle que dans le cas des tores, le groupe $H^1(\mathcal{O}, \mathcal{T}_i)$ est trivial (c'est la conjonction de \cite{Mil2}, III.3.11 a) et du th\'eor\`eme de Lang). \\
Montrons d'abord le lemme suivant :

\begin{lem}
\label{lem inj local}
Soit $\mathcal{F}$ un faisceau (\'etale) localement constant $\Z$-constructible sans torsion sur $\textup{Spec } \mathcal{O}$. \\
Alors le morphisme $H^2(\mathcal{O}, \mathcal{F}) \rightarrow H^2(K, \mathcal{F})$ est injectif.
\end{lem}

\begin{dem}
On note $K^{\textup{nr}}$ l'extension maximale non ramifi\'ee de $K$, et $\mathcal{O}^{\textup{nr}}$ son anneau des entiers.
On consid\`ere les suites spectrales de Hochschild-Serre suivantes :
$H^p(\F, H^q(\mathcal{O}^{\textup{nr}}, \mathcal{F})) \Rightarrow H^{p+q}(\mathcal{O}, \mathcal{F})$ et $H^p(\F, H^q(K^{\textup{nr}}, \mathcal{F})) \Rightarrow H^{p+q}(K, \mathcal{F})$. Ces deux suites spectrales induisent le diagramme commutatif suivant de suites exactes de bas degr\'e :
\begin{changemargin}{-1.5cm}{1cm}
\begin{displaymath}
\xymatrix{
H^0(\F, H^1(\mathcal{O}^{\textup{nr}}, \mathcal{F})) \ar[d] \ar[r] & H^2(\F, H^0(\mathcal{O}^{\textup{nr}}, \mathcal{F})) \ar[r] \ar[d] & \textup{Ker} \left( H^2(\mathcal{O}, \mathcal{F}) \rightarrow H^0(\F, H^2(\mathcal{O}^{\textup{nr}}, \mathcal{F})) \right) \ar[r] \ar[d] & H^1(\F, H^1(\mathcal{O}^{\textup{nr}}, \mathcal{F})) \ar[d] \\
H^0(\F, H^1(K^{\textup{nr}}, \mathcal{F})) \ar[r] & H^2(\F, H^0(K^{\textup{nr}}, \mathcal{F})) \ar[r] & \textup{Ker} \left( H^2(K, \mathcal{F}) \rightarrow H^0(\F, H^2(K^{\textup{nr}}, \mathcal{F})) \right) \ar[r] & H^1(\F, H^1(K^{\textup{nr}}, \mathcal{F}))
}
\end{displaymath}
\end{changemargin}
Or on sait que $H^1(\mathcal{O}^{\textup{nr}}, \mathcal{F}) = H^2(\mathcal{O}^{\textup{nr}}, \mathcal{F}) = 0$ car $\mathcal{O}^{\textup{nr}}$ est acyclique pour la topologie \'etale. De plus, $\mathcal{O}^{\textup{nr}}$ est simplement connexe (pour la topologie \'etale), donc $\mathcal{F}$ est un faisceau constant $\Z^k$ sur $\mathcal{O}^{\textup{nr}}$, et donc on a $H^1(K^{\textup{nr}}, \mathcal{F}) = 0$ puisque l'on sait que $H^1(K^{\textup{nr}}, \Z) = 0$. Finalement, le diagramme pr\'ec\'edent devient :
\begin{displaymath}
\xymatrix{
H^2(\F, H^0(\mathcal{O}^{\textup{nr}}, \mathcal{F})) \ar[rrr]^{\simeq} \ar[d] & & & H^2(\mathcal{O}, \mathcal{F}) \ar[d] \\
H^2(\F, H^0(K^{\textup{nr}}, \mathcal{F})) \ar[rrr]^{\simeq} & & & \textup{Ker} \left( H^2(K, \mathcal{F}) \rightarrow H^0(\F, H^2(K^{\textup{nr}}, \mathcal{F})) \right)
}
\end{displaymath}
Or la fl\`eche verticale de gauche est un isomorphisme, donc celle de droite \'egalement, ce qui assure le lemme.
\end{dem}

\begin{theo}
Dans la situation pr\'ec\'edente, dans la dualit\'e parfaite : $\H^0(K, C)^{\wedge} \times \H^1(K, \widehat{C}) \rightarrow \Q / \Z$, les sous-groupes $\H^0_{\textup{nr}}(K, C)$ et $\H^1_{\textup{nr}}(K, \widehat{C})$ sont les orthogonaux respectifs l'un de l'autre. De m\^eme, dans la dualit\'e $\H^1(K, C) \times \H^0(K, \widehat{C})^{\wedge} \rightarrow \Q / \Z$, les sous-groupes $\H^1_{\textup{nr}}(K, C)$ et $\H^0_{\textup{nr}}(K, \widehat{C})^{\wedge}$ sont les orthogonaux respectifs l'un de l'autre.	
\end{theo}

\begin{dem}
On montre seulement le premier point, le second est similaire. 
Puisque $\textup{Br } \mathcal{O} = 0$, les accouplements $\H^0(K, C)^{\wedge} \times \H^1(K, \widehat{C}) \rightarrow \Q / \Z$ et $\H^1(K, C) \times \H^0(K, \widehat{C})^{\wedge} \rightarrow \Q / \Z$ induisent des morphismes  $\H^1(K, \widehat{C}) / \H^1(\mathcal{O}, \widehat{\mathcal{C}}) \rightarrow \H^0(\mathcal{O}, \mathcal{C})^D$ et $\H^0(K, C)^{\wedge} / \H^0(\mathcal{O}, \mathcal{C})^{\wedge} \rightarrow \H^1(\mathcal{O}, \widehat{\mathcal{C}})^D$. Il suffit alors de montrer que ces morphismes sont injectifs. 
On consid\`ere d'abord le diagramme commutatif suivant, dont les lignes sont exactes :
\begin{displaymath}
\xymatrix{
H^1(\mathcal{O}, \widehat{\mathcal{T}_1}) \ar[r] \ar@{->>}[d] & \H^1(\mathcal{O}, \widehat{\mathcal{C}}) \ar[r] \ar[d] & H^2(\mathcal{O}, \widehat{\mathcal{T}_2}) \ar[r] \ar[d] & H^2(\mathcal{O}, \widehat{\mathcal{T}_1}) \ar[d] \\
H^1(K, \widehat{T_1}) \ar[r] \ar[d] & \H^1(K, \widehat{C}) \ar[r] \ar[d] & H^2(K, \widehat{T_2}) \ar[r] \ar[d] & H^2(K, \widehat{T_1}) \ar[d] \\
H^1(\mathcal{O}, \mathcal{T}_1)^D = 0 \ar[r] & \H^0(\mathcal{O}, \mathcal{C})^D \ar[r] & H^0(\mathcal{O}, \mathcal{T}_2)^D \ar[r] & H^0(\mathcal{O}, \mathcal{T}_1)^D
}
\end{displaymath}
En outre, la troisi\`eme colonne est exacte (voir par exemple \cite{HSz}, lemma 2.11). 
On se donne alors un \'el\'ement $\alpha \in \H^1(K, \widehat{C})$ s'envoyant sur $0$ dans $\H^0(\mathcal{O}, \mathcal{C})^D$. On note $\beta$ son image dans $H^2(K, \widehat{T_2})$. Puisque $\beta$ s'envoie sur $0$ dans $H^0(\mathcal{O}, \mathcal{T}_2)^D$, on sait que $\beta$ provient d'un $\beta^0 \in H^2(\mathcal{O}, \widehat{\mathcal{T}_2})$, d'image $\gamma^0 \in  H^2(\mathcal{O}, \widehat{\mathcal{T}_1})$. Par commutativit\'e du diagramme, $\gamma^0$ s'envoie sur $0$ dans $H^2(K, \widehat{T_1})$. Or le morphisme $H^2(\mathcal{O}, \widehat{\mathcal{T}_1}) \rightarrow H^2(K, \widehat{T_1})$ est injectif par le lemme \ref{lem inj local}, donc $\gamma^0 = 0$. Par cons\'equent, $\beta^0$ se rel\`eve en un $\delta^0 \in \H^1(\mathcal{O}, \widehat{\mathcal{C}})$, d'image $\delta \in \H^1(K, \widehat{C})$. Or $\alpha$ et $\delta$ ont tous les deux pour image $\beta \in H^2(K, \widehat{T_2})$, donc $\beta = \delta + \epsilon$, 	avec $\epsilon \in H^1(K, \widehat{T_1})$. Par surjectivit\'e de $H^1(\mathcal{O}, \widehat{\mathcal{T}_1}) \rightarrow H^1(K, \widehat{T_1})$ (voir par exemple \cite{HSz}, page 105, second diagramme; on peut aussi voir cela comme une cons\'equence de \cite{Mil}, th\'eor\`eme I.2.6 et du fait que $H^1(\mathcal{O}, \mathcal{T}_1) = 0$), $\epsilon$ se rel\`eve en $\epsilon^0 \in H^1(\mathcal{O}, \widehat{\mathcal{T}_1})$. Alors l'\'el\'ement $\delta^0 + \epsilon^0 \in \H^1(\mathcal{O}, \widehat{\mathcal{C}})$ s'envoie sur $\alpha \in \H^1(K, \widehat{C})$. Cela assure donc que le morphisme $\H^1(K, \widehat{C}) / \H^1(\mathcal{O}, \widehat{\mathcal{C}}) \rightarrow \H^0(\mathcal{O}, \mathcal{C})^D$ est injectif. \\
Montrons l'injectivit\'e de $\H^0(K, C)^{\wedge} / \H^0(\mathcal{O}, \mathcal{C})^{\wedge} \rightarrow \H^1(\mathcal{O}, \widehat{\mathcal{C}})^D$. On consid\`ere cette fois le diagramme suivant (les groupes $H^0(\mathcal{O}, \mathcal{T}_i)$ \'etant compacts) :
\begin{displaymath}
\xymatrix{
H^0(\mathcal{O}, \mathcal{T}_1) \ar[r] \ar[d] & H^0(\mathcal{O}, \mathcal{T}_2) \ar[r] \ar[d] & \H^0(\mathcal{O}, \mathcal{C}) \ar[r] \ar[d] & H^1(\mathcal{O}, \mathcal{T}_1) = 0 \ar[d] \\
H^0(K, T_1)^{\wedge} \ar[r] \ar[d] & H^0(K, T_2)^{\wedge} \ar[r] \ar[d] & \H^0(K, C)^{\wedge} \ar[r] \ar[d] & H^1(K, T_1) \ar[d] \\
H^2(\mathcal{O}, \widehat{\mathcal{T}_1})^D \ar[r] & H^2(\mathcal{O}, \widehat{\mathcal{T}_2})^D \ar[r] & \H^1(\mathcal{O}, \widehat{\mathcal{C}})^D \ar[r] & H^1(\mathcal{O}, \widehat{\mathcal{T}_1})^D
}
\end{displaymath}

Soit $\alpha \in \H^0(K, C)^{\wedge}$ s'envoyant sur $0$ dans $ \H^1(\mathcal{O}, \widehat{\mathcal{C}})^D$. On note $\beta$ son image dans $H^1(K, T_1)$. Par injectivit\'e de la fl\`eche $H^1(K, T_1) \rightarrow H^1(\mathcal{O}, \widehat{\mathcal{T}_1})^D$ (voir \`a nouveau \cite{Mil}, th\'eor\`eme I.2.6 et le fait que $H^1(\mathcal{O}, \mathcal{T}_1) = 0$ par exemple), on sait que $\beta = 0$. Donc par exactitude de la seconde ligne en $\H^0(K, C)^{\wedge}$, $\alpha$ se rel\`eve en un $\gamma \in H^0(K, T_2)^{\wedge}$. On note $\gamma'$ l'image de $\gamma$ dans $H^2(\mathcal{O}, \widehat{\mathcal{T}_2})^D$. Par exactitude de la ligne inf\'erieure, ce $\gamma'$ se rel\`eve en un \'el\'ement $\delta' \in H^2(\mathcal{O}, \widehat{\mathcal{T}_1})^D$. Par surjectivit\'e de la fl\`eche $H^0(K, T_1)^{\wedge} \rightarrow H^2(\mathcal{O}, \widehat{\mathcal{T}_1})^D$ (voir lemme \ref{lem inj local}), on rel\`eve $\delta'$ en un $\delta \in H^0(K, T_1)^{\wedge} $. L'image $\delta'$ de $\delta$ dans $H^0(K, T_2)^{\wedge}$ a m\^eme image que $\gamma$ dans $H^2(\mathcal{O}, \widehat{\mathcal{T}_2})^D$, donc par exactitude de la deuxi\`eme colonne (dualit\'e de Tate-Nakayama, voir lemme 2.11 de \cite{HSz}), $\gamma - \delta'$ se rel\`eve en un $\epsilon \in  H^0(\mathcal{O}, \mathcal{T}_2)$, dont on note $\epsilon'$ l'image dans $ \H^0(\mathcal{O}, \mathcal{C})$. Alors, puisque la deuxi\`eme ligne du diagramme est un complexe, $\epsilon'$ s'envoie sur $\alpha$ dans $ \H^0(K, C)^{\wedge}$. Cela assure l'injectivit\'e de $\H^0(K, C)^{\wedge} / \H^0(\mathcal{O}, \mathcal{C}) \rightarrow \H^1(\mathcal{O}, \widehat{\mathcal{C}})^D$, et cela conclut la preuve du th\'eor\`eme.

\end{dem}

Traitons \'egalement le cas du corps $\R$ : pour cela, on introduit les groupes de cohomologie modifi\'es \`a la Tate $\widehat{\H}^i(\R, C)$.

\begin{prop}
Soit $C$ un complexe de tores sur $\R$. Alors le cup-produit induit
une dualit\'e parfaite de groupes finis, fonctorielle en $C$ :
$$\widehat{\H}^0(\R, C) \times \widehat{\H}^1(\R, \widehat{C}) \rightarrow \Z / 2$$
\end{prop}

\begin{dem}
c'est un d\'evissage facile \`a partir du cas des tores.
\end{dem}

Pour finir, il convient de r\'ediger un analogue du premier th\'eor\`eme pour le corps des fractions d'un anneau hens\'elien, puis de comparer les groupes de cohomologie pour un corps hens\'elien et pour son compl\'et\'e, afin de remplacer les compl\'et\'es $\widehat{k}_v$ par les henselis\'es $k_v$ quand $k$ est un corps de nombres et $v$ une place de $k$ (voir la section suivante).

\begin{theo}
\label{theo hensel}
Soit $A$ un anneau hens\'elien, de corps r\'esiduel fini, $F$ son
corps des fractions. On suppose $F$ de caract\'eristique $0$. Soit $C$
un complexe de tores sur $F$. Alors le cup-produit induit des
dualit\'es parfaites, fonctorielles en $C$, entre les groupes :
\begin{itemize}
	\item $\H^{-1}(F, C)^{\wedge}$ et $\H^2(F, \widehat{C})$.
	\item $\H^0(F, C)^{\wedge}$ et $\H^1(F, \widehat{C})$.
	\item $\H^1(F, C)$ et $\H^0(F, \widehat{C})^{\wedge}$.
	\item $\H^2(F, C)$ et $\H^{-1}(F, \widehat{C})^{\wedge}$
\end{itemize}
\end{theo}

\begin{rem}
{\rm On munit $T_i(F)$ de la topologie induite par celle de $T_i(K)$, et $\H^0(F, C)$ de la topologie "naturelle". On va d'ailleurs montrer que pour le groupe ab\'elien $\H^0(F, C)$, sa compl\'etion profinie co\"incide avec sa compl\'etion pour la topologie des sous-groupes ouverts d'indice fini.}
\end{rem}

\begin{dem}

\begin{lem}
\label{lem hensel}
Soit $K$ la compl\'etion de $F$. Alors les morphismes canoniques $\H^i(F, C) \rightarrow \H^i(K, C)$ sont des isomorphismes pour $i \geq 1$, et $\H^0(F, C) \rightarrow \H^0(K, C)$ induit un isomorphisme $\H^0(F, C)^{\wedge} \rightarrow \H^0(K, C)^{\wedge}$. \\
Les morphismes $\H^i(F, \widehat{C}) \rightarrow \H^i(K, \widehat{C})$ sont des isomorphismes pour tout $i \geq 0$.
\end{lem}

\begin{dem}
On sait que $H^1(F, T_i) \cong H^1(K, T_i)$ par \cite{HSz}, lemme 2.7, et ces groupes sont finis. 
Par d\'efinition de la topologie sur $T_i(F)$, on remarque que le sous-groupe $n T_i(F)$ de $T_i(F)$, qui est d'indice fini dans $T_i(F)$, est ouvert : en effet, son adh\'erence dans $T_i(K)$ est exactement $n T_i(K)$ (th\'eor\`eme de Greenberg, voir \cite{Gre}. Voir \'egalement \cite{BLR}, section 3.6, corollaire 10 pour une formulation pr\'ecise du r\'esultat utilis\'e), qui est d'indice fini dans $T_i(K)$, donc ouvert ($K$ est de caract\'eristique $0$, voir \cite{Mil}, page 32), et il se trouve que $n T_i(F) = n T_i(K) \cap T_i(F)$ (voir \cite{HSz}, preuve du lemme 2.7), d'o\`u le fait que $n T_i(F)$ soit ouvert dans $T_i(F)$. Par cons\'equent, pour les groupes $T_i(F)$ (resp. $T_i(K)$), la compl\'etion pour les sous-groupes $n T_i(F)$ (resp. $n T_i(K)$)  co\"incide avec la compl\'etion pour les sous-groupes ouverts d'indice fini. Or $T_i(F)_{\wedge} \cong T_i(K)_{\wedge}$ (voir \cite{HSz}, preuve du lemme 2.7), donc $H^0(F,  T_i)^{\wedge} \cong H^0(K,  T_i)^{\wedge}$.

Or on peut compl\'eter la suite exacte
$$H^0(F,  T_1) \rightarrow H^0(F, T_2) \rightarrow  \H^0(F, C) \rightarrow H^1(F, T_1) \rightarrow H^1(F, T_2)$$
puisque les deux derniers groupes sont finis et le seconde fl\`eche est ouverte par d\'efinition de la toplogie sur $\H^0(F, C)$. D'o\`u le diagramme commutatif \`a lignes exactes suivant :
\begin{displaymath}
\xymatrix{
H^0(F,  T_1)^{\wedge} \ar[r] \ar[d]^{\simeq} & H^0(F, T_2)^{\wedge} \ar[d]^{\simeq} \ar[r] &  \H^0(F, C)^{\wedge} \ar[r] \ar[d] & H^1(F, T_1) \ar[r] \ar[d]^{\simeq} & H^1(F, T_2) \ar[d]^{\simeq} \\
H^0(K,  T_1)^{\wedge} \ar[r] & H^0(K, T_2)^{\wedge} \ar[r] &  \H^0(K, C)^{\wedge} \ar[r] & H^1(K, T_1) \ar[r] & H^1(K, T_2) \\
}
\end{displaymath}

On a d\'ej\`a montr\'e que les deux premi\`eres fl\`eches verticales, ainsi que les deux derni\`eres, \'etaient des isomorphismes, donc par une chasse au diagramme, la fl\`eche centrale est un isomorphisme.
Pour les groupes $\H^i(F, C)$ avec $i \geq 1$, la preuve est exactement celle du lemme 2.7 de \cite{HSz}.
Pour le complexe dual $\widehat{C}$, la preuve en degr\'e strictement positif est toujours la m\^eme; pour le degr\'e $0$, on consid\`ere le diagramme exact suivant :
\begin{displaymath}
\xymatrix{
H^0(F, \widehat{T_2}) \ar[r] \ar[d]^{\simeq} & H^0(F, \widehat{T_1}) \ar[d]^{\simeq} \ar[r] &  \H^0(F, \widehat{C}) \ar[r] \ar[d] & H^1(F, \widehat{T_2}) \ar[r] \ar[d]^{\simeq} & H^1(F, \widehat{T_1}) \ar[d]^{\simeq} \\
H^0(K,  \widehat{T_2}) \ar[r] & H^0(K, \widehat{T_1}) \ar[r] &  \H^0(K, \widehat{C}) \ar[r] & H^1(K, \widehat{T_2}) \ar[r] & H^1(K, \widehat{T_2}) \\
}
\end{displaymath}
les isomorphismes provenant du fait que $\widehat{T_i}$ est localement constant et que les corps $F$ et $K$ ont m\^eme groupe de Galois absolu. 
\end{dem}

Ce lemme \ref{lem hensel}, ainsi que la remarque sur les diff\'erentes compl\'etions dans la preuve du lemme, assure le th\'eor\`eme \ref{theo hensel} gr\^ace au th\'eor\`eme \ref{dualite locale}.

\end{dem}

\section{Dualit\'e globale : cohomologie \'etale}
\label{section etale}
Soit $k$ un corps de nombres, $\mathcal{O}_k$ son anneau des entiers. Soit $U$ un ouvert non vide de $\textup{Spec } \mathcal{O}_k$, et $\Sigma_f$ l'ensemble des places finies de $k$ correspondant \`a des points ferm\'es hors de $U$. Si $v$ d\'esigne un place de $k$, on note $k_v$ l'hens\'elis\'e de $k$ en $v$, et $\widehat{k}_v$ le compl\'et\'e de $k$ en $v$. Dans toute la suite, si $v$ est une place infinie de $k$, les groupes d'hypercohomolgie modifi\'es de Tate $\widehat{\H}^i(k_v, .)$ seront not\'es $\H^i(k_v, .)$. On note \'egalement $\Sigma := \Sigma_f \cup \Omega_{\infty}$, o\`u $\Omega_{\infty}$ d\'esigne l'ensemble des places infinies de $k$.

On renvoie au d\'ebut de la section 3 de \cite{HSz} pour la d\'efinition des groupes d'hypercohomologie \`a support compact \`a valeur dans un complexe de faisceaux ab\'eliens cohomologiquement born\'e. On utilisera \'egalement la propri\'et\'e de fonctorialit\'e covariante de la cohomologie \`a support compacte rappel\'ee au d\'ebut de la page 107 dans \cite{HSz} : si $V \rightarrow U$ est une immersion ouverte, on a un morphisme canonique $\H^i_c(V, \mathcal{C}_V) \rightarrow \H^i_c(U, \mathcal{C})$.

Suivant \cite{Mil}, on note $\D^i(U, .)$ l'image de $\H^i_c(U, .)$ dans $\H^i(U, .)$.

\begin{lem}
\label{lem struc etale}
Soit $\mathcal{C}$ un complexe de tores sur $U$.
\begin{enumerate}
	\item Les groupes $\H^i(U, \mathcal{C})$ et $\H^i(U, \widehat{\mathcal{C}})$ sont de torsion pour $i \geq 1$, ainsi que les groupes $\H^j_c(U, \widehat{\mathcal{C}})$ et $\H^j_c(U, \mathcal{C})$ pour $j \geq 2$.
	\item Pour tout $l$ inversible sur $U$, les groupes $\H^i(U, \mathcal{C})\{ l \}$ et $\H^i(U, \widehat{\mathcal{C}}) \{ l \}$ (resp. $\H^j_c(U, \widehat{\mathcal{C}}) \{ l \}$ et $\H_c^j(U, \mathcal{C})\{ l \}$ ) sont de cotype fini pour $i \geq 1$ (resp. $j \geq 2$).
	\item Les groupes $\H^{-1}(U, \mathcal{C})$, $\H^0(U, \mathcal{C})$, $\H^{-1}_c(U, \mathcal{C})$, $\H^{-1}(U, \widehat{\mathcal{C}})$, $\H^0(U, \widehat{\mathcal{C}})$, $\H^{-1}_c(U, \widehat{\mathcal{C}})$ et $\H^0_c(U, \widehat{\mathcal{C}})$ sont de type fini. 
        \item Le groupe $\H^0_c(U, \mathcal{C})$ est extension d'un groupe de type fini par un groupe profini, le groupe $\H^1_c(U, \mathcal{C})$ est extension d'un groupe de torsion (dont les sous-groupes de torsion $l$-primaire sont de cotype fini) par un groupe profini, et le groupe $\H^1_c(U, \widehat{\mathcal{C}})$ est extension d'un groupe de torsion (dont les sous-groupes de torsion $l$-primaire sont de cotype fini) par un groupe de type fini.
\end{enumerate}
\end{lem}

\begin{dem}
\begin{enumerate}
	\item Notons $\mathcal{C} := \left[ \mathcal{T}_1 \rightarrow \mathcal{T}_2 \right]$. On consid\`ere la suite exacte suivante :
$$H^i(U,\mathcal{T}_1) \rightarrow H^i(U, \mathcal{T}_2) \rightarrow \H^i(U, \mathcal{C}) \rightarrow H^{i+1}(U, \mathcal{T}_1) \rightarrow H^{i+1}(U, \mathcal{T}_2)$$
Or on sait que les $H^i(U, \mathcal{T}_r)$ sont de torsion pour $i \geq 1$ (cf \cite{HSz}, lemme 3.2.(1)), donc il en est de m\^eme pour $ \H^i(U, \mathcal{C})$. On en d\'eduit le r\'esultat pour les groupes de cohomologie \`a support compact \`a l'aide de la suite exacte 
\begin{changemargin}{-2cm}{1cm}
$$\H^i(U, \mathcal{C}) \rightarrow \bigoplus_{v \in \Sigma_f} \H^i(k_v, \mathcal{C}) \oplus \bigoplus_{v \textup{ infinie}} \widehat{\H}^i(k_v, \mathcal{C}) \rightarrow \H^{i+1}_c(U, \mathcal{C}) \rightarrow \H^{i+1}(U, \mathcal{C}) \rightarrow \bigoplus_{v \in \Sigma_f} \H^{i+1}(k_v, \mathcal{C}) \oplus \bigoplus_{v \textup{ infinie}} \widehat{\H}^{i+1}(k_v, \mathcal{C})$$
\end{changemargin}
Le raisonnement est exactement le m\^eme pour les groupes $\H^i(U, \widehat{\mathcal{C}})$, utilisant le fait que $H^i(U, \widehat{\mathcal{T}_r})$ est de torsion pour $i \geq 1$.

	\item Par d\'efinition, ces groupes sont de torsion. On utilise \`a nouveau la suite exacte :
$$H^i(U,\mathcal{T}_1)\{ l \} \rightarrow H^i(U, \mathcal{T}_2)\{ l \} \rightarrow \H^i(U, \mathcal{C})\{ l \} \rightarrow H^{i+1}(U, \mathcal{T}_1)\{ l \}$$
Pour $i \geq 1$, les trois groupes autres que $\H^i(U, \mathcal{C})\{ l \}$ sont de cotype fini (voir \cite{HSz}, lemme 3.2.(2)), donc $\H^i(U, \mathcal{C})\{ l \}$ l'est aussi (puisqu'un sous-groupe et un quotient d'un groupe de cotype fini sont de cotype fini). En ce qui concerne les groupes $\H^i(U, \widehat{\mathcal{C}})\{ l \}$, $\H^j_c(U,  \mathcal{C})\{ l \}$ et $\H^j_c(U, \widehat{\mathcal{C}})\{ l \}$, le raisonnement est le m\^eme.

	\item Pour $\H^{-1}(U, \mathcal{C})$ et $\H^0(U, \mathcal{C})$, on utilise \`a nouveau la suite exacte pr\'ec\'edente, sachant que $H^0(U, \mathcal{T}_i)$ est de type fini (voir \cite{HSz}, lemme 3.2.(3)) et que $H^1(U, \mathcal{T}_1)$ est fini (voir \cite{Mil}, th\'eor\`eme II.4.6). En ce qui concerne le groupe $\H^{-1}_c(U, \mathcal{C})$, on utilise la suite exacte reliant cohomologie \'etale et cohomologie \`a support compact. Pour les autres groupes, le raisonnement est similaire, en utilisant en outre le fait que $\H^{-1}(k_v, \widehat{\mathcal{C}})$ est de type fini pour $v \in \Sigma$.

        \item C'est un d\'evissage \`a l'aide des r\'esultats des points pr\'ec\'edents.
\end{enumerate}
\end{dem}

On va d\'esormais utiliser les r\'ealisations $n$-adiques du complexe $\mathcal{C}$ et de son dual.

On rappelle que l'accouplement canonique $\mathcal{C} \otimes^{\L} \widehat{\mathcal{C}} \rightarrow \G_m[1]$ induit $(\mathcal{C} \otimes^{\L} \Z / n) \otimes^{\L} (\widehat{\mathcal{C}} \otimes^{\L} \Z / n) \rightarrow \G_m[2]$. On en d\'eduit des accouplements en cohomologie :
\begin{eqnarray}
\label{accouplement}
\H^i(U, \mathcal{C} \otimes^{\L} \Z / n) \times \H^{1-i}_c(U, \widehat{\mathcal{C}} \otimes^{\L} \Z / n) \rightarrow H^3_c(U, \G_m) = \Q / \Z
\end{eqnarray}

On a alors le r\'esultat de dualit\'e suivant :

\begin{prop}
\label{dualite finie}
Pour tout $i \in \Z$, si l'entier $n$ est inversible sur $U$, l'accouplement (\ref{accouplement})
$$\H^i(U, \mathcal{C} \otimes^{\L} \Z / n) \times \H^{1-i}_c(U, \widehat{\mathcal{C}} \otimes^{\L} \Z / n) \rightarrow H^3_c(U, \G_m) = \Q / \Z$$ 
est une dualit\'e parfaite de groupes finis, fonctorielle en $C$. De m\^eme en inversant les r\^oles de $\mathcal{C}$ et $\widehat{\mathcal{C}}$.
\end{prop}

\begin{dem}
Gr\^ace au lemme \ref{real}, on dispose du diagramme commutatif suivant :
\begin{changemargin}{-2.3cm}{1cm}
\begin{displaymath}
\xymatrix{
H^i(U, T_{\Z / n}(\mathcal{C})) \ar[r] \ar[d] & H^{i+2}(U, {_n (\textup{Ker } \rho)}) \ar[r] \ar[d] & \H^i(U, \mathcal{C} \otimes^{\L} \Z / n) \ar[r] \ar[d] & H^{i+1}(U, T_{\Z / n}(\mathcal{C})) \ar[r] \ar[d] & H^{i+3}(U, {_n (\textup{Ker } \rho)}) \ar[d] \\
H^{3-i}_c(U, T_{\Z / n}(\widehat{\mathcal{C}}))^D \ar[r] & H^{1-i}_c(U, \widehat{_n (\textup{Ker } \rho)})^D \ar[r] & \H^{1-i}_c(U, \widehat{\mathcal{C}} \otimes^{\L} \Z / n)^D \ar[r] & H^{2-i}_c(U, T_{\Z / n}(\widehat{\mathcal{C}}))^D \ar[r] & H^{-i}_c(U, \widehat{_n (\textup{Ker } \rho)})^D
}
\end{displaymath}
\end{changemargin}

Hormis les groupes d'hypercohomologie au centre du diagramme, tous les groupes de cohomologie intervenant ici sont finis (voir par exemple \cite{Mil}, th\'eor\`eme II.3.1), ce qui assure la finitude des groupes intervenant dans l'\'enonc\'e de la proposition.

Par d\'efinition de l'accouplement (\ref{accouplement}) et des fl\`eches dans le lemme \ref{real}, les autres morphismes autres que le morphisme central sont (au signe pr\`es) les accouplements usuels qui interviennent dans la dualit\'e d'Artin-Verdier. On va d\'esormais appliquer le th\'eor\`eme de dualit\'e de Artin-Verdier pour les faisceaux finis.

Les faisceaux $T_{\Z / n}(C)$ et $_n \textup{Ker}(\rho)$ \'etant constructibles localement constants (voir lemme \ref{real}), et $n$ \'etant inversible sur $U$, les fl\`eches verticales autres que celle du milieu sont des isomorphismes (voir par exemple \cite{Mil}, corollaire II.3.3). Donc celle du milieu est aussi un isomorphisme, par lemme des cinq. En rempla\c cant $C$ par $\widehat{C}$, on obtient finalement une dualit\'e parfaite de groupes finis $\H^i(U, \mathcal{C} \otimes \Z / n) \times \H^{1-i}_c(U, \widehat{\mathcal{C}} \otimes \Z / n) \rightarrow \Q / \Z$, pour tout $n$ inversible sur $U$ et tout $i$, ce qui conclut la preuve.
\end{dem}

D\'eduisons de ce r\'esultat le th\'eor\`eme principal de ce paragraphe :

\begin{theo} Pour tout entier $i$, pour tout nombre premier $l$
  inversible sur $U$, le cup-produit induit une dualit\'e parfaite
  fonctorielle en $C$ :
$$\overline{\H^i(U, \mathcal{C})}\{l\} \times \overline{\H^{2-i}_c(U, \widehat{\mathcal{C}})}\{l\} \rightarrow \Q / \Z$$
et de m\^eme en \'echangeant les r\^oles de $\mathcal{C}$ et $ \widehat{\mathcal{C}}$.
\end{theo}

\begin{rem}
{\rm Les groupes $\H^i(U, \mathcal{C})\{l\}$ et $\H^{2-i}_c(U, \widehat{\mathcal{C}})\{l\}$ \'etant tous de cotype fini, leur sous-groupe divisible maximal co\"incide avec leur sous-groupe form\'e des \'el\'ements divisibles.
}
\end{rem}

\begin{dem}
 Fixons $l$ un nombre premier inversible sur $U$.
On d\'eduit de la proposition \ref{dualite finie} une dualit\'e parfaite entre le groupe discret $\varinjlim_n \H^{i-1}(U, \mathcal{C} \otimes \Z / {l^n})$ et le groupe profini $\varprojlim_n  \H^{2-i}_c(U, \widehat{\mathcal{C}} \otimes \Z / {l^n})$. Celle-ci induit une dualit\'e parfaite entre le quotient $\overline{\varinjlim_n \H^{i-1}(U, \mathcal{C} \otimes \Z / {l^n})}$ et le sous-groupe de torsion $\left( \varprojlim_n  \H^{2-i}_c(U, \widehat{\mathcal{C}} \otimes \Z / {l^n}) \right)_{\textup{tors}}$.

Or en tensorisant (produit tensoriel d\'eriv\'e $\otimes^{\L}$) la suite exacte $0 \rightarrow \Z \rightarrow \Z \rightarrow \Z / n \rightarrow 0$ par $\mathcal{C}$, on obtient une suite exacte de groupes finis :
$$0 \rightarrow \H^{i-1}(U, \mathcal{C}) \otimes \Z / n \rightarrow \H^{i-1}(U, \mathcal{C} \otimes^{\L} \Z / n) \rightarrow _n \H^i(U, \mathcal{C}) \rightarrow 0$$
qui permet d'identifier (apr\`es passage \`a la limite inductive) le groupe $\overline{\varinjlim_n \H^{i-1}(U, \mathcal{C} \otimes \Z / {l^n})}$ au groupe $\overline{\H^i(U, \mathcal{C})}\{ l \}$. De m\^eme, en rempla\c cant $\mathcal{C}$ par $\widehat{\mathcal{C}}$, et en passant \`a la limite projective (les groupes en question ici \'etant finis), on obtient un isomorphisme 
$$\left( \varprojlim_n  \H^{2-i}_c(U, \widehat{\mathcal{C}} \otimes \Z / {l^n}) \right)_{\textup{tors}} \cong \left( \H^{2-i}_c(U, \widehat{\mathcal{C}})^{(l)} \right)_{\textup{tors}} = \H^{2-i}_c(U, \widehat{\mathcal{C}})^{(l)} \{ l \}$$
Enfin, la structure des groupes $\H^{2-i}_c(U, \widehat{\mathcal{C}})$ (voir lemme \ref{lem struc etale}) fournit un isomorphisme $\H^{2-i}_c(U, \widehat{\mathcal{C}})^{(l)} \{ l \} \cong \overline{\H^{2-i}_c(U, \widehat{\mathcal{C}})} \{ l \}$ ($\H^{2-i}_c(U, \widehat{\mathcal{C}})^{(l)}$ est toujours de cotype fini).

D'o\`u finalement le r\'esultat, en remarquant que tout ce raisonnement fonctionne en intervertissant $\mathcal{C}$ et $\widehat{\mathcal{C}}$.

\end{dem}




On souhaite d\'eduire de ce r\'esultat un r\'esultat de dualit\'e sur les groupes $D^i(U, \mathcal{C})$ et $D^j(U, \widehat{\mathcal{C}})$. Pour cela on va utiliser le lemme suivant (voir \cite{HSz}, corrigenda) :

\begin{lem}
\label{lem erratum}
Soit $a \in \D^1(U, \mathcal{C})$, $l^r$-divisible dans $\H^1(U, \mathcal{C})$, et orthogonal au sous-groupe $_{l^r} \D^1(U, \widehat{\mathcal{C}})$ de $\H^1(U, \widehat{\mathcal{C}})$. Alors $a$ est $l^r$-divisible dans $\D^1(U, \mathcal{C})$.
\end{lem}

\begin{dem}
Notons $n = l^r$. On consid\`ere le diagramme commutatif exact suivant :
\begin{displaymath}
\xymatrix{
& & \H^1_c(U, \mathcal{C}) \ar[r] \ar[d]^{n} & \H^1(U, \mathcal{C}) \ar[d]^{n} \\
& \bigoplus_{v \in \Sigma} \H^0(k_v, \mathcal{C}) \ar[r] \ar[d] & \H^1_c(U, \mathcal{C}) \ar[r] \ar[d] & \H^1(U, \mathcal{C}) \ar[d] \\
\H^0(U, \mathcal{C} \otimes^{\L} \Z / n) \ar[r] & \bigoplus_{v \in \Sigma} \H^0(k_v, \mathcal{C} \otimes^{\L} \Z / n) \ar[r] & \H^1_c(U, \mathcal{C} \otimes^{\L} \Z / n) \ar[r] & \H^1(U, \mathcal{C} \otimes^{\L} \Z / n)
}
\end{displaymath}

Soit $a \in \D^1(U, \mathcal{C})$ comme dans l'\'enonc\'e, image d'un \'el\'ement $\tilde{a} \in \H^1_c(U, \mathcal{C})$. On \'ecrit $a = n a_1$, avec $a_1 \in \H^1(U, \mathcal{C})$. On note $a_2$ l'image de $\tilde{a}$ dans $\H^1_c(U, \mathcal{C} \otimes^{\L} \Z / n)$. Par commutativit\'e et exactitude, $a_2$ provient d'un \'el\'ement $(c_v) \in  \bigoplus_{v \in \Sigma} \H^0(k_v, \mathcal{C} \otimes^{\L} \Z / n)$.
On dispose alors des accouplements naturels : $\D^1(U, \mathcal{C}) \times \D^1(U, \widehat{\mathcal{C}}) \xrightarrow{P-T} \Q / \Z$ et $\H^1_c(U, \mathcal{C} \otimes^{\L} \Z / n) \times \H^0(U, \widehat{\mathcal{C}} \otimes^{\L} \Z / n) \xrightarrow{A-V} \Z / n$, qui sont compatibles via le diagramme commutatif :
\begin{displaymath}
\xymatrix{
\H^0(U, \widehat{\mathcal{C}} \otimes^{\L} \Z / n) \ar[d] & \times & \H^1_c(U, \mathcal{C}\otimes^{\L} \Z / n) \ar[rr]^{\textup{A-V}} & & \Q / \Z \ar[d]^{=} \\
_n \H^1(U, \widehat{\mathcal{C}}) & \times & \H^1_c(U, \mathcal{C}) / n \ar[u] \ar[d] \ar[rr] & & \Q / \Z \ar[d]^{=} \\
_n \D^1(U,  \widehat{\mathcal{C}}) \ar[u] & \times & \D^1(U, \mathcal{C}) / n \ar[rr]^{\textup{P-T}} & & \Q / \Z
}
\end{displaymath}
On en d\'eduit que $<a, a'>_{\textup{P-T}} = <a_2, b' >_{\textup{A-V}}$, si $a' \in {_n \D^1(U,  \widehat{\mathcal{C}})}$ est l'image de $b' \in \H^0(U, \widehat{\mathcal{C}} \otimes^{\L} \Z / n)$. Or par hypoth\`ese $<a_2, b' >_{\textup{A-V}} = \sum_{v \in \Sigma} <c_v, b'_v >_v$ (o\`u $b'_v$ est l'image de $b'$ dans $\bigoplus_{v \in \Sigma} \H^0(k_v,\widehat{\mathcal{C}} \otimes^{\L} \Z / n)$), et on a suppos\'e que $<a, a'>_{\textup{P-T}} = 0$ pour tout $a' \in {_n \D^1(U,  \widehat{\mathcal{C}})}$, donc $a_2$ est orthogonal \`a l'image de $S_n(U, \widehat{\mathcal{C}}) := \textup{Ker}( \H^0(U, \widehat{\mathcal{C}} \otimes^{\L} \Z / n) \rightarrow \bigoplus_{v \in \Sigma} \H^1(k_v, \widehat{\mathcal{C}}))$ dans $\bigoplus_{v \in \Sigma} \H^0(k_v, \widehat{\mathcal{C}} \otimes^{\L} \Z / n)$.

On va maintenant montrer que l'\'el\'ement $a_2$ v\'erifiant cette propri\'et\'e est n\'ecessairement nul, ce qui va permettre de conclure puisque si $a_2 = 0$, alors $\tilde{a}$ s'\'ecrit $\tilde{a} = n \tilde{a_1}$, avec $ \tilde{a_1} \in \H^1_c(U, \mathcal{C})$, ce qui conclut la preuve du lemme \ref{lem erratum}.

Pour montrer la trivialit\'e de $a_2$, on utilise l'analogue du lemme des corrigenda de \cite{HSz} : 
plus pr\'ecis\'ement, montrons que $(c_v)$ est somme d'un \'el\'ement de $\bigoplus_{v \in \Sigma} \H^0(k_v, C)$ et d'un \'el\'ement de $\H^0(U, C \otimes^{\L} \Z / n)$. Pour cela, regardons le diagramme suivant ($\bigoplus_{v \in \Sigma} \H^0(k_v, C \otimes^{\L} \Z / n)$ est fini) :
\begin{displaymath}
\xymatrix{
\bigoplus_{v \in \Sigma} \H^0(k_v, C \otimes^{\L} \Z / n) \ar[rd]^{\gamma_1} &  & & \\
\bigoplus_{v \in \Sigma} \H^0(k_v, C)^{\wedge} \ar[u] \ar[r] & \H^0(U, \widehat{\mathcal{C}} \otimes^{\L} \Z / n)^D \ar[r] & S_n(U, \widehat{\mathcal{C}})^D \ar[r] & 0
}
\end{displaymath}
On sait donc que $\gamma_1((c_v))$ est nul quand on le pousse dans $S_n(U, \widehat{\mathcal{C}})^D$, donc par exactitude de la derni\`ere ligne (elle est exacte par d\'efinition de $S_n(U, \widehat{\mathcal{C}})$), on en d\'eduit que $\gamma_1((c_v))$ se rel\`eve dans $\bigoplus_{v \in \Sigma} \H^0(k_v, C)^{\wedge}$, en un \'el\'ement not\'e $(t_v)$. On envoie alors $(t_v)$ dans $\bigoplus_{v \in \Sigma} \H^0(k_v, C \otimes^{\L} \Z / n)$, et on remarque que $(c_v) - (t_v) \in  \textup{Ker}(\gamma_1)$. Consid\'erons alors le diagramme commutatif suivant, dont la premi\`ere ligne est exacte :
\begin{displaymath}
\xymatrix{
\H^0(U, \mathcal{C} \otimes^{\L} \Z / n) \ar[r] & \bigoplus_{v \in \Sigma} \H^0(k_v, \mathcal{C} \otimes^{\L} \Z / n) \ar[r] \ar[d]^{\cong} & \H^1_c(U, \widehat{\mathcal{C}} \otimes^{\L} \Z / n)^D \ar[d]^{\cong} \\
& \bigoplus_{v \in \Sigma} \H^0(k_v, \widehat{\mathcal{C}} \otimes^{\L} \Z / n)^D \ar[r]^(.55){\gamma_1} & \H^0(U, \widehat{\mathcal{C}} \otimes^{\L} \Z / n)^D
}
\end{displaymath}
o\`u les isomorphismes verticaux proviennent du th\'eor\`eme \ref{theo hensel} de dualit\'e locale et la proposition \ref{dualite finie}.On d\'eduit donc de ce diagramme la suite exacte suivante : $\H^0(U, C \otimes^{\L} \Z / n) \rightarrow \bigoplus_{v \in \Sigma} \H^0(k_v, C \otimes^{\L} \Z / n) \xrightarrow{\gamma_1}  \H^0(U, \widehat{\mathcal{C}} \otimes^{\L} \Z / n)^D$, donc $(c_v) - (t_v)$ se rel\`eve en $\mu \in \H^0(U, C \otimes^{\L} \Z / n)$. Et finalement on a $(c_v) = (t_v) + \mu$. On conclut alors en remarquant que par finitude de $\bigoplus_{v \in \Sigma} \H^0(k_v, C \otimes^{\L} \Z / n)$, l'image de $\bigoplus_{v \in \Sigma} \H^0(k_v, C)^{\wedge}$ dans $\bigoplus_{v \in \Sigma} \H^0(k_v, C \otimes^{\L} \Z / n)$ co\"incide avec l'image de $\bigoplus_{v \in \Sigma} \H^0(k_v, C)$ dans $\bigoplus_{v \in \Sigma} \H^0(k_v, C \otimes^{\L} \Z / n)$, donc on peut bien \'ecrire $(c_v) = (t'_v) + \mu$ avec $(t'_v) \in \bigoplus_{v \in \Sigma} \H^0(k_v, C)$ et $\mu \in  \H^0(U, C \otimes^{\L} \Z / n)$. Cela ssure que $a_2 = 0$, ce qui termine la preuve du lemme \ref{lem erratum}.

\end{dem}

D\'efinissons aussi $$\D^0_{\wedge}(U, \mathcal{C}) := \textup{Ker} \left( \H^0(U, \mathcal{C}) \rightarrow \bigoplus_{v \in \Sigma} \H^0(k_v, \mathcal{C})^{\wedge} \right)$$

\begin{cor}
\label{coro D^i}
On a des accouplements fonctoriels $\D^0_{\wedge}(U, \mathcal{C}) \{l\} \times \D^{2}(U, \widehat{\mathcal{C}})\{l\} \rightarrow \Q / \Z$ et $\D^1(U, \mathcal{C}) \{l\} \times \D^{1}(U, \widehat{\mathcal{C}})\{l\} \rightarrow \Q / \Z$ (ainsi que $\D^2(U, \mathcal{C}) \{l\} \times \D^{0}(U, \widehat{\mathcal{C}})\{l\} \rightarrow \Q / \Z$), dont les noyaux \`a droite et \`a gauche sont les sous-groupes divisibles maximaux des deux groupes.
\end{cor}

\begin{dem}
Pour le premier accouplement, il suffit d'\'ecrire le diagramme
\begin{displaymath}
\xymatrix{
0 \ar[r] & \D^0_{\wedge}(U, \mathcal{C}) \{l\} \ar[r] \ar[d] & \H^0(U, \mathcal{C}) \{l\} \ar[r] \ar[d] & \bigoplus_{v \in \Sigma} \H^0(k_v, \mathcal{C})^{(l)} \ar[d] \\
0 \ar[r] & \overline{\D^2(U, \widehat{\mathcal{C}})}\{l\}^D \ar[r] & \overline{\H^2_c(U, \widehat{\mathcal{C}})}\{l\}^D \ar[r] & \bigoplus_{v \in \Sigma} \H^1(k_v, \widehat{\mathcal{C}})\{l\}^D
}
\end{displaymath}
et de remarquer que la fl\`eche verticale centrale est un isomorphisme puisque $ \H^0(U, \mathcal{C}) \{l\}$ est fini ($\H^0(U, \mathcal{C})$ est de type fini par le lemme \ref{lem struc etale}). Consid\'erons la troisi\`eme fl\`eche verticale : par le th\'eor\`eme de dualit\'e locale (th\'eor\`eme \ref{theo hensel}), le groupe $ \H^0(k_v, C)^{\wedge}$ s'identifie \`a $\H^1(k_v, \widehat{C})^D$, donc la troisi\`eme fl\`eche verticale du diagramme est injective, ce qui assure le r\'esultat. Pour le troisi\`eme accouplement, la preuve est similaire : on s'int\'eresse au diagramme exact suivant :
\begin{displaymath}
\xymatrix{
0 \ar[r] & \D^0(U, \widehat{\mathcal{C}}) \{l\} \ar[r] \ar[d] & \H^0(U, \widehat{\mathcal{C}}) \{l\} \ar[r] \ar[d] & \bigoplus_{v \in \Sigma} \H^0(k_v, \widehat{\mathcal{C}})\{l\} \ar[d] \\
0 \ar[r] & \overline{\D^2(U, \mathcal{C})}\{l\}^D \ar[r] & \overline{\H^2_c(U, \mathcal{C})}\{l\}^D \ar[r] & \bigoplus_{v \in \Sigma} \H^1(k_v, \mathcal{C})\{l\}^D
}
\end{displaymath}
Alors \`a nouveau la seconde fl\`eche verticale est un isomorphisme gr\^ace au lemme \ref{lem struc etale} ($\H^0(U, \widehat{\mathcal{C}})$ est de type fini), et la troisi\`eme fl\`eche verticale est injective (puisque $\H^0(k_v, \widehat{\mathcal{C}})$ s'injecte dans son compl\'et\'e $\H^0(k_v, \widehat{\mathcal{C}})^{\wedge}$, et ce dernier est isomorphe au dual du groupe de torsion $\H^1(k_v,\mathcal{C})$). 

Pour le second accouplement, le diagramme est le suivant :
\begin{displaymath}
\xymatrix{
0 \ar[r] & \D^1(U, \mathcal{C}) \{l\} \ar[r] \ar[d] & \H^1(U, \mathcal{C}) \{l\} \ar[r] \ar[d] & \bigoplus_{v \in \Sigma} \H^1(k_v, \mathcal{C})  \ar[d] \\
0 \ar[r] & \D^1(U, \widehat{\mathcal{C}})^D \ar[r] & \H^1_c(U, \widehat{\mathcal{C}})^D \ar[r] & \bigoplus_{v \in \Sigma} \H^0(k_v, \widehat{\mathcal{C}})^D
}
\end{displaymath}
On sait que la troisi\`eme fl\`eche verticale est injective, et que le noyau de la seconde est un sous-groupe divisible. Alors le lemme \ref{lem erratum} assure que le morphisme $\overline{\D^1(U, \mathcal{C})}\{ l \} \rightarrow \overline{\D^1(U, \widehat{\mathcal{C}})}^D$ est injectif.
Cela assure la seconde partie du corollaire \ref{coro D^i}, puisque l'argument "dual" est exactement le m\^eme : on dispose en effet des r\'esultats de finitude analogues pour $\widehat{C}$, qui permettent d'adapter la preuve du lemme \ref{lem erratum} avec les r\^oles de $\mathcal{C}$ et $\widehat{\mathcal{C}}$ \'echang\'es, pour obtenir que le morphisme $\overline{\D^1(U, \widehat{\mathcal{C}})}\{ l \} \rightarrow \overline{\D^1(U, \mathcal{C})}^D$ est lui aussi injectif, ce qui conclut la preuve de la dualit\'e parfaite entre $\overline{\D^1(U, \mathcal{C})} \{ l \}$ et $\overline{\D^1(U, \widehat{\mathcal{C}})}\{ l \}$.
\end{dem}

\begin{cor}
\label{coro D^1}
Avec les notations pr\'ec\'edentes, l'accouplement $\D^1(U,
\mathcal{C}) \{l\} \times \D^{1}(U, \widehat{\mathcal{C}})\{l\}
\rightarrow \Q / \Z$ est une dualit\'e parfaite de groupes finis,
fonctorielle en $C$.
\end{cor}

\begin{dem}
Il suffit de montrer que les groupes $\D^1(U, \mathcal{C}) \{l\}$ et $\D^1(U, \widehat{\mathcal{C}}) \{l\}$ sont finis. Pour $\D^1(U, \mathcal{C}) \{l\}$, par d\'evissage, il suffit de montrer que $D^1(U, \mathcal{T}_2)\{l\}$, $H^1(k_v, T_1)$ et $D^2(U, \mathcal{T}_1)$ sont finis. Ces trois points sont connus (voir \cite{Mil}, th\'eor\`eme II.4.6, corollaire I.2.3 et proposition II.4.14). En effet, consid\'erons le diagramme commutatif exact suivant :

\begin{displaymath}
\xymatrix{
H^1(U, \mathcal{T}_1)\{l\} \ar[r] \ar[d] & H^1(U, \mathcal{T}_2)\{l\} \ar[r] \ar[d] & \H^1(U, \mathcal{C})\{l\} \ar[r] \ar[d] & H^2(U, \mathcal{T}_1)\{l\} \ar[d] \\
\bigoplus_{v \in \Sigma} H^1(k_v, \mathcal{T}_1)\{l\} \ar[r] & \bigoplus_{v \in \Sigma} H^1(k_v, \mathcal{T}_2)\{l\} \ar[r] & \bigoplus_{v \in \Sigma} \H^1(k_v, \mathcal{C})\{l\} \ar[r] & \bigoplus_{v \in \Sigma} H^2(k_v, \mathcal{T}_1)\{l\}
}
\end{displaymath}
Par finitude de $\Sigma$ et de $H^1(k_v, T_1)$, le conoyau de la premi\`ere fl\`eche verticale est fini. Par finitude de $\D^1(U, \mathcal{T}_2)$, le noyau de la seconde est fini \'egalement. Enfin, la finitude de $D^2(U, \mathcal{T}_1)\{l\}$ assure la finitude du noyau de la quatri\`eme fl\`eche verticale $H^2(U, \mathcal{T}_1)\{l\} \rightarrow \bigoplus_{v \in \Sigma} H^2(k_v, T_1)\{l\}$. Par cons\'equent, une application imm\'ediate du lemme du serpent assure que le noyau de la troisi\`eme fl\`eche
$$\H^1(U, \mathcal{C})\{l\} \rightarrow \bigoplus_{v \in \Sigma} \H^1(k_v, \mathcal{C})\{l\}$$
 est fini, c'est-\`a-dire que $\D^1(U, \mathcal{C})\{l\}$ est fini.

En ce qui concerne la finitude de $\D^1(U, \widehat{\mathcal{C}})\{l\}$, le raisonnement est exactement le m\^eme, en utilisant les finitudes de $\Sigma$, $H^1(k_v, \widehat{\mathcal{T}_2})$, $D^2(U, \widehat{\mathcal{T}_2})\{l\}$ et $D^1(U, \widehat{\mathcal{T}_1})\{l\}$.
\end{dem}

\section{Dualit\'e globale : cohomologie galoisienne}
\label{section Galois}
L'objectif principal de cette section est de montrer les trois th\'eor\`emes de dualit\'e globale en cohomologie galoisienne, \`a savoir les th\'eor\`emes \ref{theo SH 1}, \ref{dualite SH 0} et \ref{theo SH 2}, que l'on utilisera dans la section \ref{section PT} pour obtenir la suite de Poitou-Tate.

\subsection{G\'en\'eralit\'es}

Dans cette section, $\Sigma$ d\'esigne un ensemble fini de places de $k$ contenant les places archim\'ediennes. On d\'efinit \'egalement $k_{\Sigma}$ comme \'etant l'extension maximale de $k$ non ramifi\'ee hors de $\Sigma$, et on note $\Gamma_{\Sigma} := \textup{Gal}(k_{\Sigma} | k)$.

\begin{defi}
Soit $C$ un complexe de tores sur $k$. On d\'efinit pour tout $i \geq 0$, 
$$\P^i(k, C) := \prod'_{v \in \Omega_k} \H^i(\widehat{k}_v, C)$$
o\`u le produit restreint est pris par rapport aux sous-groupes $\H^i_{\textup{nr}}(\widehat{k}_v, C)$, et
$$\textup{\cyr{SH}}^i(C) = \textup{\cyr{SH}}^i(k, C) := \textup{Ker} \left( \H^i(k, C) \rightarrow \P^i(k, C) \right)$$
On d\'efinit de m\^eme $\textup{\cyr{SH}}^i(\widehat{C})$.
Si $C$ est la fibre g\'en\'erique d'un complexe de tores $\mathcal{C}$ sur $U := \textup{Spec}(\mathcal{O}_{k, \Sigma})$, alors on note 
$$ \textup{\cyr{SH}}^i_{\Sigma}(C) = \textup{\cyr{SH}}^i_{\Sigma}(k, C) := \textup{Ker} \left( \H^i(\Gamma_{\Sigma}, C) \rightarrow \prod_{v \in \Sigma} \H^i(\widehat{k}_v, C) \right)$$
Enfin, on note
$$\textup{\cyr{SH}}^i_{\wedge}(C) := \textup{Ker} \left( \H^i(k, C)_{\wedge} \rightarrow \P^i(k, C)_{\wedge} \right)$$
\end{defi}

\begin{rem}
{\rm
On remarque que le groupe $\P^1(k, C)$ est en fait la somme directe $\P^1(k, C) = \bigoplus_v \H^1(\widehat{k}_v, C)$, puisque pour presque toute place $v$ $H^1(\widehat{\mathcal{O}}_v, \mathcal{T}_2) = H^2(\widehat{\mathcal{O}}_v, \mathcal{T}_1) = 0$.
}
\end{rem}

\begin{prop}
\label{prop etale-galois}
Soit $\mathcal{C}$ un complexe de tores d\'efini sur $U$, $C$ sa fibre g\'en\'erique, $l$ inversible sur $U$.
\begin{enumerate}
	\item Si $U$ est suffisamment petit, le morphisme canonique $\H^i(U, \mathcal{C})\{l\} \rightarrow \H^i(\Gamma_{\Sigma}, C)\{ l \}$ (resp. $\H^i(U, \widehat{\mathcal{C}})\{l\} \rightarrow \H^i(\Gamma_{\Sigma}, \widehat{C})\{ l \}$) est un isomorphisme, pour tout $i \geq 1$.
	\item Si $U$ est suffisamment petit, le morphisme canonique $\H^1(U, \mathcal{C})\{l\} \rightarrow \H^1(k, C)\{ l \}$ (resp. $\H^1(U, \widehat{\mathcal{C}})\{l\} \rightarrow \H^1(k, \widehat{C})\{ l \}$) est injectif.
\end{enumerate}
\end{prop}

\begin{dem}
\begin{enumerate}
	\item Gr\^ace aux corrigenda de \cite{HSz}, on conna\^it le r\'esultat pour les tores sur $U$ assez petit. On utilise alors la suite exacte :
$$H^i(U,\mathcal{T}_1) \rightarrow H^i(U, \mathcal{T}_2) \rightarrow \H^i(U, \mathcal{C}) \rightarrow H^{i+1}(U, \mathcal{T}_1) \rightarrow H^{i+1}(U, \mathcal{T}_2)$$
On sait que $ H^{j}(U, \mathcal{T}_r)\{ l \}  \cong  H^{j}(\Gamma_{\Sigma}, T_r)\{ l \}$, pour $j \geq 1$ (voir \cite{HSz}, proposition 4.1 et corrigenda), d'o\`u le r\'esultat par lemme des cinq.

En ce qui concerne $\widehat{\mathcal{C}}$, le raisonnement est similaire : en effet, $\widehat{\mathcal{T}_i}$ \'etant localement constant, on a un isomorphisme $H^1(U, \widehat{\mathcal{T}_i}) \xrightarrow{\simeq} H^1(\Gamma_{\Sigma}, \widehat{T_i})$. De m\^eme, par la proposition II.2.9 de \cite{Mil}, on a un isomorphisme $H^2(U, \widehat{\mathcal{T}_i}) \xrightarrow{\simeq} H^2(\Gamma_{\Sigma}, \widehat{T_i})$, d'o\`u le point 1. par d\'evissage.

	\item Il suffit de montrer que le morphisme naturel $\H^1(\Gamma_{\Sigma}, C) \rightarrow \H^1(k, C)$ (resp. $\H^1(\Gamma_{\Sigma}, \widehat{C}) \rightarrow \H^1(k, \widehat{C})$) est injectif. Pour le cas de $C$, on regarde le diagramme suivant :
\begin{displaymath}
\xymatrix{
H^1(\Gamma_{\Sigma},  T_1) \ar[r] \ar[d] & H^1(\Gamma_{\Sigma}, T_2) \ar[d] \ar[r] &  \H^1(\Gamma_{\Sigma}, C) \ar[r] \ar[d] & H^2(\Gamma_{\Sigma}, T_1) \ar[r] \ar[d] & H^2(\Gamma_{\Sigma}, T_2) \ar[d] \\
H^1(k, T_1) \ar[r] & H^1(k, T_2)  \ar[r] & \H^1(k, C) \ar[r] & H^2(k, T_1) \ar[r] & H^2(k, T_2)
}
\end{displaymath}

\'Ecrivons la suite exacte de restriction-inflation relative au quotient $\Gamma_{\Sigma}$ de $\Gamma_k$ :
$$0 \rightarrow H^1(\Gamma_{\Sigma}, T_i) \rightarrow H^1(k, T_i) \rightarrow H^1(k_{\Sigma}, T_i)^{\Gamma_{\Sigma}} \rightarrow H^2(\Gamma_{\Sigma}, T_i) \rightarrow H^2(k, T_i)$$
Quitte \`a augmenter $\Sigma$ (i.e. quitte \`a r\'eduire $U$), on peut supposer que $T_i$ est d\'eploy\'e par $k_{\Sigma}/k$, et donc par Hilbert 90, $H^1(k_{\Sigma}, T_i) = 0$. On en d\'eduit donc que les deux premi\`eres fl\`eches verticales du diagramme sont des isomorphismes, et que les deux derni\`eres sont injectives. Une chasse au diagramme assure alors que la fl\`eche verticale centrale est injective.

En ce qui concerne $\widehat{C}$, si $k_{\Sigma}$ d\'eploie $T_1$ et $T_2$, on a des isomorphisme $H^1(\Gamma_{\Sigma}, \widehat{T_i}) \xrightarrow{\simeq} H^1(k, \widehat{T_i})$ et des morphismes injectifs $H^2(\Gamma_{\Sigma}, \widehat{T_i}) \rightarrow H^2(k, \widehat{T_i})$ gr\^ace \`a la suite exacte de restriction inflation pour les $\widehat{T_i}$. Par d\'evissage, cela implique bien l'injectivit\'e de $\H^1(\Gamma_{\Sigma}, \widehat{C}) \rightarrow \H^1(k, \widehat{C})$ pour $\Sigma$ suffisamment gros. Cela conclut la preuve du point 2.
\end{enumerate}
\end{dem}

\begin{cor}
Soit $\mathcal{C}$ un complexe de tores d\'efini sur $U$, $C$ sa fibre g\'en\'erique, $l$ inversible sur $U$.
On a une dualit\'e parfaite fonctorielle de groupes finis $\textup{\cyr{SH}}^1_{\Sigma}(C)\{l\} \times \textup{\cyr{SH}}^1_{\Sigma}(\widehat{C})\{l\} \rightarrow \Q / \Z$.
\end{cor}

\begin{dem}
C'est une application directe du corollaire \ref{coro D^1}, en utilisant le point 1 de la proposition \ref{prop etale-galois}.
\end{dem}

\subsection{Dualit\'e entre $\textup{\cyr{SH}}^1(C)$ et $\textup{\cyr{SH}}^1(\widehat{C})$}

On va montrer dans cette section qu'il existe une dualit\'e parfaite de groupes finis entre $\textup{\cyr{SH}}^1(C)$ et $\textup{\cyr{SH}}^1(\widehat{C})$. On conserve les notations de la section pr\'ec\'edente : on fixe un nombre premier $l$, et on choisit un ouvert $U = \textup{Spec}(\mathcal{O}_{k, \Sigma})$ suffisamment petit pour que $l$ soit inversible sur $U$ et $C$ s'\'etende en un complexe de tores sur $U$.

\begin{lem}
Avec les notations pr\'ec\'edentes, le groupe $\D^1(U, \mathcal{C})\{ l \}$, vu comme sous-groupe de $\H^1(k, C)$, est contenu dans $\textup{\cyr{SH}}^1(C)\{ l \}$.
\end{lem}

\begin{dem}
La trivialit\'e des groupes $H^1(\mathcal{O}_v, T_i)$ et $H^2(\mathcal{O}_v, T_i)$ pour $v \notin S$ assure le r\'esultat ($i = 1$ ou $2$). La premi\`ere est une cons\'equence du th\'eor\`eme de Lang, et la seconde provient de l'isomorphisme $H^2(\mathcal{O}_v, T_i) \cong H^2(\F_v, T_i \times_{\mathcal{O}_v} \F_v)$, le second groupe \'etant trivial car $\F_v$ est de dimension cohomologique $1$ et $T_i(\overline{\F_v})$ est de torsion.
\end{dem}

\begin{lem}
Pour un ouvert $U$ suffisamment petit, le groupe $\D^1(U, \widehat{\mathcal{C}})\{ l \}$, vu comme sous-groupe de $\H^1(k, \widehat{C})$, est contenu dans $\textup{\cyr{SH}}^1(\widehat{C})\{ l \}$.
\end{lem}

\begin{dem}
C'est la m\^eme id\'ee que la preuve du lemme 4.7 de \cite{HSz}. En effet, on a l'inclusion suivante : 
$$ \bigcap_{V \subset U} \D^1(V, \widehat{\mathcal{C}}) \{ l \} \subseteq \textup{\cyr{SH}}^1(\widehat{C})\{ l \}$$
Or les groupes $\D^1(V, \widehat{\mathcal{C}}) \{ l \}$ sont finis, donc on peut se limiter \`a une intersection finie, i.e. \`a un nombre fini d'ouverts $V_i \subset U$. On consid\`ere alors un ouvert $U'$ contenu dans l'intersection des $V_i$, et la fonctorialit\'e covariante de la cohomologie \`a support compact assure alors que l'on a une inclusion $\D^1(U', \widehat{\mathcal{C}}) \{ l \} \subset \textup{\cyr{SH}}^1(\widehat{C})\{ l \}$, ce qui prouve le lemme.
\end{dem}

\begin{theo}
\label{theo SH 1}
Soit $C$ un complexe de tores sur $k$. Il existe une dualit\'e
parfaite de groupes finis, fonctorielle en $C$ :
$$\textup{\cyr{SH}}^1(C) \times \textup{\cyr{SH}}^1(\widehat{C}) \rightarrow \Q / \Z$$
\end{theo}

\begin{rem}
{\rm
Ce th\'eor\`eme est \'equivalent \`a la proposition 4.2 de
\cite{Bor3}, dont la preuve est due \`a Kottwitz. Par ailleurs,
Borovoi a donn\'e une autre preuve de ce r\'esultat r\'ecemment (voir
\cite{Bor4}), par d\'evissage au cas des groupes de type multiplicatif, \`a partir du
th\'eor\`eme I.4.12.(a) de \cite{Mil}.
}
\end{rem}

\begin{dem}
Soit $l$ un nombre premier. Il existe un ouvert $U$ de $\textup{Spec } \mathcal{O}_k$ tel que $C$ s'\'etend en un complexe de tores sur $U$ et $l$ est inversible sur $U$. On peut r\'eduire $U$ de sorte \`a \^etre dans les conditions pr\'ec\'edentes. Alors l'inclusion $\D^1(U, \mathcal{C})\{ l \} \subset \textup{\cyr{SH}}^1(C)\{ l \}$ est en fait une \'egalit\'e, puisque un \'el\'ement $\alpha \in \textup{\cyr{SH}}^1(C)\{ l \}$ s'\'etend en un \'el\'ement de $\D^1(V, \mathcal{C}) \{ l \}$, avec $V \subset U$, et donc $\D^1(V, \mathcal{C}) \{ l \} \subset \D^1(U, \mathcal{C}) \{ l \}$. De la m\^eme fa\c con, gr\^ace aux r\'esultats pr\'ec\'edents, on a \'egalement une \'egalit\'e $\D^1(U, \widehat{\mathcal{C}})\{ l \} = \textup{\cyr{SH}}^1(\widehat{C})\{ l \}$. Ces deux \'egalit\'es permettent bien de d\'efinir l'accouplement, et son exactitude d\'ecoule du corollaire \ref{coro D^1}. Pour la finitude des groupes $\textup{\cyr{SH}}^1(C)$ et $\textup{\cyr{SH}}^1(\widehat{C})$, c'est \'egalement une cons\'equence du corollaire \ref{coro D^1}, en remarquant que la preuve de ce corollaire montre en fait que les \'el\'ements de $\D^1(U, \mathcal{C})$ dont la torsion est inversible sur $U$ sont en nombre fini.
\end{dem}

\subsection{Dualit\'e entre $\textup{\cyr{SH}}^2(C)$ et $\textup{\cyr{SH}}^0(\widehat{C})$}

Dans cette section, on montre que les groupes  $\textup{\cyr{SH}}^2(C)$ et $\textup{\cyr{SH}}^0(\widehat{C})$ sont finis et duaux l'un de l'autre.

\begin{lem}
\label{lem SH 0 fini}
Soit $C$ un complexe de tores sur $k$. Alors $\textup{\cyr{SH}}^0(\widehat{C})$ est fini.
\end{lem}

\begin{dem}
\begin{itemize}
	\item On suppose d'abord que  $T_2$ est d\'eploy\'e sur $k$ : alors $\Gamma_k$ agit trivialement sur $\widehat{T_2}$, et on consid\`ere le diagramme commutatif suivant, pour une place $v$ quelconque :
\begin{displaymath}
\xymatrix{
H^0(k,  \widehat{T_2}) \ar[r] \ar[d] & H^0(k, \widehat{T_1}) \ar[d] \ar[r] &  \H^0(k, \widehat{C}) \ar[r] \ar[d] & H^1(k, \widehat{T_2}) \ar[d] \\
H^0(\widehat{k}_v, \widehat{T_2}) \ar[r] & H^0(\widehat{k}_v, \widehat{T_1})  \ar[r] & \H^0(\widehat{k}_v, \widehat{C}) \ar[r] & H^1(\widehat{k}_v, \widehat{T_2})
}
\end{displaymath}
La deuxi\`eme fl\`eche verticale est injective, le groupe $H^1(k, \widehat{T_2})$ est fini, et enfin la fl\`eche verticale de gauche est un isomorphisme, ce qui assure le r\'esultat dans ce cas.
	\item Cas g\'en\'eral : on ne suppose plus $T_2$ d\'eploy\'e sur $k$. Il existe une extension galoisienne finie $L/k$ qui d\'eploie $T_2$. Par le cas pr\'ec\'edent, on sait que $\textup{\cyr{SH}}^0(L, \widehat{C})$ est fini. Par un argument de restriction-corestriction, le groupe $\textup{\cyr{SH}}^0(k, \widehat{C})$ est donc de torsion. Or ce groupe est de type fini, puisque $\H^0(k, \widehat{C})$ l'est, donc $\textup{\cyr{SH}}^0(k, \widehat{C})$ est fini, ce qui conclut le preuve.
\end{itemize}
\end{dem}

Notons d\'esormais $\textup{\cyr{SH}}^0_{\wedge}(\widehat{C})$ le noyau du morphisme $\H^0(k, \widehat{C})_{\wedge} \rightarrow \P^0(k, \widehat{C})_{\wedge}$.

\begin{prop}
\label{prop 2}
Il existe une dualit\'e parfaite et fonctorielle
$$\textup{\cyr{SH}}^2(C) \times  \textup{\cyr{SH}}_{\wedge}^0(\widehat{C}) \rightarrow \Q / \Z$$
\end{prop}

\begin{dem}
On montre ce r\'esultat en  \'etapes : on se donne $n$ inversible sur $U$.
\begin{itemize}
  \item Par d\'evissage, on montre facilement que le morphisme $\H^0(U, \widehat{\mathcal{C}} \otimes^{\L} \Z / n) \rightarrow 
\H^0(k, \widehat{\mathcal{C}} \otimes^{\L} \Z / n)$ est injectif. En effet, on a un diagramme commutatif exact de la forme
\begin{displaymath}
\xymatrix{
0 \ar[r] & H^1(U, T_{\Z / n}(\widehat{\mathcal{C}})) \ar[r] \ar[d] & \H^0(U, \widehat{\mathcal{C}} \otimes^{\L} \Z / n) \ar[r] \ar[d] & H^0(U,\widehat{ _n \textup{Ker } \rho}) \ar[d] \\
0 \ar[r] & H^1(k, T_{\Z / n}(\widehat{C})) \ar[r] & \H^0(k, \widehat{C} \otimes^{\L} \Z / n) \ar[r] & H^0(k, \widehat{ _n \textup{Ker } \rho})
}
\end{displaymath}
Or par finitude de $T_{\Z / n}(\widehat{\mathcal{C}})$ et $_n \textup{Ker } \rho$, les premi\`ere et derni\`ere fl\`eches verticales sont injectives, donc la fl\`eche centrale \'egalement.

\item  On montre ensuite que pour tout place $v$ de $U$, $\H^1(\mathcal{O}_v,\mathcal{C} \otimes^{\L} \Z / n) = 0$. En effet, ce groupe s'int\`egre dans la suite exacte suivante :
$$H^3(\mathcal{O}_v, { _n \textup{Ker } \rho}) \rightarrow \H^1(\mathcal{O}_v,\mathcal{C} \otimes^{\L} \Z / n) \rightarrow H^2(\mathcal{O}_v,T_{\Z / n}(\mathcal{C}))$$
Or les deux groupes extr\^emes dans cette suite sont nuls par dimension cohomologique, donc le groupe central est nul.

\item La proposition \ref{dualite finie} et un d\'evissage \`a partir de la dualit\'e locale pour les modules finis assurent que l'on a une dualit\'e parfaite entre les groupes $\D^1(U, \mathcal{C} \otimes^{\L} \Z / n)$ et $\D^0(U,\widehat{\mathcal{C}} \otimes^{\L} \Z / n)$.

\item On fixe alors un nombre premier $l$ et on prend un ouvert $U$ suffisamment petit pour que $l$ soit inversible sur $U$ et pour que le complexe de tores s'\'etende sur $U$. Par le premier point, on sait que le groupe $\textup{\cyr{SH}}^0(\widehat{C} \otimes \Z / l^m)$ contient l'intersection des groupes $\D^0(V, \widehat{\mathcal{C}} \otimes^{\L} \Z / l^m)$ ($V$ d\'ecrivant les ouverts non vide de $U$) dans $\H^0(k, \widehat{\mathcal{C}} \otimes^{\L} \Z / l^m)$. Or chacun des groupes $\D^0(V, \widehat{\mathcal{C}} \otimes^{\L} \Z / l^m)$ \'etant finis, on peut se limiter \`a une intersection sur un nombre fini d'ouverts $V$, et en prenant leur intersections, on d\'eduit par fonctorialit\'e covariant de la cohomologie \`a support compact l'\'egalit\'e $\D^0(V, \widehat{\mathcal{C}} \otimes^{\L} \Z / l^m) \cong \textup{\cyr{SH}}^0(\widehat{C} \otimes \Z / l^m)$ pour tout ouvert $V$ de $U$ contenu dans un certain ouvert $U_0$ (voir par exemple \cite{HSz}, preuves du lemme 4.7 et de la proposition 4.12).

Le second point assure que pour une immersion ouverte $V \rightarrow U$, le morphisme naturel $\H^1(U, \mathcal{C} \otimes^{\L} \Z / l^m) \rightarrow \H^1(V, \mathcal{C} \otimes^{\L} \Z / l^m)$ envoie le sous-groupe $\D^1(U, \mathcal{C} \otimes^{\L} \Z / l^m)$ dans $\D^1(V, \mathcal{C} \otimes^{\L} \Z / l^m)$. On en d\'eduit donc un morphisme naturel $\varinjlim_V \D^1(V, \mathcal{C} \otimes^{\L} \Z / l^m) \rightarrow \textup{\cyr{SH}}^1(C \otimes \Z / l^m)$, qui est surjectif par fonctorialit\'e covariante de la cohomologie \`a support compact (voir \cite{HSz}, preuve de la proposition 4.12). Et ce morphisme est injectif car le morphisme $\varinjlim_V \H^1(V, \mathcal{C} \otimes^{\L} \Z / l^m) \rightarrow \H^1(C \otimes \Z / l^m)$ l'est. 

Finalement, on a donc des isomorphismes naturels 
$$\textup{\cyr{SH}}^1(C \otimes \Z / l^m) \cong \varinjlim_V \D^1(V, \mathcal{C} \otimes^{\L} \Z / l^m)$$
et
$$\textup{\cyr{SH}}^0(\widehat{C} \otimes \Z / l^m) \cong \varprojlim_V \D^0(V, \widehat{\mathcal{C}} \otimes^{\L} \Z / l^m)$$

\item On combine alors les points 3 et 4 pour en d\'eduire que l'on a une dualit\'e parfaite entre les groupes $\textup{\cyr{SH}}^1(C \otimes \Z / l^m)$ et $\textup{\cyr{SH}}^0(\widehat{C} \otimes \Z / l^m)$. Ceci \'etant valable pour tout nombre premier $l$, on en d\'eduit une dualit\'e parfaite entre les groupes $\varinjlim_n \textup{\cyr{SH}}^1(C \otimes^{\L} \Z / n)$ et $\varprojlim_n \textup{\cyr{SH}}^0(\widehat{C} \otimes^{\L} \Z / n)$

On identifie alors les groupes $\textup{\cyr{SH}}^0_{\wedge}(\widehat{C}) \cong \varprojlim_n \textup{\cyr{SH}}^0(\widehat{C} \otimes^{\L} \Z / n)$ en passant \`a la limite projective dans le diagramme
\begin{displaymath}
\xymatrix{
0 \ar[r] & \H^0(k, \widehat{C}) / n \ar[d] \ar[r] & \H^0(k, \widehat{C} \otimes^{\L} \Z / n) \ar[r] \ar[d] & {_n \H^1(k, \widehat{C})} \ar[d] \ar[r] & 0 \\
0 \ar[r] & \P^0(k, \widehat{C}) / n \ar[r] & \P^0(k, \widehat{C} \otimes^{\L} \Z / n) \ar[r] & {_n \P^1(k, \widehat{C})} \ar[r] & 0 
}
\end{displaymath}
et en utilisant la finitude du groupe $\textup{\cyr{SH}}^1(\widehat{C})$.
De m\^eme, on identifie les groupes $\textup{\cyr{SH}}^2(C) \cong \varinjlim_n \textup{\cyr{SH}}^1(C \otimes^{\L} \Z / n)$, en passant \`a la limite inductive dans le diagramme
\begin{displaymath}
\xymatrix{
0 \ar[r] & \H^1(k, C) / n \ar[d] \ar[r] & \H^1(k, C \otimes^{\L} \Z / n) \ar[r] \ar[d] & {_n \H^2(k, C)} \ar[d] \ar[r] & 0 \\
0 \ar[r] & \P^1(k, C) / n \ar[r] & \P^1(k, C \otimes^{\L} \Z / n) \ar[r] & {_n \P^2(k, C)} \ar[r] & 0 
}
\end{displaymath}
qui permet d'identifier $\varinjlim_n \H^1(k, C \otimes^{\L} \Z / n) \cong \H^2(k, C)$ et $\varinjlim_n \P^1(k, C \otimes^{\L} \Z / n) \cong \P^2(k, C)$ puisque les groupes $\H^1(k, C)$, $\P^1(k, C)$, $\H^2(k, C)$ et $\P^2(k, C)$ sont de torsion, et on conclut gr\^ace \`a l'exactitude du foncteur $\varinjlim$ qui permet d'identifier $\varinjlim_n \textup{\cyr{SH}}^1(C \otimes^{\L} \Z / n)$ au noyau de $\varinjlim_n \H^1(k, C \otimes^{\L} \Z / n) \rightarrow \varinjlim_n \P^1(k, C \otimes^{\L} \Z / n)$.
\end{itemize}
\end{dem}

\begin{prop}
\label{prop 3}
$$\textup{\cyr{SH}}^0_{\wedge}(\widehat{C}) \cong \textup{Ker}(\H^0(k, \widehat{C})^{\wedge} \xrightarrow{\beta^0} \P^0(k, \widehat{C})^{\wedge})$$
\end{prop}

\begin{dem}
$\H^0(k, \widehat{C})$ est discret de type fini, donc le morphisme $\H^0(k, \widehat{C})_{\wedge} \rightarrow \H^0(k, \widehat{C})^{\wedge}$ est un isomorphisme. 
De plus, pour tout $n$, et pour toute place $v$ ne divisant pas $n$, le morphisme $\H^0(\mathcal{O}_v, \widehat{\mathcal{C}} \otimes^{\L} \Z / n) \rightarrow \H^0(k_v, \widehat{C} \otimes^{\L} \Z / n)$ est injectif (car les morphismes $H^0(\mathcal{O}_v, \widehat{ _n \textup{Ker } \rho}) \rightarrow H^0(k_v, \widehat{ _n \textup{Ker } \rho})$ et $H^1(\mathcal{O}_v, T_{\Z / n}(\widehat{\mathcal{C}})) \rightarrow H^1(k_v, T_{\Z / n}(\widehat{C}))$ le sont). Donc pour tout $n$ le morphisme $\P^0(k, \widehat{C}) / n \rightarrow \prod_v \left( \H^0(k_v, \widehat{C}) / n \right)$ est injectif, ce qui implique que le morphisme $\P^0(k, \widehat{C})_{\wedge} \rightarrow \P^0(k, \widehat{C})^{\wedge}$ est injectif (voir preuve de la proposition 5.4 dans \cite{HSz}).
On a donc un diagramme commutatif
\begin{displaymath}
\xymatrix{
\H^0(k, \widehat{C})_{\wedge} \ar[r]^{=} \ar[d] & \H^0(k, \widehat{C})^{\wedge} \ar[d] \\
\P^0(k, \widehat{C})_{\wedge} \ar[r] & \P^0(k, \widehat{C})^{\wedge}
}
\end{displaymath}
o\`u la fl\`eche horizontale inf\'erieure est injective, ce qui permet bien de montrer que les noyaux des deux fl\`eches verticales sont les m\^emes.

\end{dem}


\begin{theo}
\label{theo SH 2} On suppose $\textup{Ker } \rho$ fini. Alors il
existe une dualit\'e parfaite de groupes finis, fonctorielle en $C$ :
$$\textup{\cyr{SH}}^2(C) \times \textup{\cyr{SH}}^0(\widehat{C}) \rightarrow \Q / \Z$$
\end{theo}

\begin{dem}
On sait que $\H^0(k, \widehat{C})$ est un groupe discret de type fini, donc $\H^0(k, \widehat{C}) \rightarrow \H^0(k, \widehat{C})^{\wedge}$ est injective. Donc on a une suite exacte (puisque $\textup{\cyr{SH}}^0(\widehat{C})$ est un sous-groupe ferm\'e de $\H^0(k, \widehat{C})$) :
\begin{equation}
\label{sec Im theta}
0 \rightarrow \textup{\cyr{SH}}^0(\widehat{C})^{\wedge} \rightarrow \H^0(k, \widehat{C})^{\wedge} \rightarrow \textup{Im}(\theta)^{\wedge}
\end{equation}
o\`u $\theta$ est le morphisme $\H^0(k, \widehat{C}) \rightarrow \P^0(k, \widehat{C})$, et $\textup{Im}(\theta)$ est muni de la topologie quotient de $\H^0(k, \widehat{C})$, c'est-\`a-dire de la topologie discr\`ete.
On consid\`ere alors le diagramme commutatif suivant :
\begin{displaymath}
\xymatrix{
\P^0(k, \widehat{C}) \ar[r] \ar[d] & \P^0(k, \widehat{C})_{\wedge} \ar[d] \\
\prod_v \H^0(k_v, \widehat{C}) \ar[r] & \prod_v \H^0(k_v, \widehat{C})_{\wedge} 
}
\end{displaymath}
La premi\`ere fl\`eche verticale est injective, ainsi que la fl\`eche horizontale inf\'erieure ($\H^0(k_v, \widehat{C})$ est de type fini), donc la fl\`eche horizontale sup\'erieure est injective. Donc le morphisme $\P^0(k, \widehat{C}) \rightarrow \P^0(k, \widehat{C})^{\wedge}$ est injectif. En particulier $\textup{Im}(\theta)$ s'injecte dans $\P^0(k, \widehat{C})^{\wedge}$. On montre alors le lemme suivant :

\begin{lem}
On suppose $\textup{Ker } \rho$ fini. Alors on a une suite exacte naturelle
$$\H^0(k, \widehat{C}) \xrightarrow{\theta} \P^0(k, \widehat{C}) \rightarrow \H^1(k, C)^D$$
\end{lem}

\begin{dem}
Tout d'abord, sous l'hypoth\`ese de finitude de $\textup{Ker } \rho$, il est clair par d\'evissage (groupe fini et tore) que $\textup{\cyr{SH}}^2(C)$ est fini.
On consid\`ere le complexe
$$\varinjlim_n \H^{-1}(k, \widehat{C} \otimes^{\L} \Z / n) \rightarrow \varinjlim_n \P^{-1}(k, \widehat{C} \otimes^{\L} \Z / n) \rightarrow \varinjlim_n \H^{1}(k, C \otimes^{\L} \Z / n)^D$$
Celui-ci est exact en utilisant le triangle exact suivant
$$T_{\Z / n}(\widehat{C})[1] \rightarrow \widehat{C} \otimes^{\L} \Z / n \rightarrow \widehat{ _n \textup{Ker } \rho } \rightarrow T_{\Z / n}(\widehat{C})[2]$$
et la suite exacte de Poitou-Tate pour le module fini $T_{\Z / n}(\widehat{C})$, puis en passant \`a la limite inductive.
On consid\`ere alors le diagramme commutatif \`a lignes exactes suivant :
\begin{displaymath}
\xymatrix{
0 \ar[r] & \varinjlim_n \H^{-1}(k, \widehat{C}) / n \ar[r] \ar[d] & \varinjlim_n \H^{-1}(k, \widehat{C} \otimes^{\L} \Z / n) \ar[r] \ar[d] & \H^0(k, \widehat{C})_{\textup{tors}} \ar[d] \ar[r] & 0 \\
0 \ar[r] & \varinjlim_n \P^{-1}(k, \widehat{C}) / n \ar[r] \ar[d] & \varinjlim_n \P^{-1}(k, \widehat{C} \otimes^{\L} \Z / n) \ar[r] \ar[d] & \P^0(k, \widehat{C})_{\textup{tors}} \ar[d] \ar[r] & 0 \\
0 \ar[r] & \varinjlim_n (\H^0(k, C)^D) / n \ar[r] & \varprojlim_n \H^{1}(k, C \otimes^{\L} \Z / n)^D \ar[r] & (\H^1(k, C)^D)_{\textup{tors}} \ar[r] & 0
}
\end{displaymath}
dont on a montr\'e que la colonne centrale \'etait exacte. On conclut alors que la troisi\`eme colonne est exacte par surjectivit\'e du morphisme 
$$ \varinjlim_n \P^{-1}(k, \widehat{C}) / n \rightarrow \varinjlim_n {(\H^0(k, C)^D) / n}$$
dont le conoyau est $\textup{\cyr{SH}}^2(C) \otimes \Q / \Z = 0$ puisque $\textup{\cyr{SH}}^2(C)$ est fini.
\end{dem}

Ce lemme assure donc que $\textup{Im}(\theta)$ est un sous-groupe ferm\'e de $\P^0(k, \widehat{C})$, donc les arguments pr\'ec\'edents impliquent que le morphisme $\textup{Im}(\theta)^{\wedge} \rightarrow \P^0(k, \widehat{C})^{\wedge}$ est injectif.
En revenant alors \`a la suite exacte (\ref{sec Im theta}), $\textup{\cyr{SH}}^0(\widehat{C})^{\wedge}$ s'identifie au noyau $\textup{Ker}(\H^0(k, \widehat{C})^{\wedge} \xrightarrow{\beta^0} \P^0(k, \widehat{C})^{\wedge})$, et donc aussi \`a $\textup{\cyr{SH}}^0_{\wedge}(\widehat{C})$ par la proposition \ref{prop 3}. On utilise alors la finitude de $\textup{\cyr{SH}}^0(\widehat{C})$ (lemme \ref{lem SH 0 fini}) et le r\'esultat de dualit\'e globale (proposition \ref{prop 2}) pour conclure.
\end{dem}

\subsection{Dualit\'e entre $\textup{\cyr{SH}}^0_{\wedge}(C)$ et $\textup{\cyr{SH}}^2(\widehat{C})$}

Dans cette partie, les r\'esultats de dualit\'e ne seront pas valable en g\'en\'eral, mais seulement pour certains complexes de tores particuliers. On est amen\'e \`a imposer deux types de conditions sur $C = [T_1 \xrightarrow{\rho} T_2]$ : 
\begin{enumerate}
  \item $\textup{Ker } \rho$ est fini. Cette hypoth\`ese sera v\'erifi\'ee dans le cadre de la cohomologie ab\'elianis\'ee des groupes r\'eductifs. On voit facilement que la condition de finitude de $\textup{Ker } \rho$ est \'equivalente \`a la condition $\varprojlim_n  { _n (\textup{Ker } \rho)} = 0$.
  \item $\rho$ est surjective. Cette hypoth\`ese implique que $C$ est quasi-isomorphe \`a un objet $M[1]$, o\`u $M := \textup{Ker } \rho$ est un $k$-groupe de type multiplicatif. R\'eciproquement, tout $k$-groupe de type multiplicatif peut s'\'ecrire comme le noyau d'un morphisme de tores surjectif. Par ailleurs, on constate facilement que la condition de surjectivit\'e de $\rho$ \'equivaut \`a la condition $\varinjlim_n T_{\Z / n}(C) = 0$.
\end{enumerate}

\subsubsection{Cas o\`u $\textup{Ker}(\rho)$ est fini}

Dans cette section, sous l'hypoth\`ese de finitude de $\textup{Ker}(\rho)$, on montre l'existence d'une dualit\'e parfaite entre les groupes $\textup{\cyr{SH}}^0_{\wedge}(C)$ et $\textup{\cyr{SH}}^2(\widehat{C})$, puis on identifie $\textup{\cyr{SH}}^0_{\wedge}(C)$ avec le noyau de $\H^0(k, C)^{\wedge} \rightarrow \P^0(k, C)^{\wedge}$.

\begin{theo}
\label{dualite SH 0}
Supposons que $\textup{Ker}(\rho)$ est fini. Alors il existe une
dualit\'e parfaite de groupes finis, fonctorielle en $C$ :
$$\textup{\cyr{SH}}^0_{\wedge}(C) \times \textup{\cyr{SH}}^2(\widehat{C}) \rightarrow \Q / \Z$$
\end{theo}

\begin{dem}
On commence par montrer la finitude du groupe
$\textup{\cyr{SH}}^2(\widehat{C})$. Pour cela, on remarque que l'on a
un triangle exact 
$$M [1] \rightarrow C \rightarrow T \rightarrow M [2]$$
o\`u $M := \textup{Ker } \rho$ est un $k$-groupe de type multiplicatif
et $T := \textup{Coker } \rho$ un $k$-tore. On obtient alors un
diagramme commutatif de suites exactes :
\begin{displaymath}
\xymatrix{
H^3(k, \widehat{T}) \ar[r] \ar[d] & \H^2(k, \widehat{C}) \ar[r] \ar[d]
& H^2(k, \widehat{M}) \ar[r] \ar[d] & H^4(k, \widehat{T}) \ar[d] \\
P^3(k, \widehat{T}) \ar[r] & \P^2(k, \widehat{C}) \ar[r] & P^2(k, \widehat{M}) \ar[r] & P^4(k, \widehat{T})
}
\end{displaymath}
On sait que les premi\`ere et derni\`ere fl\`eches verticales sont
des isomorphismes de groupes finis. La finitude de
$\textup{\cyr{SH}}^2(\widehat{C})$ r\'esulte alors de la finitude de
$\textup{\cyr{SH}}^2(\widehat{M})$ (voir th\'eor\`eme 8.6.7 de \cite{NSW} ou th\'eor\`eme \ref{dualite
  SH 0 tm}) et du lemme du serpent appliqu\'e au diagramme pr\'ec\'edent.

Poursuivons alors la preuve du th\'eor\`eme \ref{dualite SH 0} avec le lemme suivant :

\begin{lem}
\label{lem technique}
Pour tout $n > 0$, le morphisme $\P^0(C) / n \rightarrow \prod_v \H^0(\widehat{k}_v, C) / n$ est injectif, d'image le produit restreint des $\H^0(\widehat{k}_v, C) / n$ par rapport aux $\H^0_{\textup{nr}}(\widehat{k}_v, C) / n$.
\end{lem}

\begin{dem}
Il suffit de montrer que, pour presque toute place $v$, dans le diagramme suivant
\begin{displaymath}
\xymatrix{
\H^0(\widehat{\mathcal{O}}_v, C) \ar[d] \ar[r]^{n} & \H^0(\widehat{\mathcal{O}}_v, C) \ar[r] \ar[d] & \H^0(\widehat{\mathcal{O}}_v, C \otimes^{\L} \Z / n) \ar[d] \\
\H^0(\widehat{k}_v, C) \ar[r]^{n} & \H^0(\widehat{k}_v, C) \ar[r] & \H^0(\widehat{k}_v, C \otimes^{\L} \Z / n)
}
\end{displaymath}
la troisi\`eme fl\`eche verticale est injective. Pour ce faire, on utilise la d\'ecomposition de $C \otimes^{\L} \Z / n$ d\'ej\`a utilis\'ee plus haut :
on a un diagramme commutatif \`a lignes exactes :
\begin{displaymath}
\xymatrix{
H^2(\widehat{\mathcal{O}}_v, {_n (\textup{Ker} \rho)}) \ar[d] \ar[r] & \H^0(\widehat{\mathcal{O}}_v, C \otimes^{\L} \Z / n) \ar[r] \ar[d] & H^1(\widehat{\mathcal{O}}_v, T_{ \Z / n}(C)) \ar[d] \\
H^2(\widehat{k}_v, {_n (\textup{Ker} \rho)}) \ar[r] & \H^0(\widehat{k}_v, C \otimes^{\L} \Z / n ) \ar[r] & H^1(\widehat{k}_v, T_{ \Z / n}(C))
}
\end{displaymath}

On suppose que $v$ ne divise pas $n$. Alors $H^2(\widehat{\mathcal{O}}_v, {_n (\textup{Ker} \rho)}) = 0$ (puisque $H^2(\widehat{\mathcal{O}}_v, {_n (\textup{Ker} \rho)}) \cong H^2(\F_v,  _n (\textup{Ker} \rho))$ car $v$ ne divise pas $n$, et $\F_v$ est de dimension cohomologique $1$), et la troisi\`eme fl\`eche verticale est injective, puisque le groupe $H^1_v(\widehat{\mathcal{O}}_v, T_{ \Z / n}(C))$ est le dual de $H^2(\widehat{\mathcal{O}}_v, \widehat{T_{ \Z / n}(C)})$ qui est nul par dimension cohomologique. Cela assure bien l'injectivit\'e de la seconde fl\`eche verticale pour presque toute place $v$, et donc le lemme.
\end{dem}

Poursuivons la preuve du th\'eor\`eme \ref{dualite SH 0}. 
Pour toute place $v$, la suite suivante
$$0 \rightarrow \H^0(\widehat{k}_v, C) / n \rightarrow \H^0(\widehat{k}_v, C \otimes^{\L} \Z / n) \rightarrow {_n \H^1(\widehat{k}_v, C)} \rightarrow 0$$
est exacte.
De m\^eme, pour toute place $v \notin \Sigma$, la suite
$$0 \rightarrow \H^0(\widehat{\mathcal{O}}_v, \mathcal{C}) / n \rightarrow \H^0(\widehat{\mathcal{O}}_v, \mathcal{C} \otimes^{\L} \Z / n) \rightarrow {_n \H^1(\widehat{\mathcal{O}}_v, \mathcal{C})} \rightarrow 0$$
est exacte.
Consid\'erons alors le complexe suivant :
$$\P^0(k, C) / n \rightarrow \P^0(k, C \otimes^{\L} \Z / n) \rightarrow {_n \P^1(k, C)}$$
L'exactitude des deux suites pr\'ec\'edentes assure la surjectivit\'e de la seconde fl\`eche. La premi\`ere partie du lemme \ref{lem technique} assure que la premi\`ere fl\`eche est injective, et enfin la seconde partie du lemme \ref{lem technique} assure l'exactitude du complexe pr\'ec\'edent en $ \P^0(k, C \otimes^{\L} \Z / n)$. Par cons\'equent, on dispose d'un diagramme commutatif \`a lignes exactes :

\begin{equation}
\label{diag niv fini}
\xymatrix{
0 \ar[r] & \H^0(k, C) / n \ar[d] \ar[r] & \H^0(k, C \otimes^{\L} \Z / n) \ar[r] \ar[d] & {_n \H^1(k, C)} \ar[d] \ar[r] & 0 \\
0 \ar[r] & \P^0(k, C) / n \ar[r] & \P^0(k, C \otimes^{\L} \Z / n) \ar[r] & {_n \P^1(k, C)} \ar[r] & 0 
}
\end{equation}

On passe alors \`a la limite projective sur $n$ : on obtient 
\begin{eqnarray}
\label{diag lim}
\xymatrix{
0 \ar[r] & \H^0(k, C)_{\wedge} \ar[d]^{\theta_0} \ar[r] & \varprojlim_n \H^0(k, C \otimes^{\L} \Z / n) \ar[r] \ar[d]^{\theta} & Q_1 \ar[d]^{\beta} \ar[r] & 0 \\
0 \ar[r] & \P^0(k, C)_{\wedge} \ar[r] & \varprojlim_n \P^0(k, C \otimes^{\L} \Z / n) \ar[r] & Q_2 \ar[r] & 0 
}
\end{eqnarray}
Or $Q_1$ est un sous-groupe du module de Tate $T(\H^1(k, C))$, et donc $\textup{Ker } \beta$ est contenu dans $T(\textup{\cyr{SH}}^1(C))$. Or $\textup{\cyr{SH}}^1(C)$ est fini (voir th\'eor\`eme \ref{theo SH 1}), donc $\beta$ est injective. Donc $\textup{\cyr{SH}}^0_{\wedge}(C) = \textup{Ker } \theta_0$ est isomorphe \`a $\textup{Ker } \theta$.

\begin{prop}
\label{prop dualite noyau}
Supposons que $\textup{Ker}(\rho)$ est fini. Alors il existe une dualit\'e parfaite $$\textup{Ker } \theta \times \textup{\cyr{SH}}^1\left( \varinjlim_n \widehat{C} \otimes^{\L} \Z / n \right) \rightarrow \Q / \Z$$
\end{prop}

\begin{dem}
On fait l'hypoth\`ese que $\textup{Ker}(\rho)$ est fini. On fixe un nombre premier $l$, et $U$ un ouvert sur lequel $C$ s'\'etend, et de sorte que $l$ soit inversible sur $U$.

\begin{lem} 
\label{lem AV D}
On a une dualit\'e parfaite fonctorielle en $C$ : $\varprojlim_n \D^0(U, \mathcal{C} \otimes^{\L} \Z / {l^n}) \times \D^1(U, \varinjlim_n \widehat{\mathcal{C}} \otimes^{\L} \Z / {l^n}) \rightarrow \Q / \Z$. 
\end{lem}

\begin{dem}
En effet, on dispose du diagramme suivant \`a lignes exactes :
\begin{displaymath}
\xymatrix{
0 \ar[r] & \varprojlim_n \D^0(U, \mathcal{C} \otimes^{\L} \Z / {l^n}) \ar[r] \ar[d] & \varprojlim_n \H^0(U, \mathcal{C} \otimes^{\L} \Z / {l^n}) \ar[r] \ar[d]^{\simeq} &\varprojlim_n \left( \bigoplus_{v \in \Sigma} \H^0(k_v, \mathcal{C} \otimes^{\L} \Z / {l^n}) \right) \ar[d]^{\simeq} \\
0 \ar[r] & \D^1(U, \varinjlim_n \widehat{\mathcal{C}}\otimes^{\L} \Z / {l^n})^D \ar[r] & (\varinjlim_n \H^1_c(U, \widehat{\mathcal{C}} \otimes^{\L} \Z / {l^n}))^D \ar[r] & (\varinjlim_n \bigoplus_{v \in \Sigma} \H^0(k_v, \widehat{\mathcal{C}} \otimes^{\L} \Z / {l^n}))^D
}
\end{displaymath}
l'isomorphisme de la colonne centrale provient de la dualit\'e globale en cohomologie \'etale (proposition \ref{dualite finie}). Concernant la troisi\`eme fl\`eche verticale, on va montrer que c'est un isomorphisme. Pour cela, on consid\`ere une place $v$ de $\Sigma$, et on s'int\'eresse au diagramme suivant, \`a lignes exactes 
\begin{changemargin}{-1.5cm}{1cm}
\begin{displaymath}
\xymatrix{
H^0(k_v, T_{\Z / n}(C)) \ar[r] \ar[d] & H^2(k_v, {_n \textup{Ker } \rho}) \ar[r] \ar[d] & \H^0(k_v, C \otimes^{\L} \Z / n) \ar[r] \ar[d] & H^1(k_v,  T_{\Z / n}(C)) \ar[r] \ar[d] & H^3(k_v, {_n \textup{Ker } \rho}) \ar[d] \\
H^2(k_v, T_{\Z / n}(\widehat{C}))^D \ar[r] & H^0(k_v, \widehat{_n \textup{Ker } \rho})^D \ar[r] & \H^0(k_v, \widehat{C} \otimes^{\L} \Z / n)^D \ar[r] & H^1(k_v,  T_{\Z / n}(\widehat{C}))^D \ar[r] & H^{-1}(k_v, \widehat{_n \textup{Ker } \rho}) \\
}
\end{displaymath}
\end{changemargin}

Alors par le th\'eor\`eme de dualit\'e locale pour un module fini sur un corps hens\'elien  (th\'eor\`eme I.2.14.(c) et th\'eor\`eme I.2.13.(a) de \cite{Mil}), les deux premi\`eres fl\`eches verticales et les deux derni\`eres sont des isomorphismes, donc la fl\`eche centrale \'egalement. Par cons\'equent, en revenant au premier diagramme, on a bien la dualit\'e annonc\'ee.
\end{dem}

\begin{lem}
\label{lem inj lim}
Le morphisme canonique $\varprojlim_n \H^0(U, \mathcal{C} \otimes^{\L} \Z / {l^n}) \rightarrow \varprojlim_n \H^0(k, C \otimes^{\L} \Z / {l^n})$ est injectif.
\end{lem}

\begin{dem}
La finitude de $\textup{Ker}(\rho)$ assure que $\varprojlim_n H^2(U, {_{l^n} \textup{Ker}(\rho)}) = 0$.  Par cons\'equent, un d\'evissage assure que les lignes du diagramme suivant sont exactes (les groupes \'etant finis, on peut passer \`a la limite projective) :
\begin{displaymath}
\xymatrix{
0 \ar[r] & \varprojlim_n \H^0(U,  \mathcal{C} \otimes^{\L} \Z / {l^n}) \ar[r] \ar[d] & \varprojlim_n H^1(U, T_{\Z / {l^n}}(\mathcal{C})) \ar[d] \\
& \varprojlim_n \H^0(k,  C \otimes^{\L} \Z / {l^n}) \ar[r] & \varprojlim_n H^1(k, T_{\Z / {l^n}}(C))
}
\end{displaymath}
Or la derni\`ere fl\`eche verticale est injective (voir \cite{Mil}, II.2.9), donc ce diagramme assure le r\'esultat.
\end{dem}

\begin{lem}
\label{lem nul lim}
Pour toute place $v$ de $k$ hors de $\Sigma$, $ \varinjlim_n \H^1(\mathcal{O}_v, \widehat{\mathcal{C}} \otimes^{\L} \Z / {l^n}) = 0$.
\end{lem}

\begin{dem}
On \'ecrit la suite exacte de cohomologie sur $\mathcal{O}_v$ :
$$H^2(\mathcal{O}_v, T_{\Z / l^n}(\widehat{C})) \rightarrow \H^1(\mathcal{O}_v, \widehat{C} \otimes \Z / {l^n}) \rightarrow H^1(\mathcal{O}_v, \widehat{{_{l^n} \textup{Ker}(\rho)}}) \rightarrow H^3(\mathcal{O}_v, T_{\Z / l^n}(\widehat{C}))$$
Or par dimension cohomologique, les groupes $H^i(\mathcal{O}_v, T_{\Z / l^n}(\widehat{C}))$ sont nuls pour $i \geq 2$, d'o\`u un isomorphisme 
$$\H^1(\mathcal{O}_v, \widehat{C} \otimes \Z / {l^n}) \cong H^1(\mathcal{O}_v, \widehat{ _{l^n} \textup{Ker}(\rho)})$$
On prend la limite inductive, et on remarque que la finitude de $\textup{Ker}(\rho)$ implique la nullit\'e de $\varinjlim_n \widehat{ _{l^n} \textup{Ker}(\rho)}$, d'o\`u $\varinjlim_n \H^1(\mathcal{O}_v, \widehat{C} \otimes \Z / {l^n}) = 0$.
\end{dem}

L'objectif est d\'esormais de passer \`a la limite sur $V$ dans le lemme \ref{lem AV D} pour en d\'eduire la proposition \ref{prop dualite noyau}. Avant de passer \`a la limite, d\'efinissons les morphismes de transition : si $V \rightarrow U$ est une immersion ouverte, la fonctorialit\'e covariante de la cohomologie \`a support compact induit un morphisme canonique $\varprojlim_n \D^0(V,  \mathcal{C} \otimes^{\L} \Z / {l^n}) \rightarrow \varprojlim_n  \D^0(U, \mathcal{C} \otimes^{\L} \Z / {l^n})$. En ce qui concerne $\widehat{\mathcal{C}}$, on consid\`ere le morphisme naturel $\varinjlim_n \H^1(U,  \widehat{\mathcal{C}} \otimes^{\L} \Z / {l^n}) \rightarrow \varinjlim_n \H^1(V,  \widehat{\mathcal{C}} \otimes^{\L} \Z / {l^n})$, et gr\^ace au lemme \ref{lem nul lim}, le sous-groupe $\D^1(U, \varinjlim_n \widehat{\mathcal{C}} \otimes^{\L} \Z / {l^n})$ de $\varinjlim_n \H^1(U,  \widehat{\mathcal{C}} \otimes^{\L} \Z / {l^n})$ s'envoie dans $\D^1(V, \varinjlim_n \widehat{\mathcal{C}} \otimes^{\L} \Z / {l^n})$, d'o\`u un morphisme canonique $\D^1(U, \varinjlim_n \widehat{\mathcal{C}} \otimes^{\L} \Z / {l^n}) \rightarrow \D^1(V, \varinjlim_n \widehat{\mathcal{C}} \otimes^{\L} \Z / {l^n})$.

\begin{lem}
\label{lem lim SH}
On a des isomorphismes canoniques
\begin{enumerate}
	\item $\varprojlim_V \varprojlim_n \D^0(V, \mathcal{C} \otimes^{\L} \Z / {l^n}) \cong \varprojlim_n \textup{\cyr{SH}}^0(k, C \otimes^{\L} \Z / {l^n})$, o\`u les fl\`eches de transition proviennent de la fonctorialit\'e covariante de la cohomologie \`a support compact.
	\item $\varinjlim_V \D^1(V,  \varinjlim_n \widehat{\mathcal{C}} \otimes^{\L} \Z / {l^n}) \cong \textup{\cyr{SH}}^1(k,  \varinjlim_n \widehat{C} \otimes^{\L} \Z / {l^n})$, o\`u les fl\`eches de transitions sont donn\'ees par le lemme \ref{lem nul lim}.
\end{enumerate}
($V$ d\'ecrit les ouverts non vides de $U$).
\end{lem}

\begin{dem}
\begin{enumerate}
	\item Le lemme \ref{lem inj lim} assure que l'on a $\bigcap_V \varprojlim_n \D^0(V, \mathcal{C} \otimes \Z / {l^n}) = \varprojlim_n \textup{\cyr{SH}}^0(k, C \otimes^{\L} \Z / {l^n})$. On utilise alors la fonctorialit\'e covariante de la cohomologie \`a support compact pour voir l'intersection comme une limite projective : si $V' \subset V$, on a une inclusion canonique (dans $\varprojlim_n \H^0(k, C \otimes^{\L} \Z / {l^n})$) $\varprojlim_n \D^0(V', \mathcal{C} \otimes \Z / {l^n}) \subset \varprojlim_n \D^0(V, \mathcal{C} \otimes \Z / {l^n})$.
	\item Le lemme \ref{lem nul lim} implique que, si $V' \subset V$, le morphisme naturel $ \varinjlim_n \H^1(V,  \widehat{\mathcal{C}} \otimes^{\L} \Z / {l^n}) \rightarrow  \varinjlim_n \H^1(V',  \widehat{\mathcal{C}} \otimes^{\L} \Z / {l^n})$ envoie $\D^1(V, \varinjlim_n \widehat{\mathcal{C}} \otimes^{\L} \Z / {l^n})$ dans $\D^1(V', \varinjlim_n \widehat{\mathcal{C}} \otimes^{\L} \Z / {l^n})$. En consid\'erant les images $$\mathcal{D}^1(V, \varinjlim_n \widehat{\mathcal{C}} \otimes^{\L} \Z / {l^n}) := \textup{Im} \left( \D^1(V,  \varinjlim_n \widehat{\mathcal{C}} \otimes^{\L} \Z / {l^n}) \rightarrow \varinjlim_n \H^1(k, \widehat{C} \otimes^{\L} \Z / {l^n}) \right)$$
dont la r\'eunion sur tous les ouverts $V$ donne exactement $\textup{\cyr{SH}}^1(k, \varinjlim_n \widehat{C} \otimes^{\L} \Z / {l^n})$, on d\'eduit de cette propri\'et\'e un morphisme surjectif
$$\varinjlim_V \D^1(V, \varinjlim_n \widehat{\mathcal{C}} \otimes^{\L} \Z / {l^n}) \rightarrow \textup{\cyr{SH}}^1(k, \varinjlim_n \widehat{C} \otimes^{\L} \Z / {l^n})$$
Ce dernier morphisme est injectif par les principes g\'en\'eraux de cohomologie \'etale : en effet, le morphisme $\varinjlim_V \H^1(V, \varinjlim_n \widehat{\mathcal{C}} \otimes^{\L} \Z / {l^n}) \rightarrow \H^1(k, \varinjlim_n \widehat{C} \otimes^{\L} \Z / {l^n})$ est un isomorphisme par d\'evissage \`a partir du th\'eor\`eme VII.6.7 de \cite{SGA4}.
\end{enumerate}
\end{dem}

On conclut la preuve de la proposition \ref{prop dualite noyau} de la fa\c con suivante : le lemme \ref{lem AV D} fournit une dualit\'e parfaite entre le groupe profini $\varprojlim_n \D^0(V, \mathcal{C} \otimes \Z / {l^n})$ et le groupe discret $\D^1(V, \varinjlim_n \widehat{\mathcal{C}} \otimes^{\L} \Z / {l^n})$. Enfin, en passant \`a la limite sur les ouverts $V$, et en utilisant le lemme \ref{lem lim SH}, on en d\'eduit une dualit\'e parfaite entre le groupe profini $\varprojlim_n \textup{\cyr{SH}}^0(k, C \otimes^{\L} \Z / {l^n})$ et le groupe discret $\textup{\cyr{SH}}^1(k, \varinjlim_n \widehat{C} \otimes^{\L} \Z / {l^n})$. Ceci \'etant valable pour tout nombre premier $l$, on en d\'eduit la proposition \ref{prop dualite noyau}, \`a savoir une dualit\'e parfaite entre le groupe compact $\textup{Ker}(\theta)$ et le groupe discret $\textup{\cyr{SH}}^1(k, \varinjlim_n \widehat{C} \otimes^{\L} \Z / {n})$.
\end{dem}

Pour finir la d\'emonstration du th\'eor\`eme \ref{dualite SH 0}, on va identifier le groupe  $\textup{\cyr{SH}}^1( \varinjlim_n \widehat{C} \otimes^{\L} \Z / n )$ avec le groupe $\textup{\cyr{SH}}^2(\widehat{C})$. Pour ce faire, on consid\`ere les suites exactes :
$$0 \rightarrow \H^1(k, \widehat{C}) \otimes \Z / n \rightarrow \H^1(k, \widehat{C} \otimes^{\L} \Z / {n}) \rightarrow _n \H^2(k, \widehat{C}) \rightarrow 0$$
et pour chaque place $v$ :
$$0 \rightarrow \H^1(\widehat{k}_v, \widehat{C}) \otimes \Z / n \rightarrow \H^1(\widehat{k}_v, \widehat{C} \otimes^{\L} \Z / {n}) \rightarrow _n \H^2(\widehat{k}_v, \widehat{C}) \rightarrow 0$$
On passe \`a la limite inductive sur $n$, et on obtient, puisque les groupes $\H^1(k, \widehat{C})$, $\H^2(k, \widehat{C})$, $\H^1(k_v, \widehat{C})$ et $\H^2(k_v, \widehat{C})$ sont de torsion, des isomorphismes :
$$\H^1(k, \varinjlim_n \widehat{C} \otimes^{\L} \Z / {n}) \cong \H^2(k, \widehat{C})$$
et en prenant le produit sur toutes les places $v$ :
$$\prod_v \H^1(\widehat{k}_v, \varinjlim_n \widehat{C} \otimes^{\L} \Z / {n}) \cong \prod_v \H^2(\widehat{k}_v, \widehat{C})$$
On d\'eduit imm\'ediatement de ces deux isomorphismes l'identification annonc\'ee, \`a savoir 
$$\textup{\cyr{SH}}^1( \varinjlim_n \widehat{C} \otimes^{\L} \Z / n ) \cong \textup{\cyr{SH}}^2(\widehat{C})$$
ce qui conclut la preuve du th\'eor\`eme \ref{dualite SH 0}.

\end{dem}

On conclut cette section par un lemme, qui nous sera utile dans la section \ref{section PT}.
\begin{lem}
\label{lem compl}
Le groupe $\textup{\cyr{SH}}^0_{\wedge}(C)$ est canoniquement isomorphe au groupe $\textup{Ker} \left( \H^0(k, C)^{\wedge} \rightarrow \P^0(C)^{\wedge} \right)$.
\end{lem}

\begin{dem}
Voir proposition 5.4 de \cite{HSz}. On sait en effet que $n \H^0(k_v, C) \subset \H^0(k_v, C)$ est un sous-groupe ouvert d'indice fini (cas des tores et finitude de $H^1(k_v, T_1)$). Cela assure que le morphisme canonique $\P^0(C)_{\wedge} \rightarrow \P^0(C)^{\wedge}$ est injectif. On montre ensuite la surjectivit\'e du morphisme $\textup{Ker } \theta_0 \rightarrow \textup{Ker } \beta_0$ dans le diagramme suivant :
\begin{displaymath}
\xymatrix{
\H^0(k, C)_{\wedge} \ar[r]^{\theta_0} \ar[d] & \P^0(C)_{\wedge} \ar[d] \\
\H^0(k, C)^{\wedge} \ar[r]^{\beta_0} & \P^0(C)^{\wedge}
}
\end{displaymath}
Pour cela, on raisonne comme dans la proposition 5.4 de \cite{HSz}, en
tenant compte des corrigenda \`a propos de la page 120 de \cite{HSz} :
on va montrer que le morphisme $\theta_0 : \H^0(k, C)_{\wedge}
\rightarrow \P^0(k, C)_{\wedge}$ est strict. Soit $n \geq 1$,
consid\'erons le diagramme suivant :
\begin{displaymath}
\xymatrix{
H^2(k, {_n (\textup{Ker } \rho)}) \ar[r] \ar[d] & \H^0(k, C
\otimes^{\L} \Z / n) \ar[r] \ar[d]^{f} & H^1(k, T_{\Z / n}(C)) \ar[d]^{h} \\ 
P^2(k, {_n (\textup{Ker } \rho)}) \ar[r]^{g} \ar[d] & \P^0(k, C
\otimes^{\L} \Z / n) \ar[r] \ar[d] & P^1(k, T_{\Z / n}(C)) \ar[d] \\
H^0(k, \widehat{_n (\textup{Ker } \rho)})^D \ar[r] & \H^0(k, \widehat{C} \otimes^{\L} \Z / n)^D \ar[r] & H^1(k, T_{\Z / n}(\widehat{C}))^D
}
\end{displaymath}
dont les lignes, ainsi que les deux colonnes extr\^emes, sont
exactes (pour les colonnes, c'est la suite de Poitou-Tate pour les modules finis). Montrons que $\textup{Im } f$ est un sous-groupe discret de $\P^0(k, C)$ : on a une suite exacte de groupes topologiques
\begin{equation}
\label{strict}
0 \rightarrow \textup{Im } f \cap \textup{Im } g \rightarrow
\textup{Im } f \rightarrow \textup{Im } h
\end{equation}
o\`u tous les groupes sont munis de la topologie induite par les
topologies ad\'eliques sur $\P^0(k, C \otimes^{\L} \Z / n)$ et $P^1(k,
T_{\Z / n}(C))$. Le groupe $H^1(k, T_{\Z / n}(C))$ est discret, et $h$ est d'image ferm\'ee par
Poitou-Tate, donc $h$ est d'image localement compacte, et donc $h$ est stricte
par \cite{HeR}, 5.29. Donc le groupe $\textup{Im } h$ est discret. De m\^eme, le groupe $P^2(k, {_n (\textup{Ker }
  \rho)})$ est discret, $g$ est d'image ferm\'ee donc localement
compacte par Poitou-Tate, et on en d\'eduit que $\textup{Im } g$ est
\'egalement discret, toujours gr\^ace \`a \cite{HeR}, 5.29. Par
cons\'equent, dans la suite exacte (\ref{strict}), le groupe
$\textup{Im } f$ admet le groupe discret $\textup{Im } f \cap
\textup{Im } g$ comme sous-groupe ouvert (c'est l'image r\'eciproque
de l'ouvert $\{0\}$ du groupe discret $\textup{Im } h$ par le
morphisme $\textup{Im } f \rightarrow \textup{Im } h$), donc le groupe
$\textup{Im } f$ est lui-m\^eme discret.
On consid\`ere alors le diagramme (\ref{diag niv fini}) : on vient de
montrer que l'image du morphisme $f : \H^0(k, C
\otimes^{\L} \Z / n) \rightarrow \P^0(k, C
\otimes^{\L} \Z / n)$ est discr\`ete. Le diagramme (\ref{diag niv
  fini}), ainsi que \cite{HeR}, 5.29, assurent que les morphismes 
$\H^0(k, C) / n \rightarrow \H^0(k, C
\otimes^{\L} \Z / n)$ et $\P^0(k, C) / n \rightarrow \P^0(k, C
\otimes^{\L} \Z / n) $ sont stricts. Donc l'image
de $p_n : \H^0(k, C) / n \rightarrow \P^0(k, C) / n$ s'identifie (topologiquement) \`a un
sous-groupe de l'image de $\H^0(k, C
\otimes^{\L} \Z / n) \rightarrow \P^0(k, C
\otimes^{\L} \Z / n)$, laquelle est discr\`ete. Donc l'image
de $p_n$ est discr\`ete, donc
ferm\'ee, donc localement compacte, donc le morphisme $p_n$ est strict
(par \cite{HeR}, 5.29). Enfin, $\textup{Ker } p_n$ est fini et
$\H^0(k, C) / n$ est discret, donc la limite projective des morphismes
stricts $p_n$ est un morphisme strict $\theta_0 : \H^0(k, C)_{\wedge}
\rightarrow \P^0(k, C)_{\wedge}$ (voir \cite{HSz}, corrigenda).
Ce fait, joint \`a l'injectivit\'e de $\P^0(k, C)_{\wedge} \rightarrow
\P^0(k, C)^{\wedge}$, assure que $\textup{Ker } \theta_0 \rightarrow
\textup{Ker } \beta_0$ est un isomorphisme.
\end{dem}

\subsubsection{Cas o\`u $\rho$ est surjective}

Dans cette section, on montre une dualit\'e parfaite entre les groupes finis $\textup{\cyr{SH}}^0(C)$ et $\textup{\cyr{SH}}^2(\widehat{C})$ sous l'hypoth\`ese que le morphisme $\rho$ est surjectif.

\begin{rem}
{\rm
Si l'on suppose que $\rho$ est surjectif, cela implique que le complexe $C$ est quasi isomorphe \`a $(\textup{Ker} \rho)[1]$, et donc on va ici retrouver un r\'esultat de dualit\'e globale pour un groupe de type multiplicatif (voir par exemple \cite{NSW}, th\'eor\`eme 8.6.7).
}
\end{rem}

\begin{theo}
\label{dualite SH 0 tm}
Supposons que $\rho$ est surjective. Alors il existe une dualit\'e
parfaite de groupes finis, fonctorielle en $M$ :
$$\textup{\cyr{SH}}^0(C) \times \textup{\cyr{SH}}^2(\widehat{C}) \rightarrow \Q / \Z$$
\end{theo}

\begin{dem}
Fixons un nombre premier $l$ et un ouvert $U := \textup{Spec}(\mathcal{O}_{k, S})$ de $\textup{Spec}(\mathcal{O}_k)$ sur lequel $C$ s'\'etend en un complexe surjectif de $U$-tores $\mathcal{C}$, de sorte que $l$ soit inversible sur $U$. On remarque alors les faits suivants :
\begin{itemize}
  \item Le morphisme $\H^0(U, \mathcal{C})\{l\} \rightarrow \H^0(k, C)\{ l \}$ est injectif si $U$ est suffisamment petit. En effet, on consid\`ere un triangle exact de la forme
$$\mathcal{T} \rightarrow \mathcal{C}[-1] \rightarrow \mathcal{F} \rightarrow \mathcal{T}[1]$$
o\`u $\mathcal{T}$ est un $U$-tore (la composante connexe de $\textup{Ker } \rho$) et $\mathcal{F}$ un $U$-sch\'ema en groupes de type multiplicatif fini. On conclut en remarquant que $H^0(U, \mathcal{F})\{l\} \rightarrow H^0(k, F)\{ l \}$ est un isomorphisme (puisque $F$ est localement constant) et $H^1(U, \mathcal{F})\{l\} \rightarrow H^1(k, F)\{ l \}$ et $H^1(U, \mathcal{T})\{l\} \rightarrow H^1(k, T)\{ l \}$ sont injectifs (voir la proposition 2.9 de \cite{Mil} pour le premier, et la proposition 4.1 et les corrigenda de \cite{HSz}, ainsi que la preuve de la proposition \ref{prop etale-galois} pour le second).
\item De m\^eme, $\H^2(U, \widehat{\mathcal{C}})\{l\} \rightarrow \H^2(k, \widehat{C})\{ l \}$ est injectif.
  \item Pour tout ouvert $V$ de $U$, les groupes $\D^0(V, \mathcal{C})\{l\}$ et $\D^2(V, \widehat{\mathcal{C}})\{l\}$ sont finis. En effet, en utilisant \`a nouveau le d\'evissage pr\'ec\'edent, on sait que $H^1(V, \mathcal{F})\{l\}$ est fini, que $D^1(V, \mathcal{T})\{l\}$ est fini, et que $H^0(k_v, F)$ est fini pour toute place $v$, donc le lemme du serpent assure la finitude de $\D^0(V, \mathcal{C})\{l\}$. Pour le groupe $\D^2(V, \widehat{\mathcal{C}})\{l\}$, on utilise le m\^eme d\'evissage.
    \item Comme dans la preuve du lemme 4.7 de \cite{HSz}, on d\'eduit des points pr\'ec\'edents qu'il existe un ouvert $U_0$ dans $U$ tel que $\D^0(V, \mathcal{C})\{l\} = \textup{\cyr{SH}}^0(k, C)\{l\}$ et $\D^2(V, \widehat{\mathcal{C}})\{l\} = \textup{\cyr{SH}}^2(k, \widehat{C})\{l\}$ pour tout $V$ dans $U_0$.
      \item On conlut la preuve du th\'eor\`eme \ref{dualite SH 0 tm} gr\^ace au corollaire \ref{coro D^i}, en remarquant que $\H^0(k_v, C) \rightarrow \H^0(k_v, C)^{\wedge}$ est un isomorphisme car $\H^0(k_v, C)$ est fini sous l'hypoth\`ese de surjectivit\'e de $\rho$.

\end{itemize}
\end{dem}

\section{Deux suites de Poitou-Tate}
\label{section PT}
\subsection{Cas o\`u $\textup{Ker } \rho$ est fini}

\subsubsection{La suite exacte de Poitou-Tate}
\begin{theo}
\label{theo PT}
Soit $C = [T_1 \xrightarrow{\rho} T_2 ]$ un complexe de tores d\'efini
sur $k$, avec $\textup{Ker}(\rho)$ fini. On a alors une suite exacte
fonctorielle en $C$ :
\begin{displaymath}
\xymatrix{
0 \ar[r] & \H^{-1}(k, C) \ar[r] & \P^{-1}(k, C) \ar[r] & \H^2(k, \widehat{C})^D \ar[d] \\
& \H^1(k, \widehat{C})^D \ar[d] & \P^0(k, C)^{\wedge} \ar[l] & \H^0(k, C)^{\wedge} \ar[l] \\
& \H^1(k,C) \ar[r] & \P^1(k, C) \ar[r] & \H^0(k, \widehat{C})^D \ar[d] \\
0 & \H^{-1}(k, \widehat{C})^D \ar[l] & \P^2(k, C) \ar[l] & \H^2(k, C) \ar[l]
}
\end{displaymath}

On dispose \'egalement, sous les m\^emes hypoth\`eses, de la suite exacte duale :
\begin{displaymath}
\xymatrix{
0 \ar[r] & \H^{-1}(k, \widehat{C})^{\wedge} \ar[r] & \P^{-1}(k, \widehat{C})^{\wedge} \ar[r] & \H^2(k, C)^D \ar[d] \\
& \H^1(k, C)^D \ar[d] & \P^0(k, \widehat{C}) \ar[l] & \H^0(k, \widehat{C}) \ar[l] \\
& \H^1(k,\widehat{C}) \ar[r] & \P^1(k, \widehat{C})_{\textup{tors}} \ar[r] & \left( \H^0(k, C)^D \right)_{\textup{tors}} \ar[d] \\
0 & \H^{-1}(k, C)^D \ar[l] & \P^2(k, \widehat{C}) \ar[l] & \H^2(k, \widehat{C}) \ar[l]
}
\end{displaymath}

\end{theo}
	
\begin{dem}
\begin{itemize}
	\item Montrons l'exactitude de la deuxi\`eme ligne. \\
Pour cela, on commence par le lemme suivant :
\begin{lem}
\label{lem PT1}
On suppose $\textup{Ker } \rho$ fini. Alors la suite 
$$ \varprojlim_n \H^0(k, C \otimes^{\L} \Z / n) \rightarrow \varprojlim_n \P^0(k, C \otimes^{\L} \Z / n) \rightarrow  \H^0(k, \varinjlim_n \widehat{C} \otimes \Z / n)^D$$ 
est exacte.
\end{lem}

\begin{dem}
Consid\'erons le diagramme \`a lignes exactes suivant :
\begin{changemargin}{-2.3cm}{1cm}
\begin{displaymath}
\xymatrix{
H^0(k, T_{\Z / n}(C)) \ar[r] \ar[d] & H^2(k, {_n (\textup{Ker } \rho)}) \ar[r] \ar[d] & \H^0(k, C \otimes^{\L} \Z / n) \ar[r] \ar[d] & H^1(k, T_{\Z / n}(C)) \ar[r] \ar[d] & H^3(k, {_n (\textup{Ker } \rho)}) \ar[d]^{\simeq} \\
P^0(k, T_{\Z / n}(C)) \ar[r] \ar[d] & P^2(k, {_n (\textup{Ker } \rho)}) \ar[r] \ar[d] & \P^0(k, C \otimes^{\L} \Z / n) \ar[r] \ar[d] & P^1(k, T_{\Z / n}(C)) \ar[r] \ar[d] & P^3(k, {_n (\textup{Ker } \rho)})  \\
H^2(k, T_{\Z / n}(\widehat{C})^D \ar[r] & H^0(k, \widehat{_n (\textup{Ker } \rho)})^D \ar[r] & \H^0(k, \widehat{C} \otimes^{\L} \Z / n)^D \ar[r] & H^1(k, T_{\Z / n}(\widehat{C}))^D \ar[r] & 0
}
\end{displaymath}
\end{changemargin}
Les groupes de la colonne de gauche \'etant finis ou compacts, et l'image du morphisme $P^0(k, T_{\Z / n}(C)) \rightarrow P^2(k, {_n (\textup{Ker } \rho)})$ \'etant finie (car compacte et discr\`ete), le th\'eor\`eme 7.3 de \cite{Jen} assure que l'on a un diagramme \`a lignes exactes (o\`u $\varprojlim^{(1)}$ d\'esigne le foncteur d\'eriv\'e du foncteur $\varprojlim$) :

\begin{eqnarray}
\label{PT2}
\xymatrix{
\varprojlim_n H^2(k, {_n (\textup{Ker } \rho)}) \ar[r] \ar[d] & \varprojlim_n \H^0(k, C \otimes^{\L} \Z / n) \ar[r] \ar[d] & \varprojlim_n Q_1^n \ar[d] \ar[r] & \varprojlim_n^{(1)} P_1^n \ar[d] \\
 \varprojlim_n P^2(k, {_n (\textup{Ker } \rho)}) \ar[r] \ar[d] & \varprojlim_n \P^0(k, C \otimes^{\L} \Z / n) \ar[r] \ar[d] & \varprojlim_n Q_2^n \ar[d] \ar[r] & \varprojlim_n^{(1)} P_2^n \\
\varprojlim_n H^0(k, \widehat{_n (\textup{Ker } \rho)})^D \ar[r] & \varprojlim_n \H^0(k, \widehat{C} \otimes^{\L} \Z / n)^D \ar[r] & \varprojlim_n H^1(k, T_{\Z / n}(\widehat{C}))^D \ar[r] & 0
}
\end{eqnarray}
o\`u les groupes $P_1^n$, $P_2^n$, $Q_1^n$ et $Q_2^n$ sont d\'efinis par les suites exactes de groupes topologiques suivantes :
\begin{equation} \label{sec P_1}
H^0(k, T_{\Z / n}(C)) \rightarrow H^2(k, {_n (\textup{Ker } \rho)}) \rightarrow P_1^n \rightarrow 0 \end{equation}
$$0 \rightarrow P_1^n \rightarrow \H^0(k, C \otimes^{\L} \Z / n) \rightarrow H^1(k, T_{\Z / n}(C))$$

\begin{equation} \label{sec P_2}
P^0(k, T_{\Z / n}(C)) \rightarrow P^2(k, {_n (\textup{Ker } \rho)}) \rightarrow P_2^n \rightarrow 0 \end{equation}
$$0 \rightarrow P_2^n \rightarrow \P^0(k, C \otimes^{\L} \Z / n) \rightarrow P^1(k, T_{\Z / n}(C))$$

$$H^2(k, {_n (\textup{Ker } \rho)}) \rightarrow \H^0(k, C \otimes^{\L} \Z / n) \rightarrow Q_1^n \rightarrow 0$$
$$0 \rightarrow Q_1^n \rightarrow H^1(k, T_{\Z / n}(C)) \rightarrow H^3(k, {_n (\textup{Ker } \rho)})$$

$$P^2(k, {_n (\textup{Ker } \rho)}) \rightarrow \P^0(k, C \otimes^{\L} \Z / n) \rightarrow Q_2^n \rightarrow 0$$
$$0 \rightarrow Q_2^n \rightarrow P^1(k, T_{\Z / n}(C)) \rightarrow P^3(k, {_n (\textup{Ker } \rho)})$$

Or ces groupes s'int\`egrent dans le diagramme \`a lignes exactes suivant :
\begin{eqnarray}
\label{PT3}
\xymatrix{
0 \ar[r] & \varprojlim_n Q_1^n \ar[r] \ar[d] & \varprojlim_n H^1(k, T_{\Z / n}(C)) \ar[r] \ar[d] & \varprojlim_n R_1^n \ar[d] \\
0 \ar[r] & \varprojlim_n Q_2^n \ar[r] & \varprojlim_n P^1(k, T_{\Z / n}(C)) \ar[r] \ar[d] & \varprojlim_n R_2^n \\
& & \varprojlim_n H^1(k, T_{\Z / n}(\widehat{C}))^D &
}
\end{eqnarray}
o\`u $R_1^n$ et $R_2^n$ sont les noyaux respectifs des morphismes $H^3(k, {_n (\textup{Ker } \rho)}) \rightarrow \H^1(k, C \otimes^{\L} \Z / n)$ et $P^3(k, {_n (\textup{Ker } \rho)}) \rightarrow \P^1(k, C \otimes^{\L} \Z / n)$. Les groupes $H^3(k, {_n (\textup{Ker } \rho)})$ et $P^3(k, {_n (\textup{Ker } \rho)})$ \'etant isomorphes par Poitou-Tate, l'exactitude \`a gauche du foncteur $\varprojlim_n$ assure que la fl\`eche $\varprojlim_n R_1^n \rightarrow \varprojlim_n R_2^n$ est injective.

Montrons alors le lemme : soit $\alpha \in \varprojlim_n \P^0(k, C \otimes^{\L} \Z / n)$ d'image nulle dans $\varprojlim_n \H^0(k, \widehat{C} \otimes^{\L} \Z / n)^D$. On veut montrer qu'un tel \'el\'ement se rel\`eve dans $\varprojlim_n \H^0(k, C \otimes^{\L} \Z / n)$. Pour cela, on pousse $\alpha$ dans $\varprojlim_n P^1(k, T_{\Z / n}(C))$. Par fonctorialit\'e, l'\'el\'ement $\beta \in \varprojlim_n P^1(k, T_{\Z / n}(C))$ ainsi obtenu s'envoie sur $0$ dans $\varprojlim_n H^1(k, T_{\Z / n}(\widehat{C}))^D$. 

Par Poitou-Tate, la deuxi\`eme colonne du diagramme (\ref{PT3}) est exacte (on utilise ici la finitude de $\textup{\cyr{SH}}^1(k, T_{\Z / n}(C))$ pour passer \`a la limite projective dans la suite exacte de Poitou-Tate), et donc $\beta$ se rel\`eve en un \'el\'ement $\gamma \in \varprojlim_n H^1(k, T_{\Z / n}(C))$. Or $\beta$ provient de $\varprojlim_n Q_2^n$ (c'est l'image de $\alpha$), donc l'image de $\beta$ dans $\varprojlim_n R_2^n$ est nulle, et par injectivit\'e de la fl\`eche $\varprojlim_n R_1^n \rightarrow \varprojlim_n R_2^n$, l'image de $\gamma$ dans $\varprojlim_n R_1^n$ est nulle. Donc par exactitude de la premi\`ere ligne du diagramme (\ref{PT3}), $\gamma$ provient de $\varprojlim_n Q_1^n$. 

Revenons alors au diagramme initial (\ref{PT2}) : on est dans la configuration suivante : $\alpha \in \varprojlim_n \P^0(k, C \otimes^{\L} \Z / n)$, son image $\gamma \in \varprojlim_n Q_2^n$ provient d'un \'el\'ement $\gamma' \in \varprojlim_n Q_1^n$. Montrons d\'esormais que $\gamma'$ se rel\`eve dans $\varprojlim_n \H^0(k, C \otimes^{\L} \Z / n)$ : il suffit pour cela de montrer que le morphisme $\varprojlim_n^{(1)} P_1^n \rightarrow \varprojlim_n^{(1)} P_2^n$ est injectif.

Or $T_{\Z / n}(C)$ est fini, donc $\varprojlim_n^{(1)} H^0(k, T_{\Z / n}(C)) = 0$, et donc par la suite exacte (\ref{sec P_1}), on sait que 
\begin{equation} \label{iso1}
\varprojlim_n^{(1)} H^2(k, {_n (\textup{Ker } \rho)}) \cong \varprojlim_n^{(1)} P_1^n \end{equation}
De m\^eme, le groupe $P^0(k,T_{\Z / n}(C))$ est compact, donc $\varprojlim_n^{(1)} P^0(k,T_{\Z / n}(C)) = 0$ et donc 
\begin{equation} \label{iso2}
\varprojlim_n^{(1)} P^2(k, {_n (\textup{Ker } \rho)}) \cong \varprojlim_n^{(1)} P_2^n \end{equation}
gr\^ace \`a la suite exacte (\ref{sec P_2}). Si $I_n$ d\'esigne l'image de $H^2(k, {_n (\textup{Ker } \rho)})$ dans $P^2(k, {_n (\textup{Ker } \rho)}) $, puisque $\textup{\cyr{SH}}^2(k, {_n (\textup{Ker } \rho)})$ est fini, on a un isomorphisme
\begin{equation} \label{iso3}
\varprojlim_n^{(1)} H^2(k, {_n (\textup{Ker } \rho)}) \cong \varprojlim_n^{(1)} I_n \end{equation}
Or par Poitou-Tate, le conoyau de $H^2(k, {_n (\textup{Ker } \rho)}) \rightarrow P^2(k, {_n (\textup{Ker } \rho)})$ est exactement le groupe fini $H^0(k,\widehat{ _n (\textup{Ker } \rho)})^D$, donc on a une suite exacte
$$0 \rightarrow I_n \rightarrow P^2(k, {_n (\textup{Ker } \rho)}) \rightarrow H^0(k, \widehat{ _n (\textup{Ker } \rho)})^D \rightarrow 0$$
qui apr\`es passage \`a la limite projective donne la suite exacte :
$\varprojlim_n H^0(k,\widehat{ _n (\textup{Ker } \rho)})^D  \rightarrow  \varprojlim_n^{(1)} I_n \rightarrow \varprojlim_n^{(1)} P^2(k, {_n (\textup{Ker } \rho)}) \rightarrow 0$

Or par hypoth\`ese $\textup{Ker}(\rho)$ est fini, donc $\varprojlim_n H^0(k,\widehat{ _n (\textup{Ker } \rho)})^D = 0$, donc cela prouve que le morphisme 
\begin{equation} \label{iso4} 
\varprojlim_n^{(1)} I_n \rightarrow \varprojlim_n^{(1)} P^2(k, {_n (\textup{Ker } \rho)}) \end{equation}
est un isomorphisme.

Ainsi, en combinant les isomorphismes (\ref{iso1}), (\ref{iso2}), (\ref{iso3}) et (\ref{iso4}), on a bien montr\'e que le morphisme $\varprojlim_n^{(1)} P_1^n \rightarrow \varprojlim_n^{(1)} P_2^n$ \'etait un isomorphisme.

On a donc d\'esormais un \'el\'ement $\tau \in \varprojlim_n \H^0(k, C \otimes^{\L} \Z / n)$ relevant $\gamma' \in \varprojlim_n Q_1^n$. Notons $\tau'$ l'image de $\tau$ dans $\varprojlim_n \P^0(k, C \otimes^{\L} \Z / n)$. Par commutativit\'e et exactitude du diagramme (\ref{PT2}), $\alpha$ et $\tau'$ dans $\varprojlim_n \P^0(k, C \otimes^{\L} \Z / n)$ ont m\^eme image dans $\varprojlim_n Q_2^n$. Donc la diff\'erence se rel\`eve dans $\varprojlim_n P^2(k, {_n (\textup{Ker } \rho)})$. Mais l'hypoth\`ese de finitude sur $\textup{Ker } \rho$ assure que le groupe $\varprojlim_n P^2(k, {_n (\textup{Ker } \rho)})$ est trivial, donc $\tau' = \alpha$, ce qui conclut la preuve du lemme \ref{lem PT1}.
\end{dem}

On consid\`ere alors \`a nouveau le diagramme (\ref{diag lim}), commutatif, dont les deux premi\`eres lignes sont exactes :
\begin{displaymath}
\xymatrix{
0 \ar[r] & \H^0(k, C)_{\wedge} \ar[d]^{\theta_0} \ar[r] & \varprojlim_n \H^0(k, C \otimes^{\L} \Z / n) \ar[r] \ar[d]^{\theta} & Q_1 \ar[d]^{\beta} \ar[r] & 0 \\
0 \ar[r] & \P^0(k, C)_{\wedge} \ar[r] \ar[d]^{\gamma_0'} & \varprojlim_n \P^0(k, C \otimes^{\L} \Z / n) \ar[r] \ar[d] & Q_2 \ar[r] & 0 \\
& \H^1(k, \widehat{C})^D \ar[r] & \H^0(k, \varinjlim_n \widehat{C} \otimes \Z / n)^D & & 
}
\end{displaymath}
On sait que la colonne centrale est exacte par le lemme \ref{lem PT1}. Or la colonne de gauche est un complexe par la loi de r\'eciprocit\'e globale du corps de classes. De plus, le morphisme $\beta$ est injectif, donc une chasse au diagramme assure que la colonne de gauche est exacte, \`a savoir la suite 
$$ \H^0(k, C)_{\wedge} \rightarrow  \P^0(k, C)_{\wedge} \rightarrow \H^1(k, \widehat{C})^D$$
est exacte. On prend la compl\'etion profinie de cette suite, on obtient le complexe 
$$ \H^0(k, C)^{\wedge} \rightarrow  \P^0(k, C)^{\wedge} \rightarrow \H^1(k, \widehat{C})^D$$
Celui-ci est exact par les consid\'erations topologiques suivantes (voir par exemple \cite{HSz}, preuve du th\'eor\`eme 5.6) : en reprenant la preuve du lemme \ref{lem PT1}, et en utilisant le fait que 
$$\textup{Coker}\left( \varprojlim_n H^1(k, T_{\Z / n}(C)) \rightarrow \varprojlim_n P^1(k, T_{\Z / n}(C)) \right)$$ 
est profini (suite de Poitou-Tate pour les modules finis), on montre facilement que 
$$\textup{Coker}\left( \H^0(k, C)_{\wedge} \rightarrow \P^0(k, C)_{\wedge} \right)$$ est profini, donc $I := \textup{Im} \left( \P^0(k, C)_{\wedge} \xrightarrow{\gamma_0'} \H^1(k, \widehat{C})^D \right)$ est un sous-groupe ferm\'e profini de $\H^1(k, \widehat{C})^D$. Or en compl\'etant la suite exacte
$$\H^0(k, C)_{\wedge} \rightarrow \P^0(k, C)_{\wedge} \xrightarrow{\gamma_0'} I \rightarrow 0$$
on obtient une suite exacte
$$\H^0(k, C)^{\wedge} \rightarrow \P^0(k, C)^{\wedge} \xrightarrow{\gamma_0'} I^{\wedge} \rightarrow 0$$
et $I$ \'etant profini, on a $I^{\wedge} = I$ et donc $I^{\wedge} \rightarrow \H^1(k, \widehat{C})^D$ est injective, ce qui assure l'exactitude de la deuxi\`eme ligne du diagramme du th\'eor\`eme \ref{theo PT}, \`a savoir 
$$\H^0(k, C)^{\wedge} \rightarrow \P^0(k, C)^{\wedge} \rightarrow \H^1(k, \widehat{C})^D$$

	\item Pour la troisi\`eme ligne, le raisonnement est plus direct :
on consid\`ere cette fois le diagramme commutatif suivant :
\begin{eqnarray}
  \label{diag ligne 3}
\xymatrix{
H^3(k, {_n (\textup{Ker } \rho)}) \ar[r] \ar[d]^{\simeq} & \H^1(k, C \otimes^{\L} \Z / n) \ar[r] \ar[d] & H^2(k, T_{\Z / n}(C)) \ar[r] \ar[d] & H^4(k, {_n (\textup{Ker } \rho)}) \ar[d]^{\simeq} \\
P^3(k, {_n (\textup{Ker } \rho)}) \ar[r] & \P^1(k, C \otimes^{\L} \Z / n) \ar[r] \ar[d] & P^2(k, T_{\Z / n}(C)) \ar[r] \ar[d] & P^4(k, {_n (\textup{Ker } \rho)})  \\
 & \H^{-1}(k, \widehat{C} \otimes^{\L} \Z / n)^D \ar[r]^{\simeq} & H^0(k, T_{\Z / n}(\widehat{C}))^D &
}
\end{eqnarray}
On sait que la troisi\`eme colonne est exacte, et que les deux fl\`eches verticales extr\^emes sont des isomorphismes. Une chasse au diagramme assure que la deuxi\`eme colonne est exacte. En outre, les finitudes de $H^3(k, {_n (\textup{Ker } \rho)})$ et $\textup{\cyr{SH}}^2(k, T_{\Z / n}(C))$ assurent celle de $\textup{\cyr{SH}}^1(k, C \otimes^{\L} \Z / n)$, donc on peut prendre la limite projective de la deuxi\`eme colonne pour obtenir la suite exacte suivante (puisque $\varprojlim_n^{(1)} \textup{\cyr{SH}}^1(k, C \otimes^{\L} \Z / n) = 0$) :
$$\varprojlim_n \H^1(k, C \otimes^{\L} \Z / n) \rightarrow \varprojlim_n \P^1(k, C \otimes^{\L} \Z / n) \rightarrow \varprojlim_n \H^{-1}(k, \widehat{C} \otimes^{\L} \Z / n)^D$$

Or on dispose pour tout $n$ de la suite exacte naturelle,
$$0 \rightarrow \H^1(k, C) / n \rightarrow \H^1(k, C \otimes^{\L} \Z / n) \rightarrow {_n \H^2(k, C)} \rightarrow 0$$
D'o\`u un diagramme commutatif \`a lignes exactes
\begin{displaymath}
\xymatrix{
0 \ar[r] & \H^1(k, C)_ {\wedge} \ar[r] \ar[d] & \varprojlim_n \H^1(k, C \otimes^{\L} \Z / n) \ar[r] \ar[d] & Q_1 \ar[d]^{\beta'} \ar[r] & 0 \\
0 \ar[r] & \P^1(k, C)_{\wedge} \ar[r] \ar[d] & \varprojlim_n \P^1(k, C \otimes^{\L} \Z / n) \ar[r] \ar[d] & Q_2 \ar[d] \ar[r] & 0 \\
0 \ar[r] & (\H^0(k, \widehat{C})^D)_{\wedge} \ar[r] & \varprojlim_n \H^{-1}(k, \widehat{C} \otimes^{\L} \Z / n)^D \ar[r] & Q_3 \ar[r] & 0
}
\end{displaymath}
Or on a vu que la deuxi\`eme colonne \'etait exacte, et $\textup{Ker}(\beta')$ est contenu dans ${\varprojlim_n} {_n \textup{\cyr{SH}}^2(C)}$, or $\textup{\cyr{SH}}^2(C)$ est un groupe fini (voir th\'eor\`eme \ref{theo SH 2}), donc son module de Tate est nul, donc $\beta'$ est injective, donc une chasse au diagramme assure que la premi\`ere colonne est exacte.
Reste \`a montrer que l'on peut "enlever" les compl\'etions dans cette suite. Pour cela, on utilise la finitude de $\textup{Ker}(\rho)$ pour montrer que les groupes $\H^1(k, C)$ et $\P^1(k, C)$ sont de $N$-torsion pour un $N$ suffisamment grand. En effet, on peut d\'evisser $C$ dans une suite exacte de complexes de la forme
\begin{displaymath}
\xymatrix{
0 \ar[r] & \textup{Ker}(\rho) \ar[r] \ar[d]^{\rho} & T_1 \ar[r] \ar[d] & T_1' \ar[d] \ar[r] & 0 \\
0 \ar[r] &0 \ar[r] & T_2 \ar[r]^{=} & T_2 \ar[r] & 0 \\
}
\end{displaymath}
o\`u $T_1'$ est le $k$-sous-tore de $T_2$ image de $T_1$ par $\rho$. Notons alors $S := T_2 / T_1'$ le tore quotient. Le diagramme pr\'ec\'edent induit alors une suite exacte en hypercohomologie :
$$H^2(k, \textup{Ker}(\rho)) \rightarrow \H^1(k, C) \rightarrow H^1(k, S)$$
Or, par Hilbert 90 et par restriction-corestriction, le groupe $H^1(k, S)$ est de $r$-torsion, o\`u $r := [L:k]$ est le degr\'e d'une extension $L/k$ d\'eployant le tore $S$, et le groupe $H^2(k, \textup{Ker}(\rho))$ est de $r'$-torsion, o\`u $r'$ est le cardinal du groupe fini $\textup{Ker}(\rho)(\overline{k})$. Donc $\H^1(k, C)$ est de $N = r r'$-torsion. Cela assure que $\H^1(k, C) \cong \H^1(k, C)_{\wedge}$. De m\^eme, $\P^1(k, C) \cong \P^1(k, C)_{\wedge}$. Avec ces identifications, l'exactitude de la suite $\H^1(k, C)_ {\wedge} \rightarrow \P^1(k, C)_ {\wedge} \rightarrow (\H^0(k, \widehat{C})^D)_{\wedge}$ implique imm\'ediatement celle de la suite 
$$\H^1(k, C) \rightarrow \P^1(k, C) \rightarrow \H^0(k, \widehat{C})^D$$
c'est-\`a-dire l'exactitude de la ligne 3 dans le th\'eor\`eme \ref{theo PT}.

	\item Montrons l'exactitude de la derni\`ere ligne : on s'int\'eresse au diagramme commutatif suivant dont les colonnes sont exactes :
\begin{eqnarray}
\label{diag ligne 4}
\xymatrix{
0 \ar[d] & 0 \ar[d] & 0 \ar[d]	\\
\H^1(k, C) / n \ar[r] \ar[d] & \P^1(k, C) / n \ar[r] \ar[d] & \left( _n \H^0(k, \widehat{C}) \right)^D \ar[d] \\
\H^1(k, C \otimes^{\L} \Z / n) \ar[r] \ar[d] & \P^1(k, C \otimes^{\L} \Z / n) \ar[r] \ar[d] & \H^{-1}(k, \widehat{C} \otimes^{\L} \Z / n)^D \ar[d]  \\
_n \H^2(k, C) \ar[r] \ar[d] & { _n \P^2(k, C)} \ar[r] \ar[d] & \left( \H^{-1}(k, \widehat{C})/n \right)^D \ar[d] \\
0 & 0 & 0
}
\end{eqnarray}
Or on a montr\'e dans la preuve du point pr\'ec\'edent l'exactitude de la deuxi\`eme ligne :
$$\H^1(k, C \otimes^{\L} \Z / n) \rightarrow \P^1(k, C \otimes^{\L} \Z / n) \rightarrow \H^{-1}(k, \widehat{C} \otimes^{\L} \Z / n)^D$$
 On passe \`a la limite inductive sur $n$ dans le diagramme (\ref{diag ligne 4}). Puisque $\H^0(k, \widehat{C})_{\textup{tors}}$ est fini, le module de Tate $T(\H^0(k, \widehat{C}))$ est nul, et donc on obtient le diagramme suivant dont la premi\`ere ligne est exacte ($\H^2(k, C)$ et $\P^2(k, C)$ sont de torsion, car $\H^2(\widehat{\mathcal{O}}_v, \mathcal{C}) = 0$ pour presque toute place $v$ de $k$) :
\begin{displaymath}
\xymatrix{
\H^1(k, \varinjlim_n C \otimes^{\L} \Z / n) \ar[r] \ar[d] & \P^1(k, \varinjlim_n C \otimes^{\L} \Z / n) \ar[r] \ar[d] & \left(  \varprojlim_n \H^{-1}(k, \widehat{C} \otimes^{\L} \Z / n) \right)^D \ar[d]^{\simeq} \ar[r] & 0 \\
\H^2(k, C) \ar[r] & \P^2(k, C) \ar[r] & \left( \H^{-1}(k, \widehat{C})_{\wedge} \right)^D &
}
\end{displaymath}
(la surjectivit\'e de la derni\`ere fl\`eche de la premi\`ere ligne provient de la finitude du conoyau de $ \P^1(k, C \otimes^{\L} \Z / n) \rightarrow \H^{-1}(k, \widehat{C} \otimes^{\L} \Z / n)^D$, laquelle est une cons\'equence du diagramme (\ref{diag ligne 3})).

Or les deux premi\`eres fl\`eches verticales sont des isomorphismes \'egalement, puisque $ \varinjlim_n \H^1(k, C) / n = 0$ et $\varinjlim_n \P^1(k, C) / n = 0$ car $\H^1(k, C)$ et $\P^1(k, C)$ sont de torsion (on rappelle que $\H^1(\widehat{\mathcal{O}}_v, \mathcal{C}) = 0$ pour presque toute place $v$). Donc la suite 
$$\H^2(k, C) \rightarrow \P^2(k, C) \rightarrow \left( \H^{-1}(k, \widehat{C})_{\wedge} \right)^D \rightarrow 0$$
est exacte. En outre, le groupe discret de type fini $\H^{-1}(k, \widehat{C})$ a m\^eme dual que son compl\'et\'e $\H^{-1}(k, \widehat{C})_{\wedge}$. D'o\`u l'exactitude de la derni\`ere ligne du diagramme du th\'eor\`eme.

	\item Pour la premi\`ere ligne, l'exactitude de la suite 
$$0 \rightarrow \H^{-1}(k, C) \rightarrow \P^{-1}(k, C) \rightarrow \H^2(k, \widehat{C})^D$$
se d\'eduit imm\'ediatement du d\'ebut de la suite de Poitou-Tate associ\'ee au module fini $\textup{Ker}(\rho)$.

	\item Montrons enfin l'exactitude dans les "coins" du diagramme : c'est exactement la traduction des dualit\'es globales montr\'ees pr\'ec\'edemment, \`a savoir $\textup{\cyr{SH}}^0_{\wedge}(C) \cong \textup{\cyr{SH}}^2(\widehat{C})^D$ pour le coin en haut \`a droite (th\'eor\`eme \ref{dualite SH 0}), $\textup{\cyr{SH}}^1(C) \cong \textup{\cyr{SH}}^1(\widehat{C})^D$ pour le coin au milieu \`a gauche (th\'eor\`eme \ref{theo SH 1}), et enfin $\textup{\cyr{SH}}^2(C) \cong \textup{\cyr{SH}}^0(\widehat{C})^D$ pour le coin en bas \`a droite (th\'eor\`eme \ref{theo SH 2}).

          \item Consid\'erons maintenant la suite exacte duale : les raisonnements sont similaires \`a ceux des points pr\'ec\'edents. Pour la premi\`ere ligne, il suffit de dualiser la derni\`ere ligne de la premi\`ere suite de Poitou-Tate, ou alors de montrer par d\'evissage l'exactitude de la suite $\varprojlim_n \H^{-1}(k, \widehat{C} \otimes^{\L} \Z / n) \rightarrow \varprojlim_n \P^{-1}(k, \widehat{C} \otimes^{\L} \Z / n) \rightarrow \varprojlim_n \H^{1}(k, C \otimes^{\L} \Z / n)^D$ et d'utiliser la finitude de $\textup{\cyr{SH}}^0(\widehat{C})$ (voir lemme \ref{lem SH 0 fini}). Pour la deuxi\`eme ligne, on montre par d\'evissage l'exactitude de $\varinjlim_n \H^{-1}(k, \widehat{C} \otimes^{\L} \Z / n) \rightarrow \varinjlim_n \P^{-1}(k, \widehat{C} \otimes^{\L} \Z / n) \rightarrow \varinjlim_n \H^{1}(k, C \otimes^{\L} \Z / n)^D$, puis on utilise la finitude de $\textup{\cyr{SH}}^2(C)$ et le fait que $\H^0(k, \widehat{C})$ et $\P^0(k, \widehat{C})$ soient de torsion (puisque $\textup{Ker } \rho$ est fini) pour en d\'eduire l'exactitude de la deuxi\`eme ligne. En ce qui concerne la troisi\`eme ligne, on montre par d\'evissage l'exactitude de la suite $\varinjlim_n \H^{0}(k, \widehat{C} \otimes^{\L} \Z / n) \rightarrow \varinjlim_n \P^{0}(k, \widehat{C} \otimes^{\L} \Z / n) \rightarrow \varinjlim_n \H^{0}(k, C \otimes^{\L} \Z / n)^D$, et la finitude de $\textup{\cyr{SH}}^1(C)$ assure l'exactitude de le troisi\`eme ligne du diagramme. Pour la quatri\`eme ligne, on montre l'exactitude de $\varprojlim_n \H^{2}(k, \widehat{C} \otimes^{\L} \Z / n) \rightarrow \varprojlim_n \P^{2}(k, \widehat{C} \otimes^{\L} \Z / n) \rightarrow \varprojlim_n \H^{-2}(k, C \otimes^{\L} \Z / n)^D$ et on conclut par finitude de $\H^3(k, \widehat{C})$ et par le fait que $\H^2(k, \widehat{C})$ et $\P^2(k, \widehat{C})$ soient de $N$-torsion pour un certain entier $N$. 
Enfin, pour les coins de ce diagramme, on utilise \`a nouveau les th\'eor\`emes \ref{theo SH 1}, \ref{theo SH 2} et \ref{dualite SH 0}
\end{itemize}
\end{dem}

\subsubsection{Lien avec une suite de Borovoi et description explicite des accouplements}

Dans toute cette partie, l'hypoth\`ese de finitude de $\textup{Ker}(\rho)$ est essentielle. On d\'efinit le complexe $\check{C} := \left[ X_*(T_1) \xrightarrow{\rho_*} X_*(T_2) \right]$, o\`u $X_*(T_i)$ d\'esigne le modules des cocaract\`eres du tore $T_i$. Suivant Borovoi (voir \cite{BorAMS}, chapitre 4), on consid\`ere la suite exacte courte de complexes :
\begin{equation}
\label{sec Bor}
0 \rightarrow \check{C} \otimes \overline{k}^* \rightarrow \check{C} \otimes \overline{\A}^* \rightarrow \check{C} \otimes \overline{\C}^* \rightarrow 0
\end{equation}
o\`u $\overline{\A}$ d\'esigne $\A_k \otimes_k \overline{k}$ et $\overline{\C}^* := \overline{\A}^*  / \overline{k}^*$.
\'Ecrivons la suite exacte d'hypercohomologie associ\'ee :
\begin{equation}
\label{sel Bor}
\dots \rightarrow \H^i(k, \check{C} \otimes \overline{k}^*) \rightarrow  \H^i(k, \check{C} \otimes \overline{\A}^*) \rightarrow  \H^i(k, \check{C} \otimes \overline{\C}^*) \rightarrow  \H^{i+1}(k, \check{C} \otimes \overline{k}^*) \rightarrow \dots
\end{equation}

L'objectif de cette section est de comparer cette suite exacte avec la suite de Poitou-Tate du th\'eor\`eme \ref{theo PT}, et d'obtenir au passage une description explicite (en terme de cup-produit) des accouplements des th\'eor\`emes \ref{theo SH 1}, \ref{theo SH 2} et \ref{dualite SH 0}.

Borovoi identifie certains des groupes qui apparaissent dans la suite (\ref{sel Bor}) (voir \cite{BorAMS}, chapitre 4), de fa\c con fonctorielle en $C$ :
$\H^i(k,  \check{C} \otimes \overline{k}^*) \cong \H^i(k, C)$ pour tout $i$, $\H^i(k,  \check{C} \otimes \overline{\A}^*) \cong \bigoplus_v \H^i(\widehat{k}_v, C)$ pour $i \geq 1$.
On identifie \'egalement le groupe $\H^0(k, \check{C} \otimes  \overline{\A}^*) \cong \P^0(k, C)$, et on remarque que le morphisme $ \H^3(k, C) \rightarrow \P^3(k, C)$ est un isomorphisme. 
Enfin, on dispose des accouplements de dualit\'e globale suivants, induits par le cup-produit :
\begin{equation}
\label{cp Bor}
\H^i(k, \check{C} \otimes \overline{\C}^*) \times \H^{1-i}(k, \widehat{C}) \rightarrow \Q / \Z
\end{equation}
pour tout $i$, induisant des isomorphismes $\H^1(k, \check{C} \otimes \overline{\C}^*)
\cong \H^0(k, \widehat{C})^D$ et $\H^2(k, \check{C} \otimes
\overline{\C}^*) \cong \H^{-1}(k, \widehat{C})^D$ (voir par exemple \cite{Mil}, exemple I.1.11 et corollaire I.4.7). 
On utilise \'egalement les isomorphismes de dualit\'e suivants, qui s'obtiennent par d\'evissage \`a partir du cas des tores (voir \cite{Mil}, corollaire I.4.7) :
$$\H^0(k,  \check{C} \otimes \overline{\C}^*)^{\wedge} \cong \H^1(k, \widehat{C})^D$$
$$\H^{-1}(k,  \check{C} \otimes \overline{\C}^*)^{\wedge} \cong \H^2(k, \widehat{C})^D$$
Enfin, on a besoin d'identifier les fl\`eches de dualit\'e globale, c'est-\`a-dire qu'il faut comparer les morphismes $\H^i(k, \widehat{C})^D \rightarrow \H^{2-i}(k, C)$ induits respectivement par les th\'eor\`emes \ref{theo SH 1}, \ref{theo SH 2} et \ref{dualite SH 0}, et par les identifications pr\'ec\'edentes dans la suite exacte (\ref{sel Bor}). Pour cela, on va utiliser la section 6 de \cite{HSz} : on consid\`ere les trois suites exactes suivantes de complexes de $\Gamma_k$-modules :
\begin{equation}
\label{sec HSz}
0 \rightarrow C(\overline{k}) \rightarrow C(\overline{\A}) \rightarrow C(\overline{\C}) \rightarrow 0
\end{equation}
$$0 \rightarrow \widehat{C}(\overline{k}) \rightarrow \widehat{C}(\overline{\A}) \rightarrow \widehat{C}(\overline{\C}) \rightarrow 0$$
$$0 \rightarrow \overline{k}^*[1] \rightarrow \overline{\A}^*[1] \rightarrow \overline{\C}^*[1] \rightarrow 0$$
o\`u par d\'efinition $C(\overline{\C})$ est le complexe $[T_1(\overline{\A}) / T_1(\overline{k}) \rightarrow T_2(\overline{\A}) / T_2(\overline{k})]$ (et de m\^eme pour $\widehat{C}(\overline{\C})$). On remarque d'abord que l'on a un isomorphisme naturel de suites exactes courtes entre la suite (\ref{sec Bor}) et la suite (\ref{sec HSz}).
On dispose d'un accouplement naturel
$$C(\overline{\A}) \otimes^{\L} \widehat{C}(\overline{\A}) \rightarrow \overline{\A}^*[1]$$
induisant l'accouplement usuel
$$C(\overline{k}) \otimes^{\L} \widehat{C}(\overline{k}) \rightarrow \overline{k}^*[1]$$
On en d\'eduit alors comme dans la section 6 de \cite{HSz} un accouplement
\begin{equation}
\label{cup-prod}
\textup{Ker} \left( \H^i(k, C(\overline{k})) \rightarrow \H^i(k, C(\overline{\A})) \right) \times \textup{Ker} \left( \H^{2-i}(k, \widehat{C}(\overline{k})) \rightarrow \H^{2-i}(k, \widehat{C}(\overline{\A})) \right) \rightarrow H^2(k, \overline{\C}^*) \cong \Q / \Z
\end{equation}
d\'efini explicitement par le cup-produit $$\H^{i-1}(k, C(\overline{C})) \times \H^{2-i}(k, \widehat{C}(\overline{k})) \rightarrow H^2(k, \overline{\C}^*) \cong \Q / \Z$$
Or on peut identifier les deux groupes apparaissant dans l'accouplement (\ref{cup-prod}) :\\
$\textup{Ker} \left( \H^i(k, C(\overline{k})) \rightarrow \H^i(k, C(\overline{\A})) \right) \cong \textup{\cyr{SH}}^i(k, C)$ et $\textup{Ker} \left( \H^{2-i}(k, \widehat{C}(\overline{k})) \rightarrow \H^{2-i}(k, \widehat{C}(\overline{\A})) \right) \cong \textup{\cyr{SH}}^{2-i}(k, \widehat{C})$. Et la preuve de la proposition 6.1 de \cite{HSz} (adapt\'ee au contexte des complexes de tores) assure que les accouplements (\ref{cup-prod}) co\"incident avec les accouplement des th\'eor\`emes \ref{theo SH 1}, \ref{theo SH 2} et \ref{dualite SH 0}. Cette identification donne en particulier une description explicite des accouplements apparaissant dans ces th\'eor\`emes, en termes de cup-produits en cohomologie galoisienne.
Pour finir la comparaison avec la suite exacte (\ref{sel Bor}), il suffit de constater que le diagramme suivant est commutatif (pour $i = 0, 1, 2$) :
\begin{displaymath}
\xymatrix{
\H^{i-1}(k, \check{C} \otimes \overline{C}^*) \ar[r]^{\phi} \ar[d]^{\partial_{\textup{B}}} & \H^{2-i}(k, \widehat{C})^D \ar[d]^{\partial_{\textup{PT}}} \\
\H^i(k, \check{C} \otimes \overline{k}^*) \ar[r]^{\psi} & \H^i(k, C)
}
\end{displaymath}
 (en rempla\c cant $\H^i(k, C)$ par son compl\'et\'e pour $i = 0$), o\`u le morphisme $\partial_{\textup{B}}$ est le cobord provenant de la suite exacte (\ref{sec Bor}), $\partial_{\textup{PT}}$ provient des th\'eor\`emes de dualit\'e globale \ref{theo SH 1}, \ref{theo SH 2} et \ref{dualite SH 0}, $\phi$ est induit par le cup-produit (\ref{cp Bor}) et $\psi$ est l'identification naturelle $\H^i(k,  \check{C} \otimes \overline{k}^*) \cong \H^i(k, C)$ (compos\'ee avec la compl\'etion si $i = 0$). Et cette commutativit\'e r\'esulte facilement de la comparaison \'etablie pr\'ec\'edemment entre les accouplement des th\'eor\`emes \ref{theo SH 1}, \ref{theo SH 2} et \ref{dualite SH 0} et les accouplements (\ref{cup-prod}).

On peut ainsi r\'esumer les r\'esultats de cette section sous la forme suivante : on a un diagramme commutatif de suites exactes, fonctoriel en $C$, entre la suite exacte longue (\ref{sel Bor}) et la suite exacte de Poitou-Tate du th\'eor\`eme \ref{theo PT} :
\begin{changemargin}{-0.8cm}{1cm}
\begin{displaymath}
\xymatrix{
\H^{-1}(k, \check{C} \otimes \overline{k}^*) \ar@{^{(}->}[rr] \ar[rd]^{\cong} & & \H^{-1}(k, \check{C} \otimes \overline{\A}^*) \ar[rr] \ar[rd]^{\cong} & & \H^{-1}(k,  \check{C} \otimes \overline{\C}^*) \ar'[d][dd] \ar[rd] & \\
& \H^{-1}(k, C) \ar@{^{(}->}[rr] & & \P^{-1}(k, C) \ar[rr] & & \H^2(k, \widehat{C})^D \ar[dd] \\
\H^0(k, \check{C} \otimes \overline{\C}^*) \ar[rd] \ar[dd] & & \H^0(k, \check{C} \otimes \overline{A}^*) \ar[rd] \ar[ll] & & \H^0(k,\check{C} \otimes \overline{k}^*) \ar[rd] \ar[ll] & \\
& \H^1(k, \widehat{C})^D \ar[dd] &&  \P^0(k, C)^{\wedge} \ar[ll] & & \H^0(k, C)^{\wedge} \ar[ll] \\
\H^1(k, \check{C} \otimes \overline{k}^*) \ar'[r][rr] \ar[rd]^{\cong} & & \H^1(k,\check{C} \otimes \overline{\A}^*) \ar[rr] \ar[rd]^{\cong} & & \H^1(k, \check{C} \otimes \overline{C}^*) \ar[rd]^{\cong} \ar'[d][dd] & \\
& \H^1(k,C) \ar[rr] & & \P^1(k, C) \ar[rr] & & \H^0(k, \widehat{C})^D \ar[dd] \\
\H^2(k, \check{C} \otimes \overline{C}^*) \ar[rd]^{\cong} & & \H^2(k,\check{C} \otimes \overline{\A}^*) \ar@{->>}[ll] \ar[rd]^{\cong}& & \H^2(k, \check{C} \otimes \overline{k}^*) \ar[ll] \ar[rd]^{\cong} & \\
& \H^{-1}(k, \widehat{C})^D & & \P^2(k, C) \ar@{->>}[ll] & & \H^2(k, C) \ar[ll]
}
\end{displaymath}
\end{changemargin}
o\`u les fl\`eches obliques qui ne sont pas des isomorphismes sont des fl\`eches de compl\'etion profinie. On peut donc dire en quelque sorte que la suite de Poitou-Tate du th\'eor\`eme \ref{theo PT} est la compl\'etion profinie de la suite d'hypercohomologie (4.3.1) de \cite{BorAMS}. Le fait de consid\'erer la suite "compl\'et\'ee" fait appara\^itre des groupes plus facilement identifiables : par exemple, si $C$ est un complexe de tores associ\'e \`a un $k$-groupe r\'eductif (voir \cite{BorAMS} par exemple), le groupe $\H^0(k, \check{C} \otimes \overline{\C})$ peut sembler myst\'erieux dans ce contexte, alors que sa compl\'etion profinie s'identifie naturellement au dual du groupe de Brauer alg\'ebrique de $G$, ce qui sera tr\`es utile dans les applications (voir \cite{Dem2}).

\subsection{Cas o\`u $\rho$ est surjective}

Dans cette partie, on \'etablit une suite exacte de Poitou-Tate pour un groupe de type multiplicatif.

\begin{theo}
\label{theo PT tm}
Soit $M$ un $k$-groupe de type multiplicatif. On a alors une suite
exacte fonctorielle en $M$ :
\begin{displaymath}
\xymatrix{
0 \ar[r] & H^{0}(k, M)^{\wedge} \ar[r] & P^{0}(k, M)^{\wedge} \ar[r] & H^2(k, \widehat{M})^D \ar[d] & \\
& H^1(k, \widehat{M})^D \ar[d] & P^1(k, M) \ar[l] & H^1(k, M) \ar[l] & \\
& H^2(k,M) \ar[r] & P^2(k, M) \ar[r] & H^0(k, \widehat{M})^D \ar[r] & 0 
}
\end{displaymath}

On dispose \'egalement de la suite exacte duale :
\begin{displaymath}
\xymatrix{
0 \ar[r] & H^{0}(k, \widehat{M})^{\wedge} \ar[r] & P^{0}(k, \widehat{M})^{\wedge} \ar[r] & H^2(k, M)^D \ar[d] & \\
& H^1(k, M)^D \ar[d] & P^1(k, \widehat{M}) \ar[l] & H^1(k, \widehat{M}) \ar[l] & \\
& H^2(k,\widehat{M}) \ar[r] & P^2(k, \widehat{M})_{\textup{tors}} \ar[r] & \left( H^0(k, M)^D \right)_{\textup{tors}} \ar[r] & 0
}
\end{displaymath}
\end{theo}

\begin{dem}
On voit $M$ comme le noyau d'un morphisme surjectif de $k$-tores, $M := \textup{Ker} (\rho : T_1 \rightarrow T_2)$, et on note $C := [T_1 \xrightarrow{\rho} T_2]$ le complexe de tores associ\'e. On a bien un quasi-isomorphisme $C \cong M[1]$.
La preuve est tr\`es similaire \`a celle du th\'eor\`eme \ref{theo PT}, \`a la diff\'erence que l'on utilise ici $\varinjlim_n T_{\Z / n}(C) = 0$ alors que dans la preuve du th\'eor\`eme \ref{theo PT} on avait utilis\'e $\varprojlim_n { _n \textup{Ker } \rho} = 0$.
\begin{itemize}
\item Pour la premi\`ere ligne, on montre par d\'evissage \`a l'aide du triangle exact 
$${ _n T_1} \rightarrow { _n T_2} \rightarrow (C \otimes^{\L} \Z / n) [-1] \rightarrow { _n T_1}[1]$$
l'exactitude de la suite suivante
$$\varprojlim_n \H^{-1}(k, C \otimes^{\L} \Z / n) \rightarrow \varprojlim_n \P^{-1}(k, C \otimes^{\L} \Z / n) \rightarrow \varprojlim_n \H^1(k, \widehat{C} \otimes^{\L} \Z / n)$$
Pour cela, on utilise notamment le fait que $\varprojlim_n \textup{\cyr{SH}}^1(k, { _n T_i}) = \textup{Ker} \left( T_i(k)_{\wedge} \rightarrow P^0(k, T_i)_{\wedge} \right) = 0$, ainsi que la finitude du groupe $\textup{\cyr{SH}}^0(C)$.
\item Pour la deuxi\`eme ligne, on d\'eduit l'exactitude de la suite $\varinjlim_n \H^{-1}(k, C \otimes^{\L} \Z / n) \rightarrow \varinjlim_n \P^{-1}(k, C \otimes^{\L} \Z / n) \rightarrow \varinjlim_n \H^1(k, \widehat{C} \otimes^{\L} \Z / n)$
du fait que $\varinjlim_n T_{\Z / n}(C) = 0$, et on conclut par finitude de $\textup{\cyr{SH}}^2(\widehat{C})$.
\item Pour la troisi\`eme ligne, la suite $\varinjlim_n \H^{0}(k, C \otimes^{\L} \Z / n) \rightarrow \varinjlim_n \P^{0}(k, C \otimes^{\L} \Z / n) \rightarrow \varinjlim_n \H^0(k, \widehat{C} \otimes^{\L} \Z / n)$ est exacte, et on conclut par finitude de $\textup{\cyr{SH}}^1(\widehat{C})$ et par le fait que $\H^1(k, C)$ et $\P^1(k, C)$ sont de torsion.
  \item Pour la suite duale, on utilise des arguments similaires.
\end{itemize}
\end{dem}

\addcontentsline{toc}{section}{R\'ef\'erences}
\bibliography{tores}
\bibliographystyle{amsalpha}

Cyril Demarche \\
Laboratoire de Math\'ematiques, B\^atiment 425, Universit\'e de Paris-Sud, F-91405 Orsay, France \\
e-mail : cyril.demarche@math.u-psud.fr

\end{document}